\documentclass{article}

\usepackage{arxiv}
\usepackage[english]{babel}
\usepackage{tabularx}
\usepackage{multirow}
\usepackage{amsmath,amsfonts}
\usepackage{graphicx}
\usepackage{xcolor}
\usepackage{subfigure}
\usepackage[ruled,vlined]{algorithm2e}
\usepackage[utf8]{inputenc} 
\usepackage[T1]{fontenc}    
\usepackage{hyperref}       
\usepackage{url}            
\usepackage{booktabs}       
\usepackage{nicefrac}       
\usepackage{microtype}      
\usepackage{graphicx}
\usepackage[square,sort,comma,numbers]{natbib}
\usepackage{doi}
\usepackage{setspace}
\usepackage{afterpage}
\usepackage{float}

\def\bp{\boldsymbol{p}}

\def\bw{\boldsymbol{w}}

\newcommand{\R}{\mathbb{R}}
\newcommand{\C}{\mathbb{C}}
\newcommand{\qed}{\hfill \ensuremath{\Box}}

\newtheorem{definition}{Definition}[section]
\newtheorem{theorem}{Theorem}[section]
\newtheorem{lemma}{Lemma}[section]

\renewcommand{\shorttitle}{Numerical Solution of Stiff ODEs with Physics-Informed RPNNs}

\begin{document}

\title{Numerical Solution of Stiff ODEs with Physics-Informed Random Projection Neural Networks
}

\author{Evangelos Galaris\\
	Dipartimento di Matematica e Applicazioni ``Renato Caccioppoli''\\
Universit\`a degli Studi di Napoli Federico II\\
	Napoli, Italy\\
	\texttt{evangelos.galaris@unina.it} \\
	\And
	Gianluca Fabiani \\
	Scuola Superiore Meridionale \\
	Universit\`a degli Studi di Napoli Federico II\\
	Napoli, Italy\\
	\texttt{gianluca.fabiani@unina.it} \\
	\And
	Francesco Calabr\`o \\
	Dipartimento di Matematica e Applicazioni ``Renato Caccioppoli''\\
        Universit\`a degli Studi di Napoli Federico II\\
	Napoli, Italy\\
	\texttt{calabro@unina.it} \\
	\And
	Daniela di Serafino\\
	Dipartimento di Matematica e Applicazioni ``Renato Caccioppoli''\\
        Universit\`a degli Studi di Napoli Federico II\\
	Napoli, Italy\\
	\texttt{daniela.diserafino@unina.it} \\
	\And
	Constantinos Siettos\thanks{Corresponding author} \\
	Dipartimento di Matematica e Applicazioni ``Renato Caccioppoli''\\
        Universit\`a degli Studi di Napoli Federico II\\
	Napoli, Italy\\
	\texttt{constantinos.siettos@unina.it} \\
}
	

\date{November 16, 2021}

\maketitle

\begin{abstract}
We propose a numerical method based on physics-informed Random Projection Neural Networks for the solution of Initial Value Problems (IVPs) of Ordinary Differential Equations (ODEs) with a focus on stiff problems. We address an Extreme Learning Machine with a single hidden layer with radial basis functions having as widths uniformly distributed random variables, while the values of the weights between the input and the hidden layer are set equal to one. The numerical solution of the IVPs is obtained by constructing a system of nonlinear algebraic equations, which is solved with respect to the output weights by the Gauss-Newton method, using a simple adaptive scheme for adjusting the time interval of integration. To assess its performance, we apply the proposed method for the solution of four benchmark stiff IVPs, namely the Prothero-Robinson, van der Pol, ROBER and HIRES problems.
Our method is compared with an adaptive Runge-Kutta method based on the Dormand-Prince pair, and a variable-step variable-order multistep solver based on numerical differentiation formulas, as implemented in the \texttt{ode45} and \texttt{ode15s} MATLAB functions, respectively. We show that the proposed scheme yields good approximation accuracy, thus outperforming \texttt{ode45} and \texttt{ode15s}, especially in the cases where steep gradients arise. Furthermore, the computational times of our approach are comparable with those of the two MATLAB solvers for practical purposes.
\end{abstract}

\section{Introduction\label{sec:intro}}
The idea of using Artificial Neural Networks (ANNs) for the numerical solution of Ordinary Differential Equations (ODEs) dates back to the '90s. One of the first works in the field was that of Lee and Kang~\cite{lee1990neural}, who addressed a modified Hopfield Neural Network to solve a first-order nonlinear ODE. This method is based on the discretization of the differential operator with finite differences and the minimization of a related energy function. Following up this work, Meade and Fernadez \cite{meade1994numerical} used a non-iterative scheme based on Feedforward Neural Networks (FNN) for the solution of linear ODEs, where the estimation of the weights of the FNN is based on the Galerkin weighted-residuals method. In 1998, Lagaris et al. \cite{lagaris1998artificial} introduced a numerical method based on FNNs for the solution of nonlinear ODEs and PDEs. The method constructs appropriate trial functions with the aid of a single-hidden layer FNN that is trained to minimize the error between the FNN prediction and the right-hand side of the differential equations. The proposed approach is demonstrated through both initial and boundary-value problems and a comparison with Galerkin finite elements is also provided. Based on the work of Lagaris et al. \cite{lagaris1998artificial}, Filici \cite{filici2008neural} proposed a single-layer multiple linear output perceptron providing a proof for the error bound. Network training was performed with the LMDER optimization solver from MINPACK \cite{more1980user,more1984minpack}. The performance of the approach was tested through simple ODE problems, including the van der Pol equation, but without stiffness. Tsoulos et al. \cite{tsoulos2009solving} employed an FNN trained by grammatical evolution and a local optimization procedure to solve second-order ODEs. More recently, Dufera \cite{dufera2021deep} addressed a deep learning network with back propagation. The performance of the approach was tested with the non-stiff ODEs presented in Lagaris et al. \cite{lagaris1998artificial}. The authors report that their scheme outperforms the non-adaptive 4-th order Runge-Kutta method in terms of numerical approximation accuracy when considering short time intervals and a small number of collocation points. A review and presentation of various ANN schemes for the solution of ODEs can be found in Yadav et al. \cite{yadav2015introduction}. In all the above procedures, a computationally demanding optimization algorithm is required for the evaluation of the parameters of the network and the ODEs problems are non-stiff.

Yang et al. \cite{yang2018novel} aimed to solve ODEs using the so-called Extreme Learning Machine (ELM) concept \cite{Huang,huang2015extreme}: the weights between the input and the hidden layer as well as the biases of the hidden nodes, were chosen randomly, and the remaining unknown weights between the hidden and the output layer were computed in one step solving a least squares problem with regularization. The authors addressed a single-layer Legendre neural network to solve non-stiff ODEs up to second order, including the Emden-Fowler equation. The performance of the scheme was compared against other machine learning schemes, including a cosine basis function neural network trained by the gradient descent algorithm, the explicit Euler scheme, the Suen 3rd-order and the classical 4th-order Runge-Kutta methods. The authors reported high numerical accuracy and computing times comparable or smaller than the other methods.

The above studies deal with non-stiff ODE problems. One of the first papers that dealt with the solution of stiff ODEs and Differential-Algebraic Equations using ANNs was that of Gerstberger and Rentrop \cite{gerstberger1997feedforward}, where a FNN architecture was proposed to implement the implicit Euler scheme. The performance of that architecture was demonstrated using scalar ODEs and the van der Pol equations considering however only mild stiffness, setting the parameter that controls the stiffness to values up to five.

More recently, theoretical and technological advances have renewed the interest for developing and applying physics-informed machine learning techniques for learning and solving differential equations \cite{raissi2019physics,jagtap2020conservative,meng2020ppinn,schiassi2021extreme,dong2021local,wang2021understanding,karniadakis2021physics}. Physics-Informed Neural Networks (PINN) are trained to solve supervised learning tasks while respecting the given laws of physics described by nonlinear differential equations. This is also part of the so-called scientific machine learning, which is emerging as a potential alternative to classical scientific computing. Due to the fact that (large-scale systems of) stiff ODEs arise in modelling an extremely wide range of problems, from biology and neuroscience to engineering processes, material science and chemical kinetics, and from financial systems to social science and epidemiology, there is a re-emerging interest in developing new methods for their efficient numerical solution \cite{hadjinicolaou1998asymptotic,valorani2001explicit,goussis2006efficient,ji2021stiff,kim2021stiff,de2021physics}. Within this framework, Mall and Chakraverty \cite{mall2016hermite} proposed a single-layer Hermite polynomial-based neural network trained with back-propagation to approximate the solutions of the van der Pol-Duffing and Duffing oscillators with low to medium values of the stiffness parameter. Following up this work, Chakraverty and Mall \cite{chakraverty2020single} proposed a single-layer Chebyshev neural network with regression-based weights to solve first- and second-order nonlinear ODEs and in particular the nonlinear Lane-Emden equation. The training was achieved using back propagation. Budkina et al.~\cite{budkina2017neural} and Famelis and Kaloutsa~\cite{famelis2021parameterized} proposed a single-layer ANN to solve highly stiff ODE problems. The networks were trained for different values of the parameter that controls stiffness over relatively short time intervals using the Levenberg-Marquardt algorithm.

The above studies on the solution of stiff ODEs report good numerical approximation accuracy of the proposed schemes, in general for relatively short time intervals and small sizes of the grid. However, no indication is given about the computational time required for the training versus the time required by widely used solvers for low-to-medium and medium-to-high stiff problems such as the \texttt{ode45} and \texttt{ode15s} functions available from the MATLAB ODE suite \cite{shampine1997matlab}. Recently, Wang et al.~\cite{wang2021understanding} revealed the difficulties of PINNs to solve gradient flow problems when the dynamics are stiff. These difficulties or even the failure of PINNs to deal with such problems are due to the arising unbalanced back-propagated gradients during model training. To deal with this problem, the authors addressed an annealing algorithm that utilizes gradient statistics during model training to balance between different terms in composite loss functions. To deal with the presence of stiff dynamics in problems of chemical kinetics dynamics, Ji et al.~\cite{ji2021stiff} proposed a quasi-steady-state approximation to reduce the stiffness of the ODE systems, and then applied PINNs for their solution. For the demonstration of the proposed approach, they used two classical stiff chemical kinetics problems, namely the ROBER~\cite{Robertson1966} and POLLU \cite{verwer1994gauss} systems.

\subsection{Our Contribution}

We propose a numerical scheme based on physics-informed Random Projection Neural Networks (RPNNs), and in particular on ELMs, for the solution of Initial Value Problems (IVPs) of ODEs with a focus on stiff problems. RPNNs are a family of networks including randomized and Random Vector Functional Link Networks (RVFLNs) \cite{schmidt1992feed,pao1994learning,husmeier1999random}, Echo-State Neural Networks and Reservoir Computing \cite{jaeger2002adaptive,jaeger2004harnessing,ozturk2007analysis,sakemi2020model,paquot2012optoelectronic}, and Extreme Learning Machines \cite{Huang,huang2015extreme,huang2014insight}. The keystone idea behind all these approaches is to use a fixed-weight configuration between the input and the hidden layer, fixed biases for the nodes of the hidden layer, and a linear output layer. Hence, the output is projected linearly onto the functional subspace spanned by the nonlinear basis functions of the hidden layer, and the only unknowns that have to be determined are the weights between the hidden and the output layer. Their estimation is done in one step by solving a (nonlinear) least squares problem.

Here, for the solution of stiff problems of ODEs, we propose a single-layer FNN with Radial Basis Functions (RBFs) with parameters that are properly uniformly-distributed random variables. Thus, the proposed neural network scheme constitutes a Lipschitz embedding constructed through the random projection. The feasibility of the scheme is guaranteed by the celebrated Johnson and Lindenstrauss Lemma \cite{johnson1984extensions} and the universal approximation theorems proved for random projection networks, in particular for ELMs \cite{huang2006universal}.
By combining this choice of the underlying functional space, we can also explicitly compute the derivatives, and hence the Jacobian matrix for the Gauss-Newton method used for solving the resulting system of nonlinear algebraic equations, which is equivalent to an under-determined least squares problem. Thus, we obtain an efficient way to calculate by collocation a neural network function that approximates the exact solution in any point of the domain. To deal with the presence of steep gradients that may arise, we propose a simple adaptive strategy for adjusting the time interval of integration. In previous works, we have shown that ELMs may efficiently deal with differential systems with solutions containing steep gradients \cite{calabro2021extreme,fabiani2021numerical}.

To demonstrate the efficiency of the proposed method, we have chosen four benchmark stiff ODE problems of different dimensions, namely the Prothero-Robinson ODE \cite{prothero1974stability}, for which an analytical solution exists, the well-known van der Pol equations \cite{van1926lxxxviii}, the ROBER problem, a stiff system of three nonlinear ODEs describing the kinetics of an autocatalytic reaction \cite{Robertson1966}, and the HIRES problem, a system of eight nonlinear stiff ODEs describing morphogenesis in plant physiology \cite{schafer1975new} (see also \cite{mazzia2012test}). The performance of the proposed scheme is assessed in terms of both approximation accuracy and computational times. We show that the proposed scheme outperforms \texttt{ode15s} in terms of numerical approximation accuracy, especially in the cases where steep gradients arise, while \texttt{ode45} in some cases completely fails or needs many points to satisfy a specific tolerance.  Moreover, our approach provides the solution directly as a function that can be evaluated at every point of the domain, in contrast to the classical numerical methods, where extra interpolation steps are usually performed for this purpose (for example, for this task the approximate solutions can be obtained using the MATLAB function \texttt{deval}, which uses an interpolation technique). It is important to mention that once the random parameters are chosen in a proper interval, the underlying space is capable to catch the steep gradient behaviours of the solution, providing good numerical accuracy especially for stiff ODE problems, being more robust and thus outperforming \texttt{ode45} and \texttt{ode15s} in terms of numerical accuracy, while resulting in comparable computational times. Moreover, when the solution has to be calculated in high-density grids of points, the proposed method results in smaller computational times when compared with \texttt{ode45} and \texttt{ode15s}.

Finally, we remark that in this work, we present a new machine learning numerical algorithm for the solution of stiff IVPs of ODEs with solutions that may also contain steep gradients. We aim neither to perform an exhaustive comparison of the proposed scheme with all the ``classical'' numerical methods that have been proposed for the solution of stiff problems of ODEs nor to extensively assess its performance on the many benchmark problems that have been suggested over the years (see for example the test data set in \cite{mazzia2012test}). Such an extensive analysis can be the subject of future work(s).

The structure of the paper is as follows. In Section 2, we give briefly a coarse definition of stiff IVPs of ODEs. In Section~3, we provide some preliminaries on the solution of ODEs with FFNs based on the ``classical'' optimization approach, and the basic concept behind the use of RPNNs; we also discuss briefly the celebrated Johnson and Lindenstrauss Lemma \cite{johnson1984extensions}. In Section 4, we describe our approach for the solution of IVPs of ODEs with the use of RPNNs, providing also a pseudo-code of the corresponding algorithm, and discuss its approximation properties within the framework of the universal approximation theorem of ELMs \cite{huang2006universal}. In section 5, we present the numerical results obtained by applying the proposed approach to the above-mentioned stiff ODE problems along with a comparison with \texttt{ode45} and \texttt{ode15s}. Conclusions are given in Section~\ref{sec:discussion}.

\section{Stiff IVPs of ODEs\label{sec:stiff-odes}}

We aim to solve IVPs of ODEs of the following form:
\begin{equation}\label{eq:ODE_PB}
   \left\{ \begin{array}{lll}
       \displaystyle \dfrac{d y_i}{dx} & = & \displaystyle f_i(x, y_1,y_2,...,y_m), \\[7pt]
       \displaystyle y_i(x_0) & = & \displaystyle \alpha_i, 
   \end{array} \right.
   \quad i = 1,2,\ldots,m,
\end{equation}
where the functions $f_i$ and the initial values $\alpha_i$ are known, and the functions $y_i(x)$ are the unknowns. In order to simplify the notation, we group the functions $y_i(x)$ in a vector function $\boldsymbol{y}(x): \R \to \R^m$, and the functions $f_i(x, y_1,y_2,...,y_m)$ in $\boldsymbol{f}(x,\boldsymbol{y})$.

The unknown solution may exhibit a complex behaviour, including steep gradients, which results in difficulties for the design and implementation of numerical methods. These difficult problems are usually the ones where the so-called stiffness arises. Until now, there is no complete and precise definition of what is called ``stiffness''. Following Lambert \cite{JDLambert1992}, one has to consider stiffness as a phenomenon exhibited by the system rather than a property of it. Generally speaking, an ODE problem is called stiff if there exists a solution that varies very slowly, but there are nearby solutions that vary rapidly, so that (explicit) numerical algorithms need extremely many small steps to obtain accurate and reliable results.

A widely used definition for stiffness is the following, given by Lambert (see also~\cite{brugnano2011fifty}):
\begin{quote}
\textit{If a numerical method with a finite region of absolute stability, applied to a system with any initial conditions, is forced to use in a certain interval of integration a step length which is excessively small in relation to the smoothness of the exact solution in that interval, then the system is said to be stiff in that interval.}
\end{quote}
In the case of constant coefficients, one can introduce the stiffness ratio, thus defining stiffness as follows. Consider the constant coefficient linear inhomogeneous system in the unknown $\boldsymbol{y}: \R \to \R^m$:
\begin{equation}\label{st}
    \dfrac{d \boldsymbol{y}}{d{x}} = A \boldsymbol{y}+\overline{\boldsymbol{f}}({x}), 
\end{equation}
where $\overline{\boldsymbol{f}} : \R \to \R^m$ is a nonlinear term and $A \in \R^{m \times m}$ is a constant and diagonalizable matrix with eigenvalues $\lambda_l \in \C$, $l=1,2,...,m$, (assumed distinct) and corresponding eigenvectors $\boldsymbol{c}_l \in \C^m$. 
Suppose also that:
\begin{equation}\label{def_hyp}
    Re(\lambda_l)<0, \quad l=1,2,...,m.
\end{equation}
\begin{definition}[Stiffness ratio]
With the above notations, let $\overline {\lambda },\underline {\lambda }\in \{\lambda _l\,,\ l=1,2,\ldots ,m\}$ be such that
\begin{equation*}
  |Re(\overline {\lambda })|\geq |Re(\lambda _l)|\geq |Re(\underline {\lambda })|,\quad l=1,2,\ldots ,m.
\end{equation*}
The stiffness ratio of system \eqref{st} is defined as
\begin{equation*}
  \frac{|Re(\overline{\lambda})|}{|Re(\underline{\lambda})|}.
\end{equation*} 
\end{definition}
The above definition is motivated by the explicit general solution of \eqref{st}, that is
\begin{equation}
    \boldsymbol{y}(x)=\sum_{l=1}^m\kappa_l exp(\lambda_l x) \boldsymbol{c}_l+\boldsymbol{g}(x),
\end{equation}
where $\boldsymbol{g}(x)$ is a particular integral taking into account the forcing term. Assuming that \eqref{def_hyp} holds true, one has that each of the terms $\exp(\lambda _{l} x){\boldsymbol{c}}_{l}$ vanishes as $x \rightarrow \infty$, and hence the solution $\boldsymbol{y}(x)$ approaches $\boldsymbol{g}(x)$ asymptotically as $x\rightarrow \infty$. In this case, the term $\exp(\lambda _{l} x)\boldsymbol{c} _l$ decays monotonically if $\lambda_l$ is real, and sinusoidally if $\lambda_l$ is complex. Interpreting $x$ to be time (as it is often in physical problems), $\sum _{l=1}^{n}\kappa _l\exp(\lambda _l x)\boldsymbol{c} _l$ is called the transient solution and $\boldsymbol{g}(x)$ the steady-state solution. If $|Re(\lambda _{l})|$ is large, then the corresponding term $\kappa _l\exp(\lambda _l x)\boldsymbol{c} _l$ decays quickly as $x$ increases and thus it is called a fast transient; if $|Re (\lambda _l)|$ is small, then $ \kappa _l\exp(\lambda_l x)\boldsymbol{c} _t$ decays slowly and is called a slow transient. 
Letting $\overline l$ and $\underline l$ be the indices identifying $\overline\lambda$ and $\underline\lambda$, respectively, we have that $\kappa_{\overline l} \exp(\overline {\lambda }x) \boldsymbol{c}_{\overline l}$ is the fastest transient and $\kappa_{\underline l} \exp(\underline{\lambda }x) \boldsymbol{c}_{\underline l}$ the slowest. The ratio between these two terms gives the stiff behaviour of the linear ODE system.

Many examples of stiff problems exhibit also other features, but for each feature there maybe other stiff problems not exhibiting that particular feature. However, we note that Lambert refers to these features as ``statements'' rather than definitions. A few of them are the following: (1) a linear constant coefficient system is stiff if all of its eigenvalues have negative real part and the stiffness ratio is large; (2) stiffness occurs when stability requirements, rather than those of accuracy, constrain the step length; and (3) stiffness occurs when some components of the solution decay much more rapidly than others.
One could also say that stiffness is a numerical efficiency issue since non-stiff methods can in theory solve stiff problems if they use many points (thus requiring very large computing times). In practice, many extremely stiff problems bias ``classical'' numerical methods to tiny steps, thus making the computational cost huge and the methods impractical.

\section{Preliminaries\label{sec:preliminaries}}

\subsection{Feedforward Neural Networks\label{sec:fnn}}

An FNN is a biologically inspired regression and/or classification algorithm. It consists of processing units, called neurons or nodes, organised in layers. There are three types of layers: the input layer, the output layer and possibly one or more hidden layers. The units of the output and every hidden layer are connected to all the units of the previous layer. These connections may have different strengths, called weights, and encode the ``knowledge'' of the network. The data enter the input layer and pass through each layer during the forward phase. Each unit accepts a signal consisting of a (usually linear) combination of the outputs of the previous layer units and, using an activation function, creates an output that is transmitted to the next layer up to the final output layer. Different activation functions and/or numbers of units can be used for each layer. A simple example of FNN is the Single Layer FNN (SLFNN), consisting of $d$ input units, a single layer with $h$ hidden units with biases,
and $k$ output units (without biases). Given an input $\boldsymbol{x}\in\R^d$, the output $\boldsymbol{N}(\boldsymbol{x})\in \R^k$ of this SLFNN reads:
\begin{equation} \label{eqn:net}
   \boldsymbol{N}(\boldsymbol{x}) = W^{o}\boldsymbol{\Phi}(W\boldsymbol{x}+\boldsymbol{b}, \boldsymbol{q}),
\end{equation}
where $W^{o} = [w^{o}_{jl}] \in \R^{k \times h}$ is the matrix containing the weights from the hidden layer to the output layer, $\boldsymbol{\Phi}: \R^h \times \R^s \to \R^h$ is a vector function with components the $h$ activation functions $\Phi_j$, $W = [w_{ij}] \in \R^{h\times d}$ is the matrix containing the weights from the input to the hidden layer, $\boldsymbol{b} \in \R^h$ is the vector of the biases of the hidden nodes, and $\boldsymbol{q} \in \R^s$ is a vector containing the hyperparameters of the neural network, such as the parameters of the activation functions (e.g., for radial basis functions, the biases of the hidden neurons and the variances of the Gaussian functions), the learning rate, and the batch size.

Many results are available regarding the approximation properties of FFNs. The most important one from the numerical point of view is the Universal Approximation Theorem, for which we refer to the original papers \cite{Cybenko1989,Hornik1989,hornik1990universal,park1991universal} and the recent review \cite{kratsios2019universal}. Next, we present some results concerning our proposed machine learning scheme. What can be summarized here is that an FNN is capable of approximating uniformly any (piecewise-)continuous (multivariate) function, to any desired accuracy. This implies that any failure of a function mapping by a (multilayer) network arises from an inadequate choice of weights and biases and/or an insufficient number of hidden nodes. Moreover, it has been shown that in the univariate case only one hidden layer is needed.

Given a training set of $n$ input-output pairs, the estimation of the appropriate weights and biases is usually attempted by using an optimization procedure aiming at minimizing a cost function. In particular, the ``classical way'' (see e.g. \cite{lagaris1998artificial}) to solve differential equations in a $d$-dimensional domain with FNNs involves the solution of a minimization problem as follows. One first defines a trial solution in the form of
\begin{equation} \label{eqn:trial_vec}
    \boldsymbol{\Psi}(\boldsymbol{x}, P,\boldsymbol{q}) = \boldsymbol{\Omega}(\boldsymbol{x}, \boldsymbol{N}(\boldsymbol{x}, P,\boldsymbol{q})),
\end{equation}
where $\boldsymbol{\Psi}(\boldsymbol{x}, P,\boldsymbol{q})$ has $m$ components, each associated with a component $y_i$ of the solution $\boldsymbol{y}$, $\boldsymbol{\Omega}: \R^{d} \times \R^{k} \to \R^m$ is sufficiently smooth, $\boldsymbol{N}(\boldsymbol{x}) = \boldsymbol{N}(\boldsymbol{x}, P,\boldsymbol{q})$ has $k$ components $N_i(\boldsymbol{x},\bp_i,\boldsymbol{q})$, and $P$ is a matrix containing the network parameters $W^o$, $W$ and the vector of biases $\boldsymbol{b}$ associated with the $i$-th component of $\boldsymbol{\Psi}(\boldsymbol{x}, P, \boldsymbol{q})$. Then, we consider a discretization of the equation that has to be solved, so as to settle the conditions that the optimization has to manage.

In the case of ODEs, according to the above notation $d=1$, i.e. the input domain is one-dimensional (representing for example time). By discretizing the domain of the ODE problem~\eqref{eq:ODE_PB} with $n$ points $x_l \in \R$, so as to collocate the equations, one can determine the values of the network parameters in $P$ by solving the following optimization problem:
\begin{equation}  \label{eqn:cost}
    \min_{P} E(P) := \sum_{j=1}^n \left\| \frac{d \boldsymbol{\Psi}}{dx} (x_j, P,\boldsymbol{q}) - \boldsymbol{f} (x_j, \boldsymbol{\Psi}(x_j, P,\boldsymbol{q})) \right\|^2.
\end{equation}
More details on the construction of trial solutions and the minimization problem concerning our approach are given in Section~\ref{sec:our_method}.
In order to deal with problem~\eqref{eqn:cost}, one usually needs quantities such as the derivatives of $N_i(x, \bp_i, \boldsymbol{q})$ with respect to the input $x$ and the weights and biases represented by $\bp_i$. These can be obtained in several ways, e.g. by computing analytical derivatives, by using finite difference or other numerical approximations, by symbolic differentiation or by automatic differentiation (see, e.g., \cite{lu2021deepxde} and the references therein).

Although a specialized form of automatic differentiation, known as back-propagation, is widely used in ANNs, in this work we can easily compute analytical first-order derivatives (see also \cite{lagaris1998artificial}), and then apply to problem~\eqref{eqn:cost} a variety of numerical methods that use those derivatives \cite{kelley1999iterative}. For a Single-Input-Single-Output (SISO) neural network (i.e. $d=1$ and $k=1$), with $n$ points and one hidden layer with $h$ elements in \eqref{eqn:net},
%
%
it is straightforward to show that the first-order derivative with respect to $x_{l}, l=1,2,\dots n$, is given by
\begin{equation}
    \frac{\partial N}{\partial x_l} =  \sum_{j=1}^{h} w^{o}_j \frac{\partial \Phi_j}{\partial x_l},
    \label{eqn:gradnet}
\end{equation}
while the derivatives of the network with respect to the training parameters, which are in general the input and output weights $w_j$ and $w^o_j$, $j=1,2,\dots h$, are given by
\begin{equation}
    \frac{\partial N}{\partial w_{j}} = \sum_{i=1}^{h} w^{o}_i \frac{\partial \Phi_i}{w_{j}}, \quad
    \frac{\partial N}{\partial w^o_{j}} =  \Phi_j.
    \label{dngdw}
\end{equation}
(Note that we use the lightface notation $N$ to indicate the SISO version of the neural network $\boldsymbol{N}$).
Once we have computed the above derivatives of the neural network $N$, we can compute the first-order derivatives of the cost function in \eqref{eqn:cost} with respect to the parameters of $N$, and thus we are in principle able to train the network by any minimization technique using those derivatives. However, the training phase may require a very large computing time. Here, we address a scheme that optimizes the parameters via the solution of a suitable nonlinear least squares problem.
Our proposed FNN is based on radial basis transfer functions, introduced in the next section. We remark that another natural choice for the transfer function, widely used in FNN, is the sigmoidal function, which is a bounded monotone S-shaped function with a bell-shaped first derivative (see, e.g. \cite{pinkus}).

\subsection{Radial Basis Function Neural Networks\label{sec:rbfnn}}

Radial Basis Function Neural Networks (RBFNNs) are special cases of ANNs where each hidden unit is associated with an RBF as transfer function. RBFs were addressed in 1985 for multivariate function interpolation \cite{Powell1987,buhmann_2000}. As transfer functions in ANN, RBFs were first addressed by Broomhead and Lowe \cite{Broomhead1988}. An RBF is a function $\phi_{\boldsymbol{c}}: \R^n \to \R$ whose value depends only on the distance between its argument $\boldsymbol{x}$ and some fixed point $\boldsymbol{c}$ called center:
\begin{equation}\label{eq:RBF}
    \phi_{\boldsymbol{c}}(\boldsymbol{x}) = \phi(||\boldsymbol{x}-\boldsymbol{c}||) ,
\end{equation}
where $||\cdot||$ denotes a vector norm. A widely used RBF is the Gaussian RBF:
\begin{equation}\label{eq:Gaussian_RBF}
    G(\boldsymbol{x},\boldsymbol{c},\sigma) = \exp \left(-\dfrac{||\boldsymbol{x}-\boldsymbol{c}||^2}{\sigma^2}\right),
\end{equation}  
employed with the Euclidean norm. This RBF is characterized by the width (otherwise called scale parameter) $\sigma$.

In RBFNNs, each hidden unit is more ``sensitive'' to data points close to its center. Furthermore, for Gaussian RBFNNs, this sensitivity can be tuned by adjusting the width parameter $\sigma$: a larger $\sigma$ implies less sensitivity. Thus, for most of the problems, the choice of the centers is crucial. In Section 4, we  suggest a choice of these parameters in the context of solving IVPs of ODEs.
For the specific case of RBFNNs with the same $\sigma$ in all hidden nodes, Park and Sandberg proved that those networks are capable of universal approximation \cite{park1991universal}. Furthermore, Liao et al. extended this result to
very mild hypotheses on the activation function \cite{Liao2003}, thus stating the following theorem.
\begin{theorem}\label{th:shift}
Let $\phi_{\boldsymbol{c}}$ be as in equation \eqref{eq:RBF}. If $\phi_{\boldsymbol{c}}$ is continuous almost everywhere, locally essentially bounded, and not a polynomial, then for any compact set $K \subset \R^n,\; \Sigma = \mbox{span}\{ \phi_{\boldsymbol{c}}(w\boldsymbol{x}+\boldsymbol{b}) : w\in\R,\ \boldsymbol{b}\in\R^n \}$ is dense in $C(K)$ with respect to the uniform norm, i.e., given any $g \in C(K)$ and $\epsilon>0$, there exists $\gamma \in \Sigma$ such that
\begin{equation*}
    ||\gamma - g ||_{L_{\infty}(K)}\leq\epsilon. 
\end{equation*} 
\end{theorem}

RBFNNs have gained popularity because of a number of advantages compared with other types of ANNs, including their simpler structure and faster learning algorithms, and have been applied to problems in many different areas, e.g. image processing \cite{Er2002, Montazer2015, Basha2021}), speech recognition \cite{Venkateswarlu2011, Renals1989, Balaji2021}, time series analysis \cite{Zemouri2003} and adaptive equalization \cite{Cha1995, Chen1993}.

\subsection{Random Projection RBFNNs\label{sec:rand_proj_rbfnn}}

The training of an ANN requires the minimization of a cost function. But, even for the simple case of SLFNNs, this task may become challenging. Many optimization algorithms used for training ANNs apply stochastic gradient-based approaches which back-propagate the error and adjust the weights through specific directions~\cite{bottou:2018}. More recently, second-order stochastic optimization methods have been widely investigated to get better performances than first-order methods, especially when ill-conditioned problems must be solved (see, e.g., \cite{diserafino2021lsos} and the references therein). Nevertheless, there are still difficulties in using these approaches, such as the setting of the so-called hyperparameters, the significant increase of computing time when the number of data or the number of nodes in the hidden layer grows, and the high non-convexity stemming from the use of nonlinear activation functions, which can lead the algorithms to local minima.

A way to deal with the ``curse of dimensionality'' in training ANNs is to apply the concept of random projection. The idea behind random projections is to construct a Lipschitz mapping that projects the data into a random subspace. The feasibility of this approach can been justified by the celebrated Johnson and Lindenstrauss (JL) Theorem: \cite{johnson1984extensions}:
\begin{theorem}[Johnson and Lindenstrauss]
Let $\mathcal{X}$ 
be a set of $n$ points in $\R^d$. Then, $\forall \, \epsilon \in (0,1)$ and $k \in \mathbb{N}$ such that $k \geq O (\frac{\ln{n}}{\epsilon^2})$, there exists a map $\boldsymbol{F}: \mathbb{R}^d \rightarrow \mathbb{R}^k$ such that
\begin{equation} \label{eqn:jl}
    (1-\epsilon)\lVert \boldsymbol{u}-\boldsymbol{v} \rVert ^2  \leq \lVert \boldsymbol{F}(\boldsymbol{u})-\boldsymbol{F}(\boldsymbol{v}) \rVert ^2 \leq (1+\epsilon) \lVert \boldsymbol{u}-\boldsymbol{v} \rVert ^2 \quad \forall \, \boldsymbol{u},\boldsymbol{v} \in \mathcal{X}.
\end{equation}
\end{theorem} 
Note that while the above theorem is deterministic, its proof relies on probabilistic techniques combined with Kirszbraun's theorem to yield a so-called extension mapping~\cite{johnson1984extensions}.  In particular, it can be shown that one of the many such embedding maps is simply a linear projection matrix
with suitable random entries. Then, the JL Theorem may be proved using the following lemma.
\begin{lemma}
Let $\mathcal{X}$ be a set of $n$ points in $\R^d$, $\epsilon \in (0,1)$ and $\boldsymbol{F}(\boldsymbol{u})$ the random projection defined by
\begin{equation*}
    \boldsymbol{F} (\boldsymbol{u}) = \frac{1}{\sqrt{k}} R \, \boldsymbol{u}, \quad \boldsymbol{u} \in \R^d,
\end{equation*}
where $R = [r_{ij}] \in \R^{k \times d}$ has components that are i.i.d.~random variables sampled from a normal distribution.
Then, $\forall \, \boldsymbol{u} \in \mathcal{X}$
\begin{equation*}
(1-\epsilon) \lVert \boldsymbol{u} \rVert ^2  \leq \lVert \boldsymbol{F} (\boldsymbol{u}) \rVert ^2 \leq (1+\epsilon) \lVert \boldsymbol{u}\rVert ^2 
\end{equation*}
is true with probability $p \geq 1-2\exp \left(-(\epsilon^2-\epsilon^3)\frac{k}{4} \right)$.
\end{lemma}
Similar proofs have been given for distributions different from the normal one (see, e.g., \cite{achlioptas2003database,dasgupta2003elementary,vempala2005random,wang2012geometric}).
In general, the proof of the JL Theorem is based on the fact that inequality~\eqref{eqn:jl} is true with probability 1 if $k$ is large enough. Thus, the theorem states that there exists a projection (referred to as encoder) of $\mathcal{X}$ into a random subspace of dimension $k \geq O (\frac{\ln{n}}{\epsilon^2})$, where the distance between any pair of points in the embedded space $F(\mathcal{X})$ is bounded in the interval $[1-\epsilon, 1+\epsilon]$. Moreover, in \cite{dasgupta2003elementary} it was proved that if the random projection is of Gaussian type, then a lower bound of the embedding dimension is given by $k\geq 4 (\epsilon ^2 /2 - \epsilon ^3 / 3)^{-1} \ln{n}$.

We note that the above mapping is a feature mapping, which in principle may result in a dimensionality reduction ($k<d$) or a projection into a higher-dimensional space ($k>d$) in which one seeks a linear manifold (in analogy to the case of kernel-based manifold learning methods). We also note that while the above linear random projection is but one of the choices for constructing a JL embedding (and proving it), it was experimentally demonstrated and/or theoretically proven that appropriately constructed nonlinear random embeddings may outperform simple linear random projections. For example, in \cite{giryes2016deep} it was shown that deep learning networks with random weights for each layer result in even better approximation accuracy than the simple linear random projection. Based on the concept of random projection, Schmidt et al. \cite{schmidt1992feed} performed computational experiments using FNNs with sigmoidal transfer functions in the hidden layer and a bias on the output layer, thus showing that by fixing the weights between the input and the hidden layer at random values, and by training the output weights via the solution of linear equations stemming from the Fisher minimization problem, the approximation accuracy is equivalent to that obtained with the standard back-propagation iterative procedure where all the weights are calibrated.

In the same year, RVFLNs were addressed in \cite{pao1992functional}. In RVFLNs the input layer is directly connected also to the output layer, the internal weights are chosen randomly in $[-1, 1]$ and the output weights are estimated in one step by solving a system of linear equations. By formulating a limit-integral of the function to be approximated with Monte Carlo simulations, in \cite{igelnik1995stochastic} it was proved that RVFLNs are universal approximators for continuous functions on bounded finite-dimensional sets and that the rate at which the approximation error vanishes is of the order of $\frac{1}{n}$, where $n$ is the number of basis functions parametrized by random variables. Extensive computational experiments performed in \cite{zhang2016comprehensive} showed that the direct links employed in RVFLNs between the input and the output layer play an important role for the performance, while the bias of the output neuron is not so important.\par
Reservoir Computing (otherwise called Eco State Networks) is another approach to network training based on the concept of random projection \cite{jaeger2002adaptive, verstraeten2007experimental, lukovsevivcius2009reservoir, butcher2013reservoir}. The basic structure of Reservoir Computing consists of a Recurrent Neural Network (RNN) with a large number of hidden units, an extra input and an extra output unit. Internal connection weights form a sparse random connectivity pattern with fixed values. The estimation of the optimal output weights is achieved by solving a system of linear equations arising from a mean-squares minimization problem.

Furthermore, it has been shown that single-layer FNNs with randomly assigned input weights and biases of the hidden layer and with infinitely  differentiable functions at the hidden layer, called ELMs, can universally approximate any continuous function on any compact input set \cite{huang2006extreme,Huang,huang2014insight,huang2015extreme}. For the case of an FNN with a single hidden layer of $h$ units, the random projection of the input space can be written as
\begin{equation}
    Y = \Phi(X),
\end{equation}
where the columns of the matrix $X \in \R^{d \times n}$ represent a set of $n$ points in the input $d$-dimensional space, the columns of $Y \in \R^{k \times n}$ are the corresponding random projections, and $\Phi : \R^{d \times n} \to \R^{k \times n}$ acts as an encoder, i.e.~a family of transfer functions whose parameters are sampled from a certain random distribution function. If the values of the weights $w_{ij}$ 
between the input and the hidden layer are fixed, then $\forall \boldsymbol{x} \in \mathbb{R}^d$ the random projection can be written as a linear map:
\begin{equation} \label{eqn:lin_map}
    Y = {W^{o}}\Phi, \quad Y \in \mathbb{R}^{k\times n},
\end{equation}
where $\Phi \in \mathbb{R}^{h \times n}$ is a random matrix containing the outputs of the nodes of the hidden layer as shaped by the $h$ random distribution functions (e.g., RBFs or sigmoidal functions)  and the $d$-dimensional inputs.
Thus, the so-called ELMs can be seen as underdetermined linear systems in which the output weights are estimated by solving minimum-norm least squares problems
\cite{huang2006extreme,huang2014insight,huang2015extreme,bjorck1996book}. Moreover, for ELMs, an interpolation-like theorem has been proven by Huang et al. in \cite{huang2006extreme}:
\begin{theorem}
Let us consider a single-hidden-layer FNN with $h$ hidden units and an infinitely differentiable transfer function $\phi: \R \rightarrow \R$, \ $n$ distinct input-output pairs ($\boldsymbol{x}_i, \boldsymbol{y}_i) \in \R^d \times \R^k$, and randomly chosen values from any continuous probability distribution for the internal weights $w_{ij}$ and for the values of the biases of the $h$ neurons of the hidden layer, grouped in $\boldsymbol{b} \in \R^{h}$. Let us also denote by $W \in \R^{h \times d}$ and $W^{o} \in \R^{k \times h}$ the matrices containing the internal and output weights, and by $X \in \R^{d \times n}$ and $Y \in \R^{k \times n}$ the matrices with columns $\boldsymbol{x}_i$ and $\boldsymbol{y}_i$. Then, the hidden layer output matrix $\Phi \in \R^{h \times n}$, whose elements are determined by the action of the transfer functions on $W X + \boldsymbol{b} \, \boldsymbol{1}^T$ (where  $\boldsymbol{1} \in \R^n$ has all the entries equal to 1), has full rank and
$$
   \lVert {W^{o}} \Phi - Y \rVert = 0
$$
with probability  1.
\end{theorem}
A review on neural networks with random weights can be found in \cite{cao2018review}.

\section{The Proposed Method\label{sec:our_method}}
Here we focus on the numerical solution of problem \eqref{eq:ODE_PB} in an interval, say $[x_0, x_{end}]$. According to the previous notation, for this problem we have $d=1$ and $k=m$. Thus, we denote by $\Psi_i(x,\bw^{o}_i,\boldsymbol{p}_i)$ the $i$-th component of the trial solution $\boldsymbol{\Psi}(x,W^o,P)$ defined in \eqref{eqn:trial_vec}, where $W^o \in \R^{m \times h}$ is the matrix containing the weights $w^{o}_{ij}$ between the hidden and the output layer. Note that we separate $W^0$ from the other network parameters consisting of the input weights, the biases of the hidden layer and the parameters of the transfer function, and denote those parameters again by $P$, with a small abuse of the notation. Furthermore, the vector of hyperparameters $\boldsymbol{q}$ in \eqref{eqn:trial_vec} is neglected because we do not use it.
Following \cite{lagaris1998artificial} and taking into account that the trial solution must satisfy the initial value conditions of the problem, i.e. $y_i(x_0)=\alpha_i$, $i=1,2,\dots,m$, we set
\begin{equation}
    \Psi_i(x,\bw^{o}_i,\boldsymbol{p}_i) = \alpha_i+(x-x_0) N_i(x,\bw^{o}_i,\boldsymbol{p}_i),
    \label{trsolsys}
\end{equation}
where $N_i(x,\bw^{o}_i,\boldsymbol{p}_i)$ is a single-output FNN with parameters the output weights $\bw^{o}_i=[w^{o}_{1i} \; w^{o}_{2i} \; \ldots \; w^{o}_{hi}]^T \in \R^{h}$, and $\boldsymbol{p}_i$ contains the remaining parameters associated with that network.
Then, if one considers the numerical solution based on $n$ collocation points $x_1, x_2,\ldots, x_n$, then the error function that we seek to minimize for training the FNN is given by
\begin{equation}
\begin{aligned}
    & E(W^o,P) = \\
    & \sum_{i=1}^m \sum_{j=1}^n \left( \frac{d\Psi_i}{dx}
    (x_j,\bw^{o}_i,\boldsymbol{p}_i) - f_i(x_j, \Psi_1(x_j,\bw^{o}_1,\boldsymbol{p}_1),
    \ldots, \Psi_m(x_j,\bw^{o}_m,\boldsymbol{p}_m) \right)^2 \!\!.
\end{aligned}
\label{errorsys}
\end{equation}
Here we propose a machine learning method based on random projections for the solution of IVPs of ODEs in $n$ input collocation points. In particular, we employ $m$ (one for each unknown variable of the systems of ODEs) SISO neural sub-networks with a linear output transfer function, with a single hidden layer having $h$ nodes with Gaussian RBFs.
In particular, we consider each sub-network $N_i$ to be a linear combination of RBFs, 
\begin{equation}
\label{eq:Ni_sum}
 N_i(x,\bw^{o}_i,\boldsymbol{p}_i)=\sum_{j=1}^h w^{o}_{ji} G_{ji}(x),
 \quad i=1,2, \dots,m ,
\end{equation}
where
\begin{align}
\label{eq:Gji}
   G_{ji}(x) & = G (w_{ji}x + b_{ji}, c_{j}, \sigma_{ji}) = \exp \left( -\dfrac{( w_{ji}x+b_{ji}-c_{j})^2}{\sigma_{ji}^2} \right), \\
   & \hspace*{45pt} j=1,2,\ldots h, \; i=1,2,\ldots m, \nonumber
\end{align}
with $w_{ji}=1$ and $c_{j} = x_0 + (j-1) \bar{s}$ with $\bar{s} = (x_{end}-x_0)/(h-1)$ and $j=1,\ldots,h$.

Under the above assumptions, the derivative of the $i$-th component $\Psi_i$ of the trial solution with respect to the collocation point $x_l$
is given by (see \eqref{eqn:gradnet})
\begin{equation}
    \frac{\partial \Psi_i}{\partial x_l} =
    N_i(x_l,\bw^{o}_i, \boldsymbol{p}_i)-(x_l-x_0)\sum_{j=1}^{h} \frac{2}{\sigma_{ji}^2} w^{o}_{ji} (x_l+b_{ji}-c_{j}) \exp \left( -\dfrac{( x_l+b_{ji}-c_{j})^2}{\sigma_{ji}^2} \right),
\label{eq:dirGaussian_RBF}
\end{equation}  
while, for any fixed $x_l$, the derivative of $\Psi_i$ with respect to the only unknown parameter $w^{o}_{ji}$ is given by (see \eqref{trsolsys} and \eqref{dngdw})
\begin{equation}
    \frac{\partial \Psi_i}{\partial w^{o}_{ji}} = (x_l-x_0) \exp \left( -\dfrac{( x_l+b_{ji}-c_{ji})^2}{\sigma_{ji}^2} \right).
    \label{eqn:gradnetrbf2}
\end{equation}

For the determination of the randomized parameters of the Gaussian RBFs, we fix reference intervals where we look for these parameters. Our requests are that the functions are neither too steep nor too flat in the reference interval, and at each collocation point there are at least two basis functions giving values that are not too small. Then, based on numerical experiments,
the biases $b_{ji}$ of the hidden units and the parameters $1/\sigma_{ji}^2$ are taken to be uniformly randomly distributed in the intervals
$$
   \left[ -\frac{(x_{end}-x_0)}{6}, 0 \right] \mbox{ and } \left[ \dfrac{3}{8(x_{end}-x_0)^2}, \dfrac{81}{2(x_{end}-x_0)^2} \right],
$$
respectively. We emphasise that this choice appears to be problem independent. Therefore, the only parameters that have to be determined by training the network are the output weights $w^o_{ji}$.
Hence, for the $n$ collocation points $x_l$, the outputs of each network ${N}_i$, $i=1,2,\dots m$, are given by:
\begin{equation}
    \boldsymbol{N}_i(x_1,x_2, \dots x_n, \bw^{o}_{i}, \boldsymbol{p}_i) = R_i \bw^{o}_{i},
    \label{eqn:rpnn}
\end{equation}
where $\boldsymbol{N}_i(x_1,x_2, \dots x_n, \bw^{o}_{i}, \boldsymbol{p}_i) \in \R^{n}$ is the vector with $l$-th component the output of $N_i$ corresponding to $x_l$, and $R_i=R_i(x_1,\ldots,x_n,\boldsymbol{p}_i) \in \R^{n\times h}$ is defined as
\begin{equation}
   R_{i}(x_1,\ldots,x_n,\boldsymbol{p}_i)=
   \begin{bmatrix}
       G_{1i}(x_1) & \cdots & G_{hi}(x_1)\\
       \vdots        & \vdots  & \vdots \\
       G_{1i}(x_n) & \cdots & G_{hi}(x_n)
   \end{bmatrix} .
\end{equation}
The minimization of the error function given in \eqref{errorsys} is performed by a Gauss-Newton scheme (see, e.g., \cite{kelley1999iterative}) over 
$nm$ nonlinear residuals $F_q$, with $q=l + (i-1)n$, $i=1,2,\dots m$, $l=1,2,\dots n$, given by
\begin{equation}
    F_q(\boldsymbol{W}^o)=
\frac{d\Psi_i}{dx_l} (x_l,\bw^{o}_i,\boldsymbol{p}_i) - f_i(x_l, \Psi_1(x_l,\bw^{o}_1,\boldsymbol{p}_1),
\ldots, \Psi_m(x_l,\bw^{o}_m,\boldsymbol{p}_m)),
\label{eqn:Fq}
\end{equation}
where $\boldsymbol{W}^0 \in \R^{m h}$ is the column vector stacking the values of all the $m$ vectors $\boldsymbol{w}^{o}_i \in \mathbb{R}^h$,
\begin{align*}
\boldsymbol{W}^o =\begin{bmatrix}
\boldsymbol{w}^{o}_1\\
\boldsymbol{w}^{o}_2\\
\vdots \\
\boldsymbol{w}^{o}_m\\
\end{bmatrix} \in \R^{m h}.
\end{align*}
Thus, by setting $\boldsymbol{F}(\boldsymbol{W}^o) = [F_1(\boldsymbol{W}^o) \cdots F_q(\boldsymbol{W}^o) \cdots F_{(nm)}(\boldsymbol{W}^o)]^T$, 
the Gauss-Newton method reads:
\begin{equation}
\begin{aligned}
\boldsymbol{W}^{o(\nu+1)} & =  \displaystyle \boldsymbol{W}^{o(\nu)} + d\boldsymbol{W}^{o(\nu)}, \\
\displaystyle d\boldsymbol{W}^{o(\nu)}  & =  \displaystyle \arg\min_{\boldsymbol{W}^{o(\nu)}} \| ( \nabla_{\boldsymbol{W}^{o}} \boldsymbol{F}(\boldsymbol{W}^{o(\nu)}))^T d\boldsymbol{W}^{o(\nu)} + \boldsymbol{F}(\boldsymbol{W}^{o(\nu)}) \|,
\end{aligned}
\label{gauss-newton}
\end{equation}
where $(\nu)$ denotes the current iteration and $\nabla_{\boldsymbol{W}^{o(\nu)}}\boldsymbol{F} \in \R^{n m \times m h }$ is the Jacobian matrix of $\boldsymbol{F}$ with respect to $\boldsymbol{W}^{o (\nu)}$:
\begin{equation}
\nabla_{\boldsymbol{W}^{o(\nu)}} \boldsymbol{F}(\boldsymbol{W}^{o(\nu)}) =
\begin{bmatrix}
\frac{\partial F_1}{\partial W^{o}_{1}} & \frac{\partial F_1}{\partial W^{o}_{2}} & \dots & \frac{\partial F_1}{\partial W^{o}_{p}} & \dots &\frac{\partial F_1}{\partial W^{o}_{(mh)}}\\
\frac{\partial F_2}{\partial W^{o}_{1}} & \frac{\partial F_2}{\partial W^{o}_{2}} & \dots & \frac{\partial F_2}{\partial W^{o}_{p}} & \dots &\frac{\partial F_2}{\partial W^{o}_{(mh)}}\\
\vdots  &\vdots & \ddots & \vdots & \ddots &\vdots\\
\frac{\partial F_q}{\partial W^{o}_{1}} & \frac{\partial F_q}{\partial W^{o}_{2}} & \dots & \frac{\partial F_q}{\partial W^{o}_{p}} & \dots &\frac{\partial F_q}{\partial W^{o}_{(mh)}}\\
\vdots  &\vdots & \ddots & \vdots & \ddots &\vdots\\
\frac{\partial F_{(nm)}}{\partial W^{o}_{1}} & \frac{\partial F_{(nm)}}{\partial W^{o}_{2}} & \dots & \frac{\partial F_{(nm)}}{\partial W^{o}_{p}} & \dots &\frac{\partial F_{(nm)}}{\partial W^{o}_{(mh)}}
\end{bmatrix}_{\big|{({\boldsymbol{W}^o}^{(\nu)})}}.
\label{jac}
\end{equation}
Note that the residuals depend on the derivatives $\frac{\partial \Psi_i}{\partial x_l}$ in \eqref{eq:dirGaussian_RBF} and the trial functions $\Psi_i$ in \eqref{trsolsys}, while the elements of the Jacobian matrix depend on the derivatives of $\frac{\partial \Psi_i}{\partial w^{o}_{ji}}$ in \eqref{eqn:gradnetrbf2} as well as on the mixed derivatives $\frac{\partial^2 \Psi_i}{\partial x_l \partial w^{o}_{ji}}$. Based on \eqref{eq:dirGaussian_RBF}, the latter are given by
\begin{equation}
     \frac{\partial^2 \Psi_i}{\partial x_l \partial w^{o}_{ji}}=\frac{\partial N_i(x_l,\bw^{o}_i, \boldsymbol{p}_i)}{\partial w^{o}_{ji}}-(x_l-x_0) \frac{2}{\sigma_{ji}^2} (x_l+b_{ji}-c_{j}) \exp \left( -\dfrac{(x_l+b_{ji}-c_{j})^2}{\sigma_{ji}^2} \right),
    \label{mixedGaussian_RBF}
\end{equation}  
where 
\begin{equation}
\frac{\partial N_i(x_l,\bw^{o}_i, \boldsymbol{p}_i)}{\partial w^{o}_{ji}}=\exp \left( -\dfrac{(x_l+b_{ji}-c_{j})^2}{\sigma_{ji}^2} \right).
\label{derivN}
\end{equation}
Note also that if the Jacobian of system~\eqref{eq:ODE_PB}, with elements $\frac{\partial f_i(x_l)}{\partial y_k}$, is given, then the elements of the Jacobian matrix $\nabla_{\boldsymbol{W}^{o(\nu)}}\boldsymbol{F}$ can be computed easily:
\begin{equation}
    \frac{\partial F_p}{\partial W^o_q}=\frac{\partial^2 y_i}{\partial x_l\partial w^{o}_{jk}}-\frac{\partial f_i(x_l)}{\partial w^{o}_{jk}}=\frac{\partial^2 \Psi_i}{\partial x_l \partial w^{o}_{jk}}-(x-x_0)\frac{\partial f_i(x_l)}{\partial y_{k}}\frac{\partial N_k(x_l,\bw^{o}_i, \boldsymbol{p}_i)}{\partial w^{o}_{jk}},
    \label{eq:jac_J}
\end{equation}
where, as before, $q=l+(i-1)n$ and $p=j+(k-1)h$.

As in general $h>n$, i.e. the minimization problem in \eqref{gauss-newton} is an underdetermined linear system,
we compute its solution with minimum $L_2$ norm. This can be performed by computing the Moore-Penrose pseudoinverse of the Jacobian with the Singular Value Decomposition. More precisely, we estimate the pseudoinverse by cutting off all the singular values (and their corresponding singular vectors) below a small tolerance $\epsilon$. This allows us to deal with the case where the Jacobian is rank deficient, which may happen because the RBFs are not injective functions. Furthermore, this choice also allows us to cope with the difference between the exact rank and the numerical rank of the Jacobian matrix. Following \cite{golub1996matrix}, for some small $\epsilon$ the $\epsilon$-rank of a matrix $M \in \R^{mn \times mh}$ is defined as follows:
\begin{equation}
    r_{\epsilon} = \min \{ \mbox{rank}(B) \in \R^{mn \times mh} : \| M - B \|_{L_2} \le \epsilon \}.
\end{equation}
Then, if $\epsilon$ is ``small enough'', we neglect all the singular values below $\epsilon$ and approximate the pseudoinverse of the Jacobian matrix as
\begin{equation}
 (\nabla_{W^o}\boldsymbol{F})^{+} = V_{r_{\epsilon}} \Sigma_{r_{\epsilon}}^{+} U_{r_{\epsilon}}^T,
    \label{eq:Moore-Penrose_PseudoInv}
\end{equation}
where $\Sigma_{r_{\epsilon}} \in \R^{r_{\epsilon} \times r_{\epsilon}} $ is the diagonal matrix with the ${r_{\epsilon}}$ largest singular values of $\nabla_{W^o}\boldsymbol{F}$, $\Sigma_{r_{\epsilon}}^{+}$ is its pseudoinverse, and $U_{r_{\epsilon}} \in \mathbb{R}^{mn \times {r_{\epsilon}}}$ and $V_{r_{\epsilon}} \in \mathbb{R}^{mh \times {r_{\epsilon}}}$ are the matrices with columns the corresponding ${r_{\epsilon}}$ left and right eigenvectors, respectively. The value of $\epsilon$ used in our experiments is specified at the beginning of Section~\ref{sec:results}. Thus, the direction $dW^{o(\nu)}$ at each Gauss-Newton iteration is given by
\begin{equation*}
    dW^{o(\nu)} = -V_{r_{\epsilon}} \Sigma_{r_{\epsilon}}^+ U_{r_{\epsilon}}^T \boldsymbol{F}(W^{o(\nu)}).
\end{equation*}
Let us also observe that the Jacobian matrix $J \in \mathbb{R}^{m\times m}$ of the system of ODEs \eqref{eq:ODE_PB} is given:
\begin{equation}
    J=\begin{bmatrix}
\frac{\partial f_1}{\partial y_1} & \frac{\partial f_1}{\partial y_{2}} & \dots &\frac{\partial f_1}{\partial y_m}\\
\frac{\partial f_2}{\partial y_1} & \frac{\partial f_2}{\partial y_{2}} & \dots &\frac{\partial f_2}{\partial y_m}\\
\vdots  &\vdots & \ddots &\vdots\\
\frac{\partial f_m}{\partial y_1} & \frac{\partial f_m}{\partial y_{2}} & \dots &\frac{\partial f_m}{\partial y_m}
\end{bmatrix}.
\label{jac_ode}
\end{equation}


\subsection{Convergence}

We note that in the above configuration we fix all the internal weights to 1 and the centers to be equidistant in the domain, while we set randomly, from appropriately chosen uniform distributions, the biases $b_{ji}$ and the width parameters $\sigma_{ji}$ of the RBFs. This configuration is slightly different from the classical ELMs with RBFs, where the internal weights are set equal to 1 and the biases equal to 0, while the centers are uniformly randomly distributed in $[-1, 1]$ (see, e.g., \cite{huang2006universal}). 
However, it is straightforward to show that in our scheme we still get universal approximation properties, as stated in the next theorem.
\begin{theorem} \label{th:univ-approx}
Let us fix a continuous (target) function $\varphi$ and consider the system of functions $\{ N_i(x,\bw^{o}_i,\boldsymbol{p}_i) \}_{i=1,\dots,m}$ in \eqref{eq:Ni_sum}, with $G_{ji}$, $w_{ji}$ and $c_{j}$ defined in \eqref{eq:Gji} and the subsequent lines. Then,
for every sequence of randomly chosen parameters $\boldsymbol{p}_i$ there exists a choice of $\bw^{o}_i$ such that
\[
   \lim_{m \to \infty} \| N_i(x,\bw^{o}_i,\boldsymbol{p}_i)-\varphi(x) \| = 0 \text{ with probability 1.}
\]
\end{theorem}
\noindent \emph{Proof.}
Since $w_{ji}$ and $c_{j}$ are fixed, one can manipulate \eqref{eq:Gji} and write it as
\[
   G_{ji}(x)= \exp \left( -\dfrac{( x+ \alpha_{ji} )^2}{\beta_{ji}} \right),
\]
where the parameters $\alpha_{ji}=b_{ji}-c_{j}$ and $\beta_{ji}=\sigma_{ji}^2$ are random variables (because $b_{ji}$ and $\sigma_{ji}$ are random variables sampled from continuous probability distributions). Therefore, the network fits the hypotheses of an FNN with random hidden nodes. Because the considered RBFs are sufficiently regular, Theorem II.1 in \cite{huang2006universal} holds and hence the thesis follows in a straightforward manner.\qed

\subsection{The Adaptive Scheme} \label{Sec:adapt}
For the numerical solution of IVPs of ODEs over long intervals in the presence of steep gradients and/or stiff dynamics, the proposed method can be iterated by dividing the whole interval $[x_0,x_{end}]$ into sub-intervals, i.e. $[x_0,x_{end}] = [x_0,x_1] \cup [x_1,x_2] \cup \dots,\cup[x_k,x_{k+1}]\cup \dots \cup [x_{end-1},x_{end}]$, where $x_1,x_2,\dots,x_k,\dots,x_{end-1}$ are determined in an adaptive way as explained below. This decomposition of the interval leads to the solution of several consecutive IVPs:
\begin{equation}\label{eq:ODE_PB_iter}
   \left\{ \begin{array}{lll}
       \displaystyle \dfrac{d y_i^{(k)}}{dx} & = & \displaystyle f_i(x, y_1^{(k)},y_2^{(k)},...,y_m^{(k)}), \\[7pt]
       \displaystyle y_i^{(k)}(x_k) & = & \displaystyle \alpha_i^{(k)}, 
   \end{array} \right.
   \quad i = 1,2,\ldots,m.
\end{equation}
In the $k$-th interval $[x_k,x_{k+1}]$, the initial value $\alpha_i^{(k)}$ of $y_i^{(k)}$ is set equal to the final value of $y_i^{(k-1)}$ in the $(k-1)$-st interval:
\begin{equation}
    \alpha_i^{(k)}=y_i^{(k-1)}(x_k).
    \label{eq:iter_initialvalues}
\end{equation}
According to \eqref{trsolsys}, the $i$-th component of the trial solution in the $k$-th interval reads: 
\begin{equation*}
    \Psi_i^{(k)}(x,\bw^{o}_i,\boldsymbol{p}_i) = y_i^{(k-1)}(x_k)+(x-x_k) N_i^{(k)}(x,\bw^{o}_i,\boldsymbol{p}_i).
\end{equation*}
Thus, in each $k$-th interval, the approximation of $y_i^{(k)}$ is computed by adjusting the external weights of the RPNNs $N_i^{(k)}$ using the Gauss-Newton scheme in \eqref{gauss-newton}-\eqref{jac}.

In order to describe the adaptive scheme, let us assume we have solved the problem up to the interval $[x_{k-1},x_{k}]$, hence we have found $y_i^{(k-1)}$ and we are seeking $y_i^{(k)}$ in the current interval $[x_k,x_{k+1}]$ with a given width $\Delta x_k = x_{k+1}-x_k$. If in the current interval the Gauss-Newton scheme does not achieve a prescribed accuracy within a number of iterations (here set to 4), then the interval width is halved, thus redefining a new guess $x^*_{k+1}$ for $x_{k+1}$:
\begin{equation}
x_{k+1}\leftarrow x^*_{k+1}=x_{k}+\dfrac{\Delta x_k}{2},
\end{equation}
and the scheme is applied again in the interval $[x_k,x^*_{k+1}]$. Otherwise, if the Gauss-Newton method converges, the initial guess $\Delta x_{k+1}$ of the width of the next interval is set by doubling the width of the current interval, so that
\begin{equation}
x_{k+2}=x_{k+1}+2\Delta x_k.
\end{equation}
We summarize the proposed method in the pseudo-code shown in Algorithm~\ref{alg:rpnn}, where $\mathcal{U}(a,b)$ denotes the uniform random distribution in the interval $[a,b]$ and $tol$ is a tolerance used in the Gauss-Newton stopping criterion.
In all our calculations, the initial guess $\Delta x_0$ for the first interval is the width of the whole interval $[x_0,x_{end}]$.

\SetKwComment{Comment}{/* }{ */}
\begin{algorithm}[htp!]
\setstretch{1.15}
\caption{solving an IVP of ODE systems using an RPNN with RBFs\label{alg:rpnn}}
\SetAlgoLined
\smallskip
\textbf{IVP problem} \\
$\displaystyle \frac{dy_i}{dx}=f_i(x, y_1,y_2,\ldots,y_m)$ in $[x_0, x_{end}]$\;
$y_i(x_0)=\alpha_i, \;\, i=1,2,\ldots,m$\; 
\textbf{INITIALIZATION}\\
$h \gets 40; \; n \gets 20$\Comment*[r]{set \# neurons and collocation points}
$maxiter\gets4$\Comment*[r]{set max \# iters for Gauss-Newton method}
$\Delta x \gets (x_{end}-x_0); \;\, x^{*} \gets x_0$\Comment*[r]{set initial interval}
\textbf{ADAPTIVE SCHEME}\\
\Repeat{$x^{*}=x_{end}$}{
Select $x_l \in [x^{*}, x^{*}+\Delta x], \;\, l=1,\ldots,n$\Comment*[r]{set collocation points}
\textbf{STEP 1: fix parameters, weights and biases of the RPNN} \\ 
$c_{j} \gets x^{*} + j \dfrac{\Delta x}{h-1}, \;\, j=1,\dots,h$\Comment*[r]{set RBF centers}
$\displaystyle b_{ji} \sim \mathcal{U} \left(-\dfrac{\Delta x}{6}, 0 \right)\!, \;\, j=1,\ldots,h, \;\, i=1,\ldots,m$\Comment*[r]{set biases}
$\displaystyle \dfrac{1}{\sigma_{ji}^2} \sim \mathcal{U} \left( \frac{3}{8\Delta x^2},\frac{81}{2\Delta x^2} \right)\!, \;\, j, i \mbox{ as above}$\Comment*[r]{set RBF inverse widths}
\textbf{STEP 2: set trial functions} \\
$N_i(x,\boldsymbol{w}^0_i,\bp_i) \gets \sum_{j=1}^h w^{o}_{ji}\text{exp}\left(-\frac{(x-c_j)^2}{\sigma_{ji}^2}-b_{ji}\right)$\Comment*[r]{see \eqref{eq:Ni_sum}-\eqref{eq:Gji}}
$\Psi_{i}(x,\boldsymbol{w}^0_i,\bp_i)=\alpha_i+(x-x^{*})N_i(x,\boldsymbol{w}^0_i,\bp_i)$\;
\textbf{STEP 3: compute output weights $W^0$}\\
$iter \gets 0$\;
\Repeat{$(err \le tol)$ or $(iter \ge maxiter)$}{
\For{$j=1,\dots,n, \;\, i=1,\dots,m$}{
$q \gets j + (i-1)n$\;
$\displaystyle F_q(W^o) \gets \dfrac{d \Psi_i}{d x_l}(x_l,\bw^{o}_i,\boldsymbol{p}_i)-$ \\ 
$ \qquad f_i(x_l, \Psi_1(x_l,\bw^{o}_1,\boldsymbol{p}_1), \ldots, \Psi_m(x_l,\bw^{o}_m,\boldsymbol{p}_m))$\; 
}
Set $\boldsymbol{F}(W^o)=[F_1(W^o) \; F_2(W^o) \; \dots \; F_{(nm)}(W^o)]^T$\;
Compute Jacobian matrix $\nabla_{W^o} \boldsymbol{F}(W^o)$\Comment*[r]{see \eqref{jac}-\eqref{eq:jac_J}}
$(\nabla_{W^o}\boldsymbol{F})^{+} \gets V_{r_{\epsilon}} \Sigma_{r_{\epsilon}}^{+} U_{r_{\epsilon}}^T$\Comment*[r]{compute pseudo-inverse}
$W^o \gets W^o - (\nabla_{W^o}\boldsymbol{F}(W^o))^{\dagger}\boldsymbol{F}(W^o)$\Comment*[r]{update weights}
$err \gets ||\boldsymbol{F}(W^o)||_2$\Comment*[r]{compute error}
$iter \gets iter + 1$\;
} 
\eIf{$err>tol$\Comment*[r]{if Gauss-Newton does not converge}}
{$\Delta x \gets\dfrac{\Delta x}{2}$\Comment*[r]{halve interval width}}
{$\Delta x \gets 2 \Delta x$\Comment*[r]{double width of next time interval}
$x^{*} \gets x^{*}+\Delta x$\Comment*[r]{update $x^{*}$}
} 
} 
\end{algorithm}

\section{Numerical Results\label{sec:results}}

We implemented Algorithm~\ref{alg:rpnn} using MATLAB 2020b on an Intel Core i7-10750H CPU @ 2.60GHz with up to 3.9 GHz frequency and a memory of 16 GBs. The Moore-Penrose pseudoinverse of $\nabla_{\boldsymbol{w}^o}\boldsymbol{F}$ was computed with the MATLAB built-in function \texttt{pinv}, with the default tolerance. For our simulations, we chose four well-known and challenging stiff ODE problems: Prothero-Robinson \cite{prothero1974stability}, van der Pol \cite{van1926lxxxviii}, ROBER \cite{Robertson1966} and HIRES problem \cite{schafer1975new}. For comparison purposes, we also solved the ODE problems with two widely-used MATLAB built-in functions, namely \texttt{ode45}, implementing an adaptive-step Runge-Kutta method based on the Dormand-Prince pair, and \texttt{ode15s}, implementing a variable-step variable-order multistep method based on numerical differentiation formulas. In order to estimate the error in the approximate solution, we used as reference solution the one computed by \texttt{ode15s} with absolute and relative error tolerances equal to 1e$-$14. To this aim, we computed the $L_2$ and $L_{\infty}$ norms of the differences between the computed and the reference solutions, as well as the Mean Absolute Error (MAE). These performance metrics were evaluated at equidistant collocation points, selecting their number according to the problem (see below). Finally, we ran each solver 10 times and computed the median, maximum and minimum execution times in seconds. The time of each run was measured by using the MATLAB commands \texttt{tic} and \texttt{toc}.

Henceforth, we use $t$ instead of $x$, since the independent variable in the test problems represents time. We also assume $t_0 = x_0 = 0$. In our simulations, when implementing the proposed machine learning method, we fixed as $n=20$ the number of points where the solution was sought in a specific interval. The initial time interval was the whole interval, as described in Section \ref{Sec:adapt}.

\subsection{Case Study 1: Prothero-Robinson problem}

The Prothero-Robinson ODE \cite{prothero1974stability} is given by
\begin{equation}
    \frac{dy}{dt}=\lambda(y-\phi(t))+\phi'(t), \quad \lambda<0.
    \label{pr}
\end{equation}
Its solution is $\phi(t)$ and the initial condition is $y(0)=\phi(0)$. The problem becomes stiff for $\lambda\ll0$.
For our numerical simulations, we choose $\phi(t)=sin(t)$, $u(0)=\phi(0)=sin(0)=0$, and $[0,2\pi]$ as the time interval where the solution is sought, while the parameter $\lambda$ controlling the stiffness is set equal to 1e$-$05. For the implementation of the proposed approach, we use the following (initial) trial solution:
\begin{equation}
    \Psi (t, \bw^o) = \alpha^{(0)} + t N(t,\bw^o), \quad \alpha^{(0)}=y(0)=sin(0)=0.
    \label{trsolpr}
\end{equation}
\begin{figure}[ht]
    \centering
    \subfigure[]{
        \includegraphics[width=0.45\textwidth]{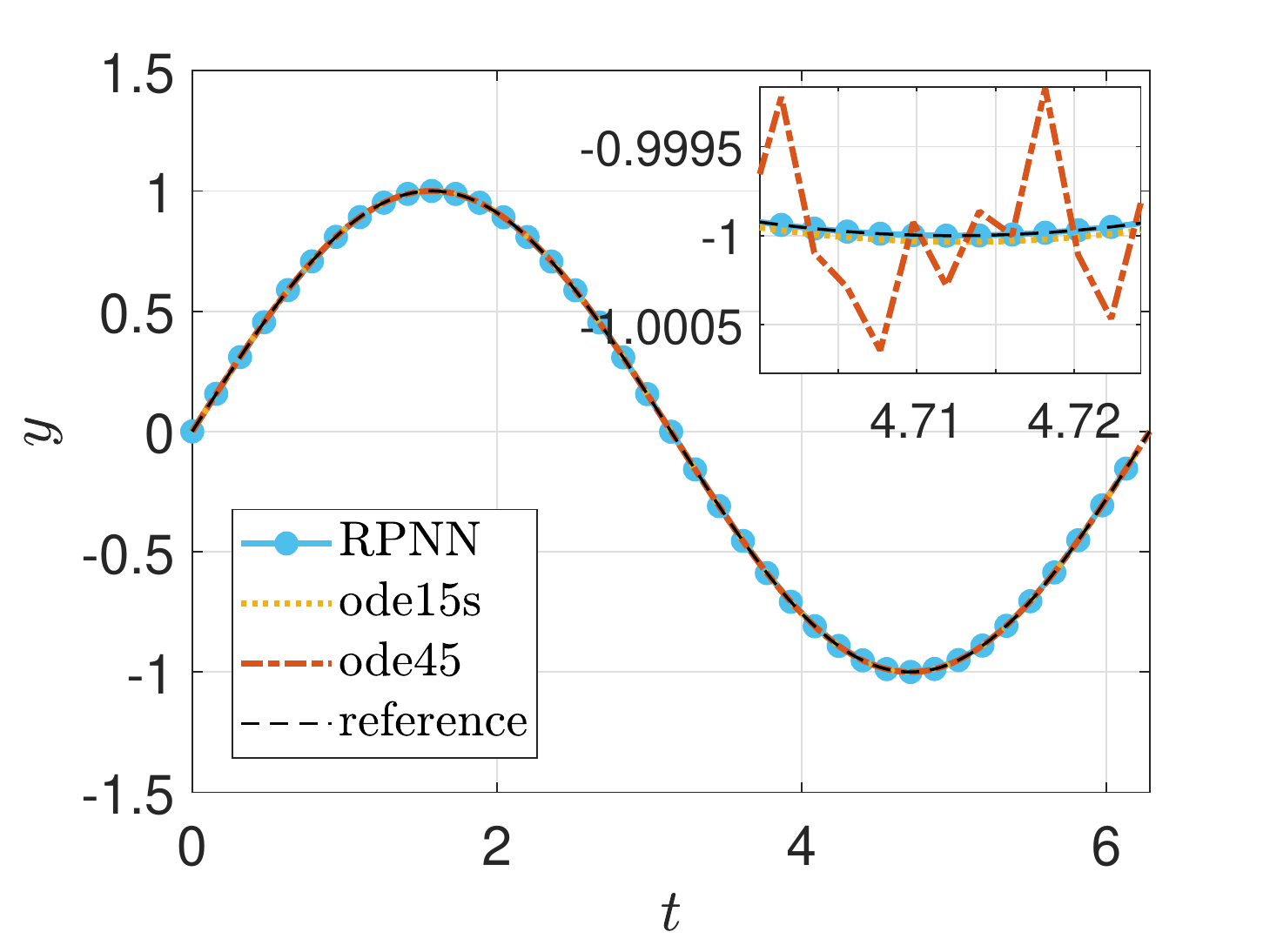}
    }
    \subfigure[]{
        \includegraphics[width=0.45\textwidth]{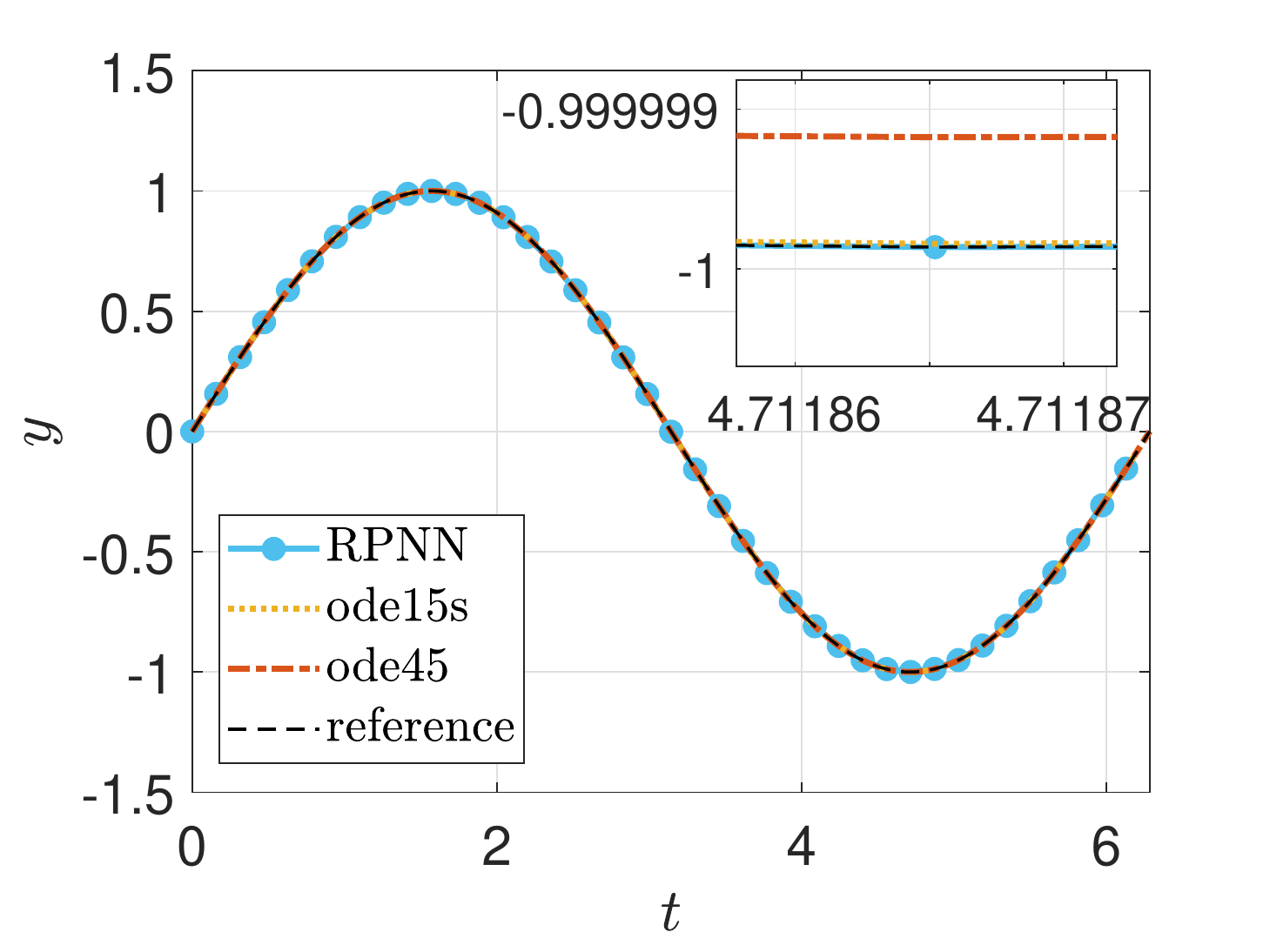}
    }
    \caption{Numerical solutions of the Prothero-Robinson problem with $\lambda =$ 1e$-$05 in the interval $[0,2\pi]$: (a) $tol = 1$e$-$03, (b) $tol = 1$e$-$06. The reference (analytical) solution is $sin(t)$. The insets depict a zoom around the reference solution.\label{fig:prothero-robinson}}
\end{figure}

In Figure~\ref{fig:prothero-robinson}, we show the approximate solutions obtained with tolerances 1e$-$03 and 1e$-$06. The insets depict a zoom around the reference solution $sin(t)$. As it is shown, when the absolute and relative tolerances of \texttt{ode45} are equal to 1e$-$03, this function ``oscillates'' around the reference solution, while the solutions provided by \texttt{ode15s} and the proposed RPNN are indistinguishable from the reference solution (see Figure~\ref{fig:prothero-robinson}(a)). When the tolerance is 1e$-$06, the proposed RPNN outperforms both \texttt{ode45} and \texttt{ode15s}, providing a solution that is indistinguishable from the reference solution, as can be seen in the inset.
\begin{table}[ht]
\begin{center}
\caption{Prothero-Robinson problem with $\lambda = 1$e$-$05 in the interval $[0,2\pi]$. Absolute errors ($L_{2}$-norm, $L_{\infty}$-norm and MAE) for the solutions computed with tolerances 1e$-$03 and 1e$-$06. The reference solution is the analytical solution $sin(t)$.\label{tab:prothero_accuracy}}
{\small
\begin{tabular}{|l|lll|lll|}
\hline
& \multicolumn{3}{c|}{$tol=$ 1e$-$03} & \multicolumn{3}{c|}{$tol=$ 1e$-$06} \\
\cline{2-7}
& $L_2$ & $L_{\infty}$ & MAE & $L_2$ & $L_{\infty}$ & MAE\\
\hline
RPNN   & 5.95e$-$08 & 6.12e$-$09 & 3.31e-10 & 2.23e$-$08 & 2.84e$-$09 & 9.36e-11\\
\texttt{ode45}  & 2.24e$-$02 & 1.13e$-$03 & 2.97e$-$04 & 2.19e$-$05 & 1.14e$-$06 & 2.89e$-$07\\
\texttt{ode15s} & 7.36e$-$03 & 5.86e$-$04 & 1.01e$-$04 & 5.65e$-$06 & 6.14e$-$07 & 6.08e$-$08 \\
\hline
\end{tabular}
}
\end{center}
\end{table}

In Table~\ref{tab:prothero_accuracy}, we report the numerical approximation accuracy obtained with the various methods in terms of the $L_2$-norm and $L_{\infty}$-norm of the error and of MAE with respect to the reference solution. In order to compute these errors, we evaluated the corresponding solutions taking 3000 equidistant points in $[0,2\pi]$. For \texttt{ode45} and \texttt{ode15s}, the values of the approximate solutions at these points were obtained using the MATLAB function \texttt{deval}, which uses an interpolation technique.
\begin{table}[ht]
\begin{center}
\caption{Prothero-Robinson problem with $\lambda = 1$e$-$05  in the interval $[0,2\pi]$. Computational times in seconds (median, minimum and maximum times over 10 runs) and number of points required with tolerances~1e$-$03 and~1e$-$06. The reference solution is the analytical solution $sin(t)$.\label{tab:prothero_time_points}}
{\small
\setlength{\tabcolsep}{3pt}
\begin{tabular}{|l|lllr|lllr|}
\hline
& \multicolumn{4}{c|}{$tol=$ 1e$-$03} & \multicolumn{4}{c|}{$tol=$ 1e$-$06} \\
\cline{2-9}
& median & min & max & \multicolumn{1}{l|}{\# pts} & median & min & max & \multicolumn{1}{l|}{\# pts} \\
\hline
RPNN & 5.52e$-$04 & 4.56e$-$04 & 1.88e$-$03 &  20  & 6.21e$-$04 & 5.27e$-$04 & 6.87e$-$03 & 40\\
\texttt{ode45}  & 4.96e$-$01 & 4.91e$-$01 & 5.52e$-$01 & 189298 & 5.07e$-$01  & 4.93e$-$01 & 5.16e$-$01  & 189308\\
\texttt{ode15s} & 7.87e$-$04 & 7.45e$-$04 & 1.26e$-$03 & 28 & 1.38e$-$03 & 1.22e$-$03 & 1.85e$-$03 & 73\\
\hline
\end{tabular}
}
\end{center}
\end{table}
As it is shown, the proposed numerical scheme outperforms both \texttt{ode45} and \texttt{ode15s} in all metrics for both tolerances.

Finally, in Table~\ref{tab:prothero_time_points} we report the computational times and the number of points required by each method, for both tolerance values. Notably, our method outperforms both \texttt{ode45} and \texttt{ode15s}, since it results in significantly smaller computational times and number of points than \texttt{ode45}, and (on average) in slightly smaller computational times and number of points than \texttt{ode15s}. We note that \texttt{ode45} shows ``poor'' performance due to the very high stiffness, thus requiring a huge number of points and yet resulting in relatively large approximations errors.

\subsection{Case Study 2: van der Pol problem}

The van der Pol model and the concept of the so-called relaxation oscillations was introduced by Balthazar van der Pol \cite{van1926lxxxviii} to describe the nonlinear oscillations observed in a triode circuit characterized by a slow charge of a capacitor followed by a very rapid discharge. The model is given by the following equation:
\begin{equation}
    y_1''-\mu(1-y_1^2)y_1'+y_1=0,
    \label{vdp}
\end{equation}
where $\mu>0$ is a scalar parameter. The problem becomes stiff for $\mu \gg 1$. By setting $y_1' = y_2$, we obtain the following first-order ODE system:
\begin{equation}
\begin{aligned}
y_1' & = y_2, \\
y_2' & = \mu(1-y_1^2)y_2-y_1.
\end{aligned}
\label{vdpsys}
\end{equation}
As initial conditions, we consider $y_1(0)=2$ and $y_2(0)=0$. It can be shown that when $\mu\gg 1$, the period of the relaxation
oscillations is
$$T\approx \mu (3-2\ln{2}),$$
i.e. of the order of $\mu$ (see, e.g., \cite{strogatz2018nonlinear}). Thus, for our computations we set the time intervals to be $[0,3\mu]$, i.e. approximately three times the period of oscillation, with the exception of the case $\mu=1$, where we consider the interval $[0,30]$ (we note that in other studies, in order to assess the performance of numerical schemes this time interval is usually set to $[0,2\mu]$).

For the implementation of the proposed method, we defined a trial solution satisfying the initial conditions as follows:
\begin{equation}
\begin{aligned}
\Psi_{1}(t,\bw^o_1) & = \alpha^{(0)}_1+tN_1(t,\bw^o_1),\\
\Psi_{2}(t,\bw^o_2) & = \alpha^{(0)}_2+tN_2(t,\bw^o_2),
\end{aligned}
\label{trsolfs_vdp}
\end{equation}
where $\alpha^{(0)}_1=2$ and $\alpha^{(0)}_2=0$.  
We carried out numerical experiments within a wide range of values of the parameter $\mu$ that controls the stiffness. For our illustrations, we report the results for $\mu=1,10,100,1000$.

\begin{figure}[p]
    \centering
    \par \vskip -30pt
    \subfigure[]{
    \includegraphics[width=0.45 \textwidth]{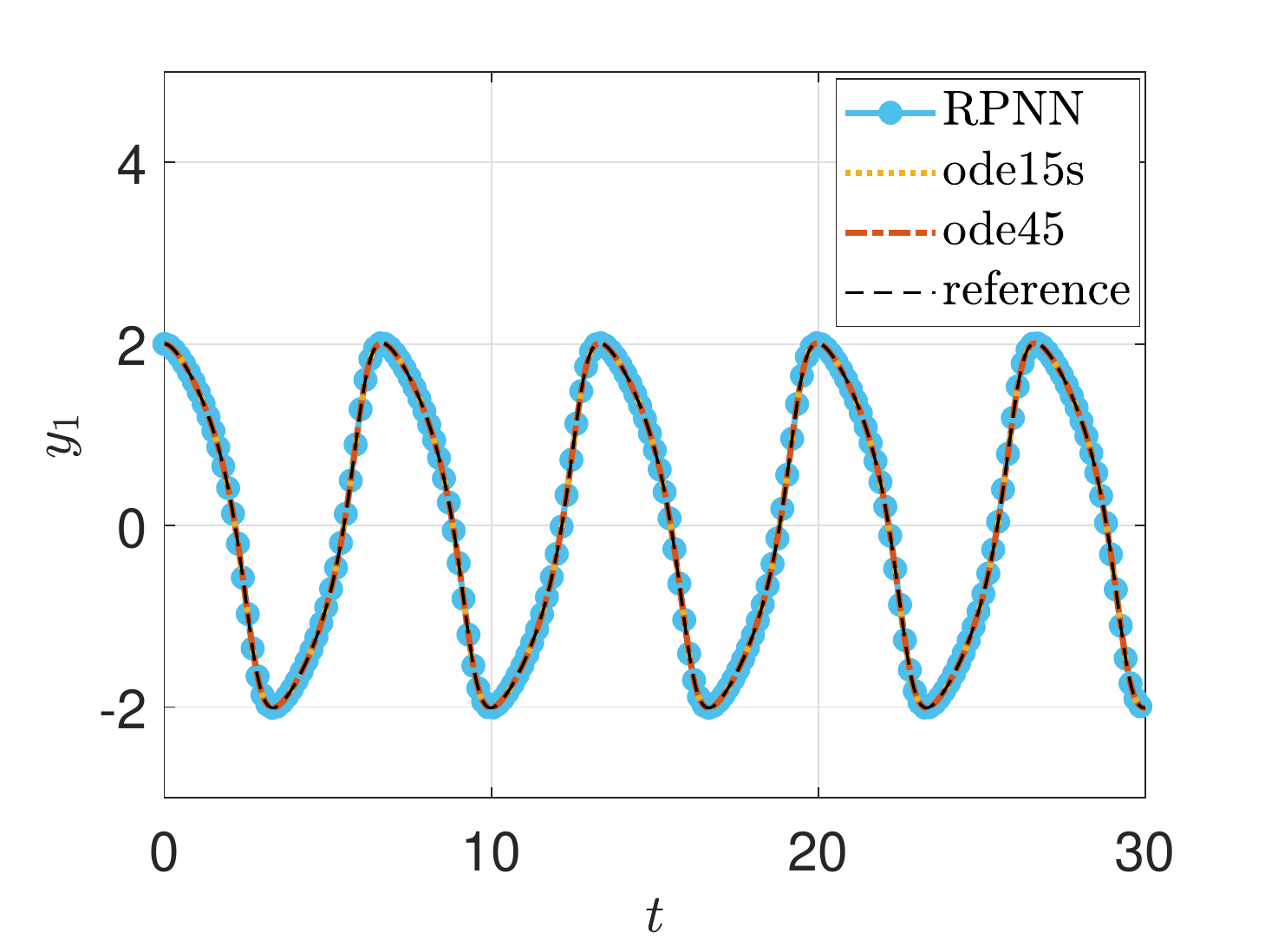}
    }
    \subfigure[]{
    \includegraphics[width=0.45 \textwidth]{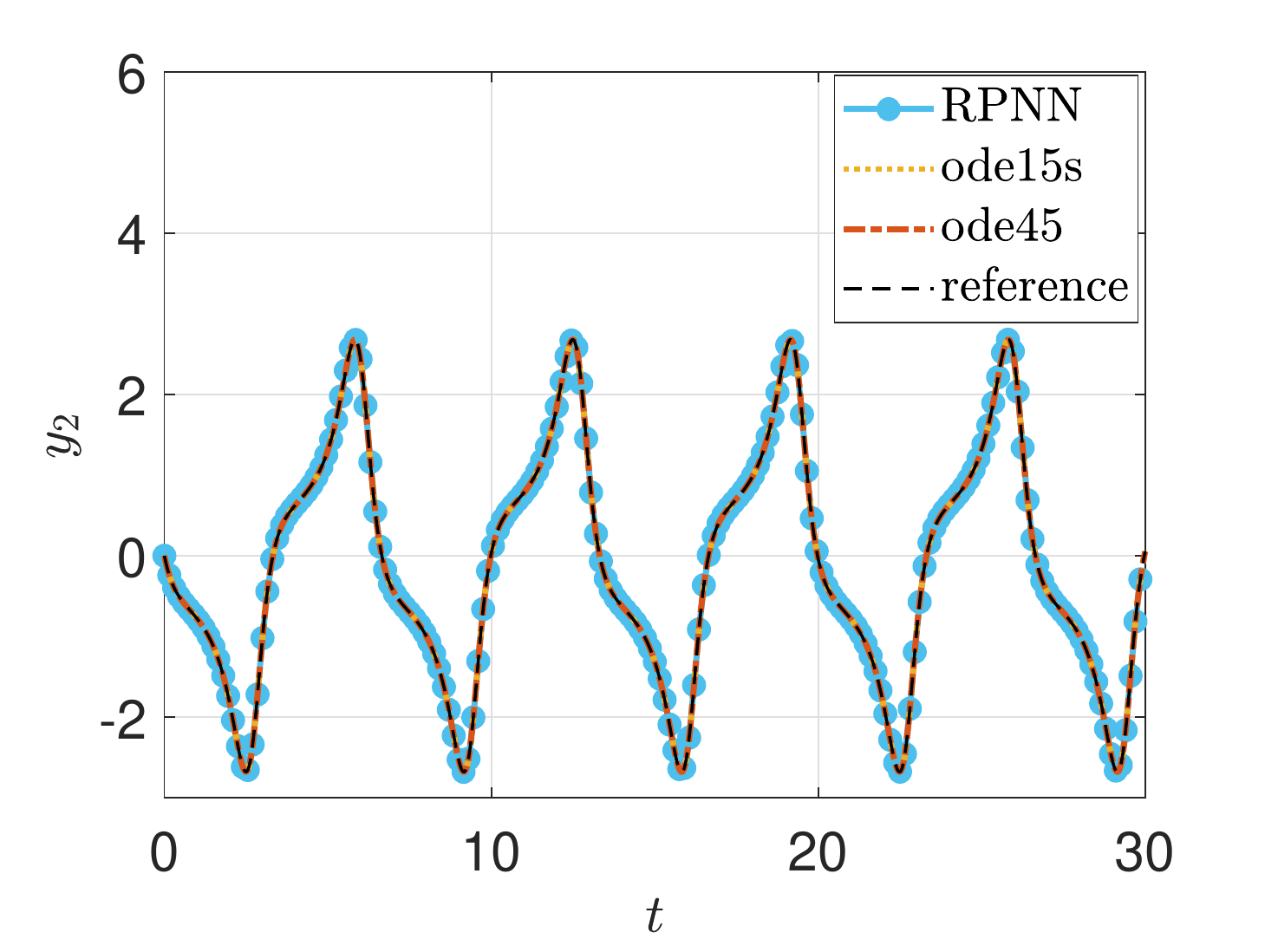}
    }
    \par \vskip -10pt
    \subfigure[]{
    \includegraphics[width=0.45 \textwidth]{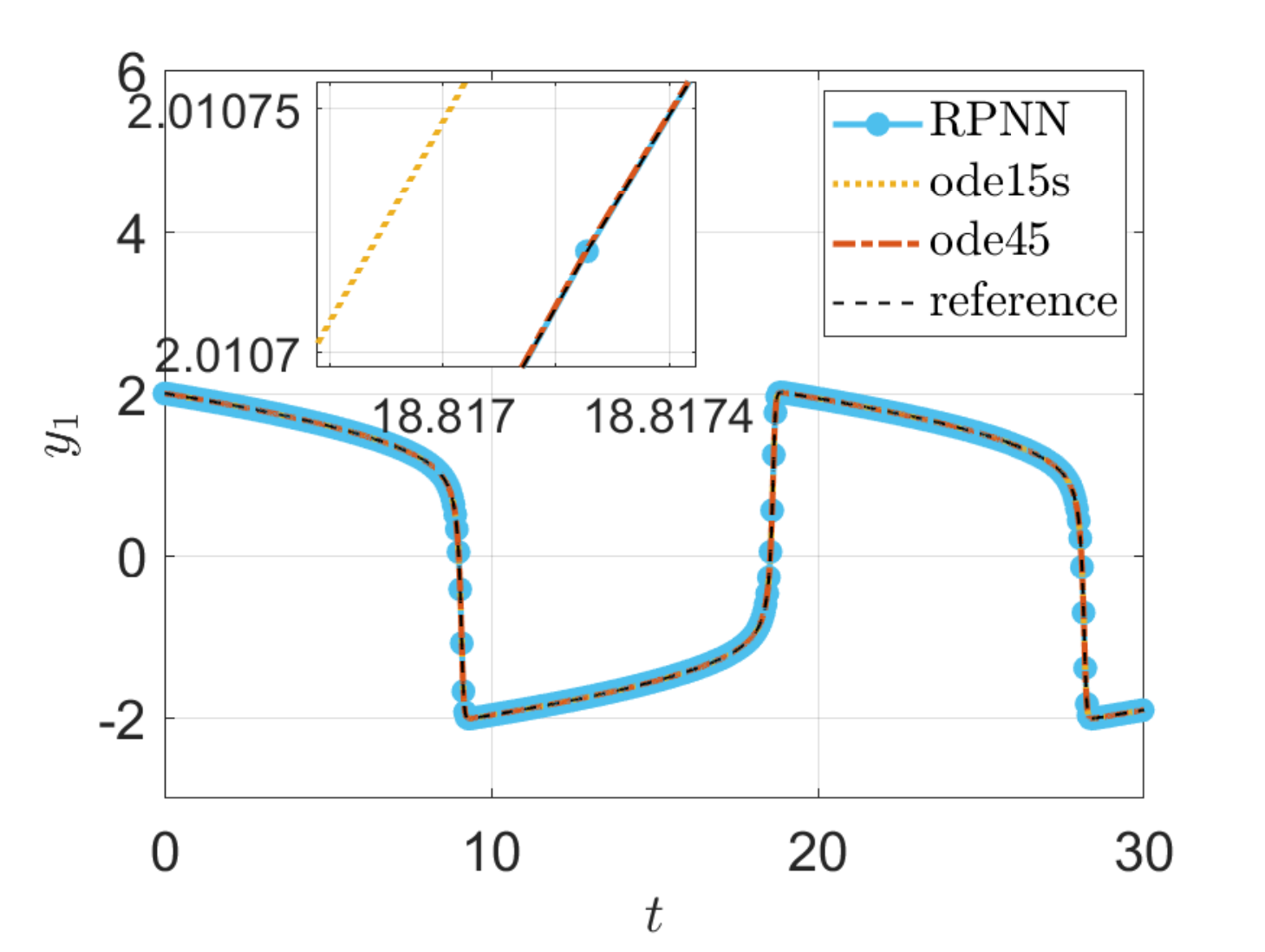}
    }
    \subfigure[]{
    \includegraphics[width=0.45 \textwidth]{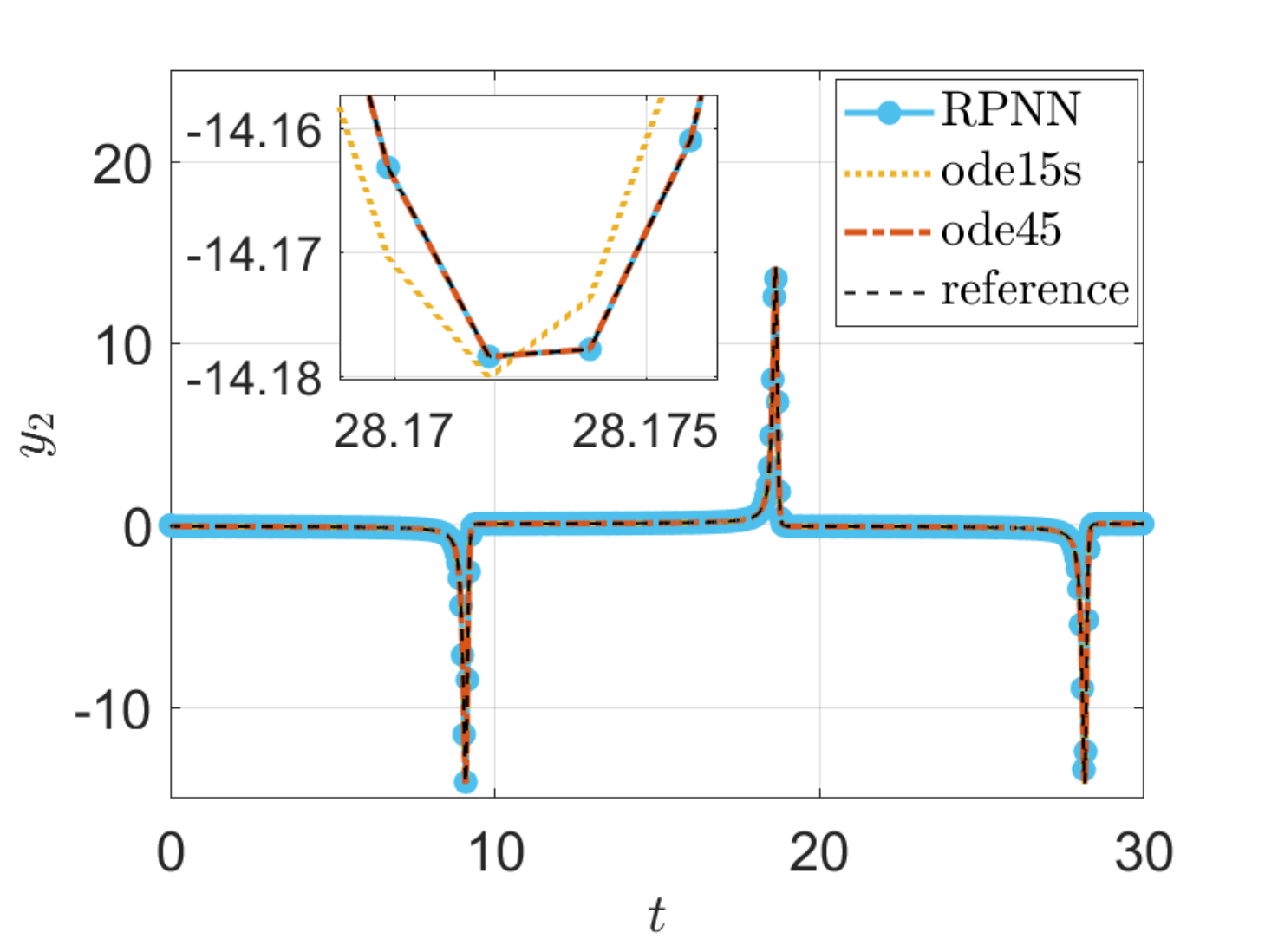}
    }
    \par \vskip -10pt
    \subfigure[]{
    \includegraphics[width=0.45 \textwidth]{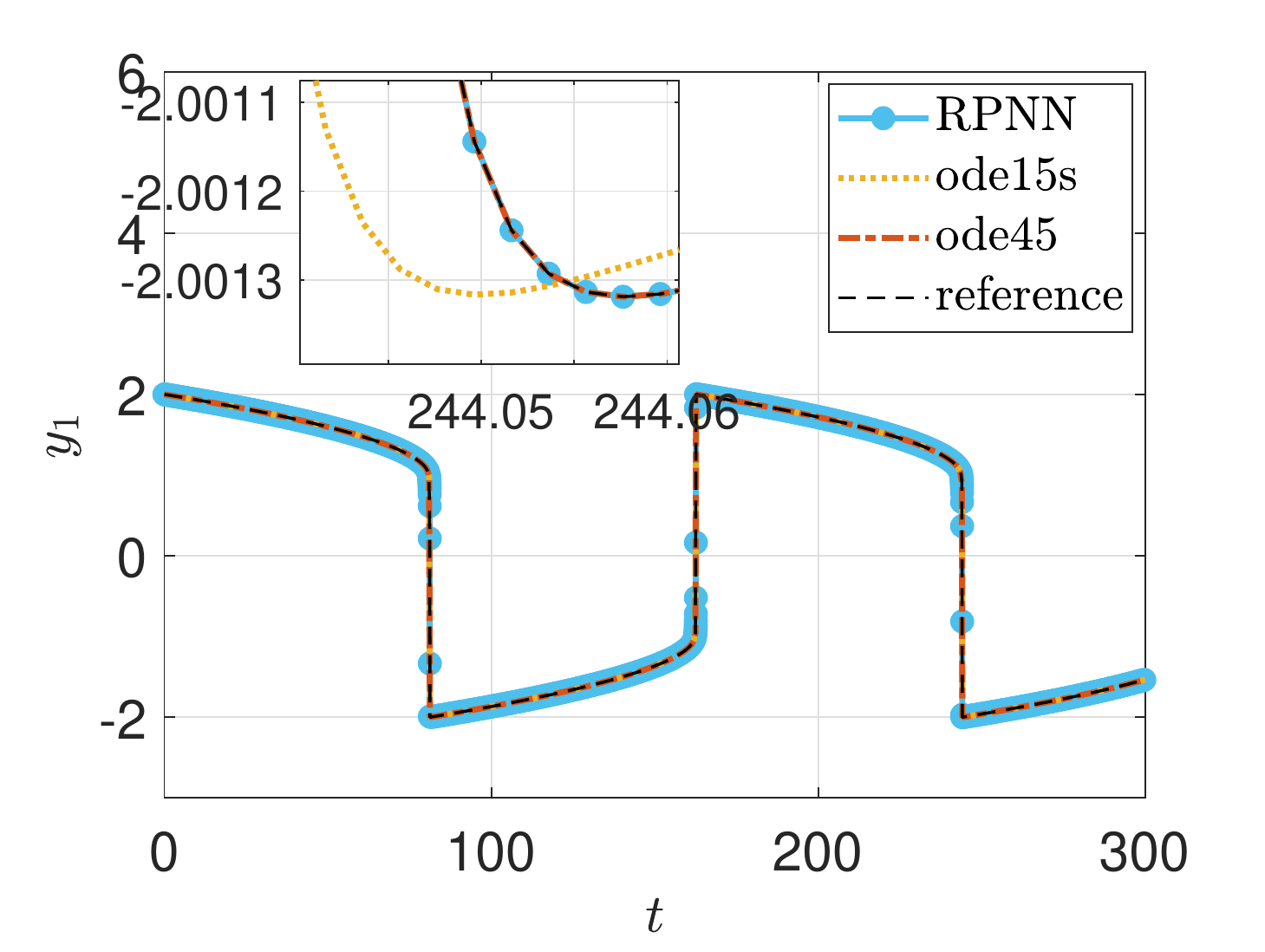}
    }
    \subfigure[]{
    \includegraphics[width=0.45 \textwidth]{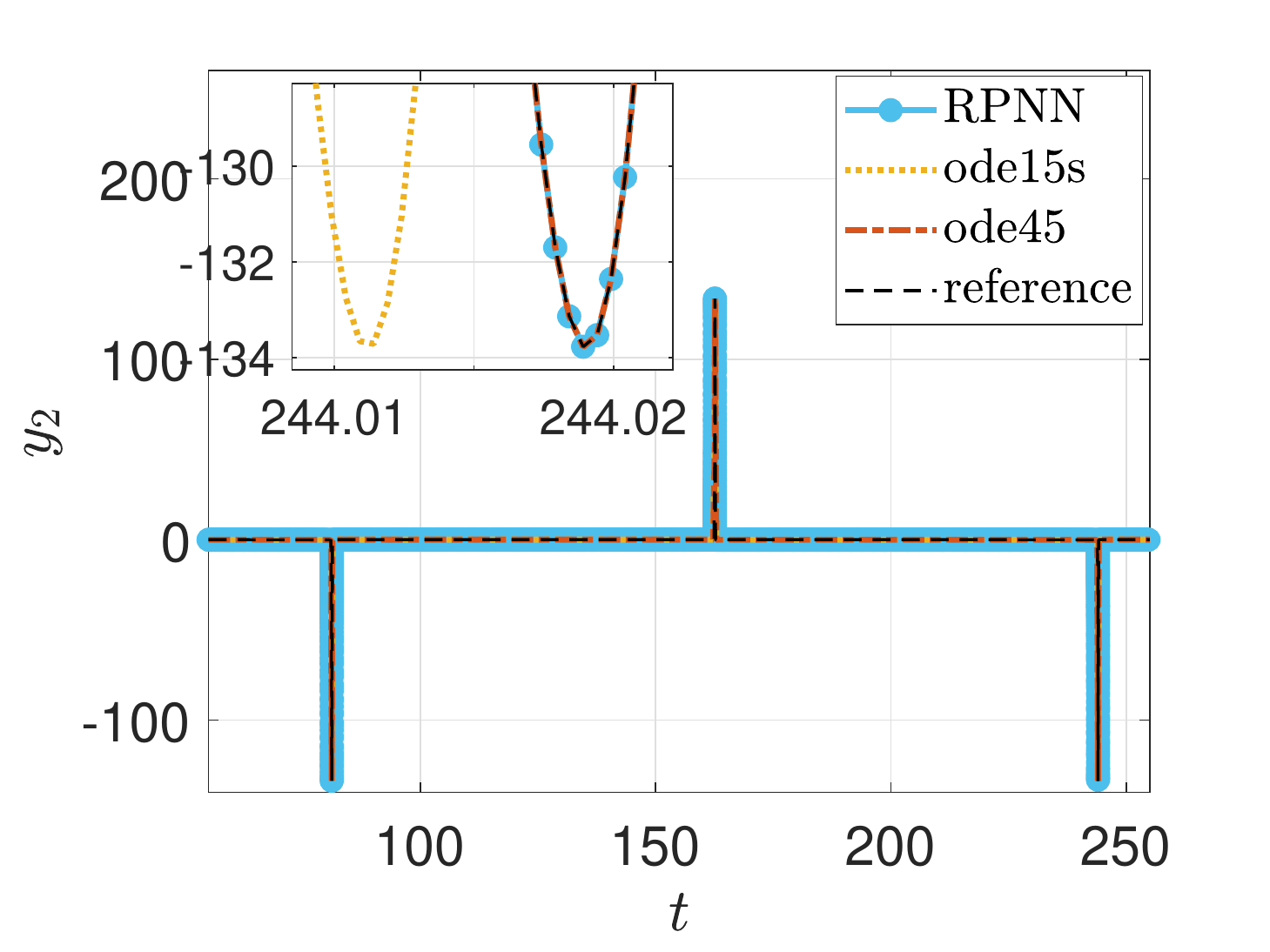}
    }
    \par \vskip -10pt
    \subfigure[]{
    \includegraphics[width=0.45 \textwidth]{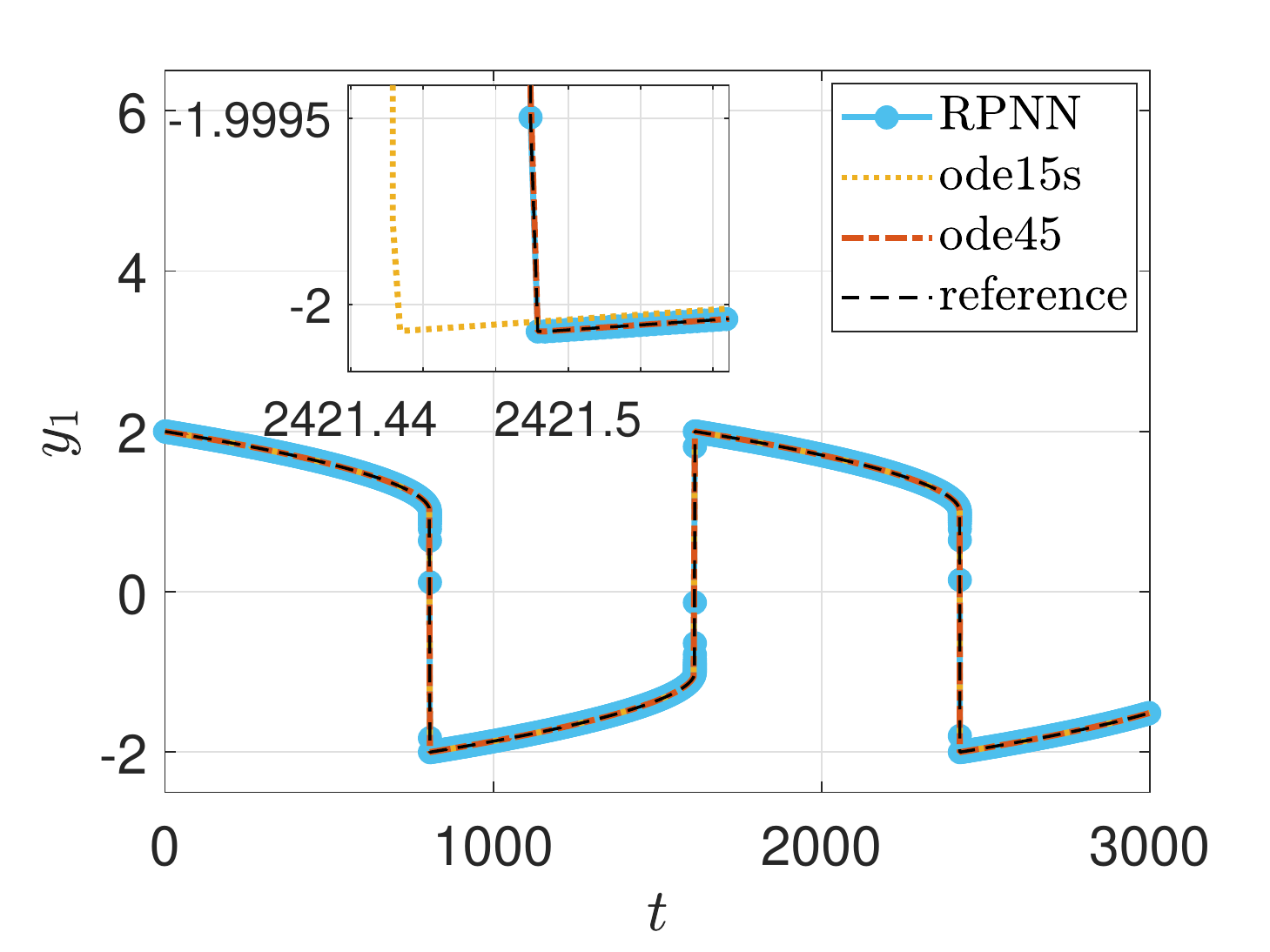}
    }
    \subfigure[]{
    \includegraphics[width=0.45 \textwidth]{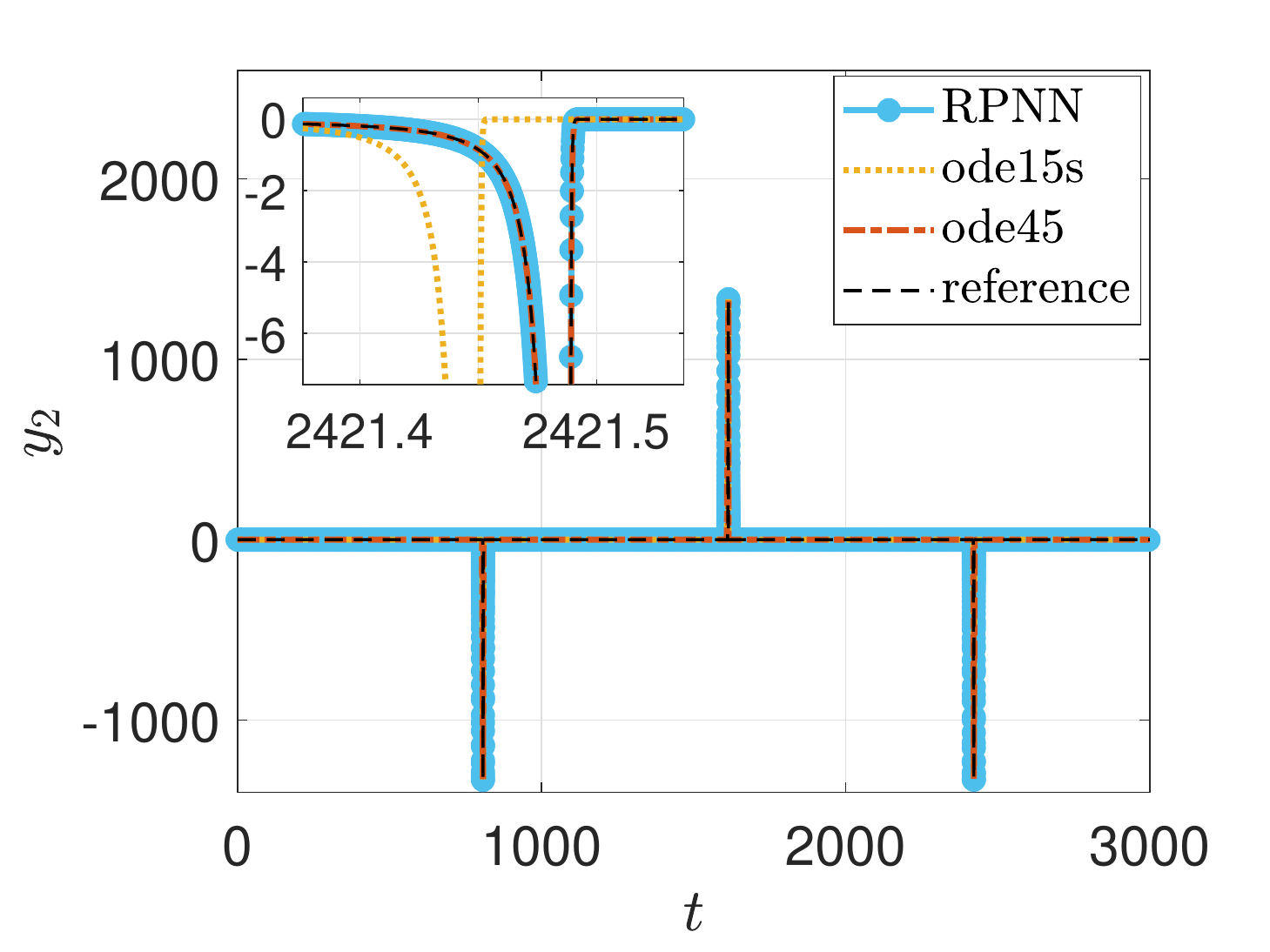}
    }
    \par \vskip -10pt
    \caption{Numerical solutions of the van der Pol problem in the interval $[0,3\mu]$ (and in $[0,30]$ for $\mu=1$) with tolerances set to $tol=1e-6$: (a)-(b) $\mu=1$, (c)-(d) $\mu=10$, (e)-(f) $\mu=100$, and (g)-(h) $\mu=1000$. The $y_1$ component is shown on the left ((a),(c),(e) and (g)) and the $y_2$ component on the right ((b),(d),(f) and (h)). The reference solution was obtained with \texttt{ode15s} with tolerances set to 1e$-$14. The insets depict a zoom around the reference solution.\label{fig:vdp_tol_6}}
\end{figure}

The numerical performance of the proposed method was validated against \texttt{ode45} and \texttt{ode15s}, using the same value for the absolute and relative tolerances. Two values for both tolerances were considered, namely 1e$-$03 and 1e$-$06. The reference solution was computed using \texttt{ode15s} with tolerances equal to 1e$-$14.
\begin{figure}[p]
    \centering
    \par \vskip -30pt
    \subfigure[]{
    \includegraphics[width=0.45 \textwidth]{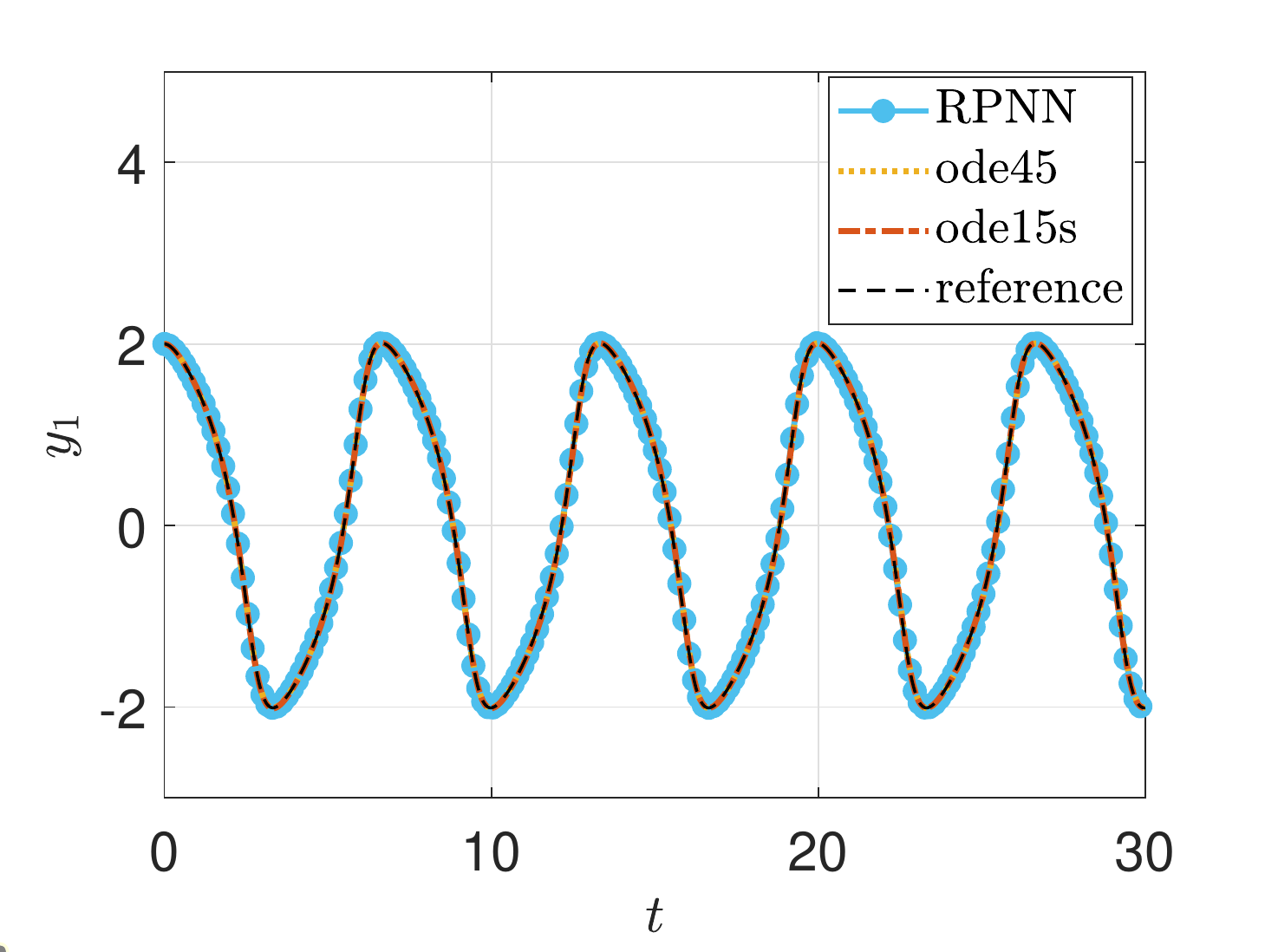}
    }
    \subfigure[]{
    \includegraphics[width=0.45 \textwidth]{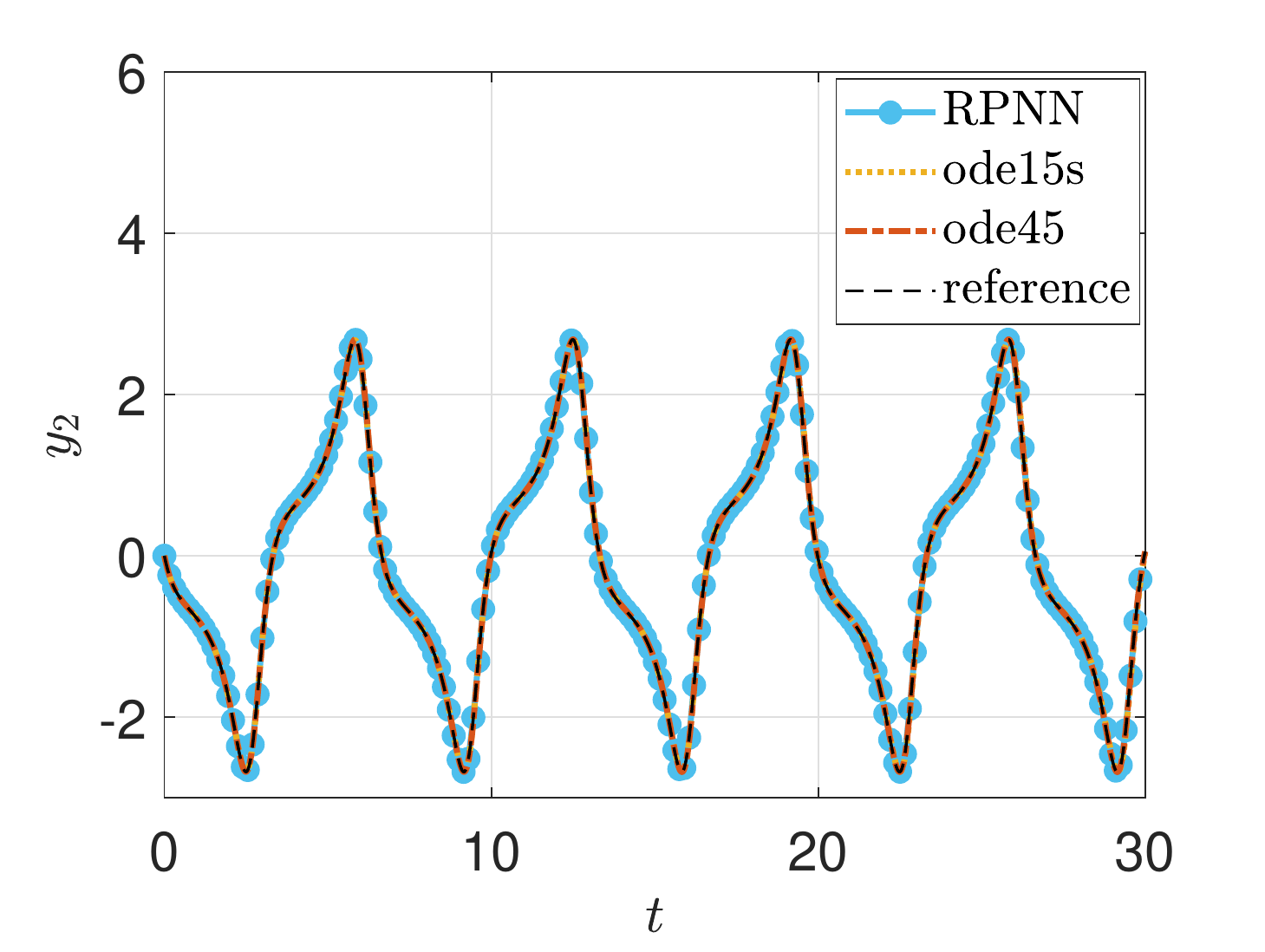}
    }
    \par \vskip -10pt
    \subfigure[]{
    \includegraphics[width=0.45 \textwidth]{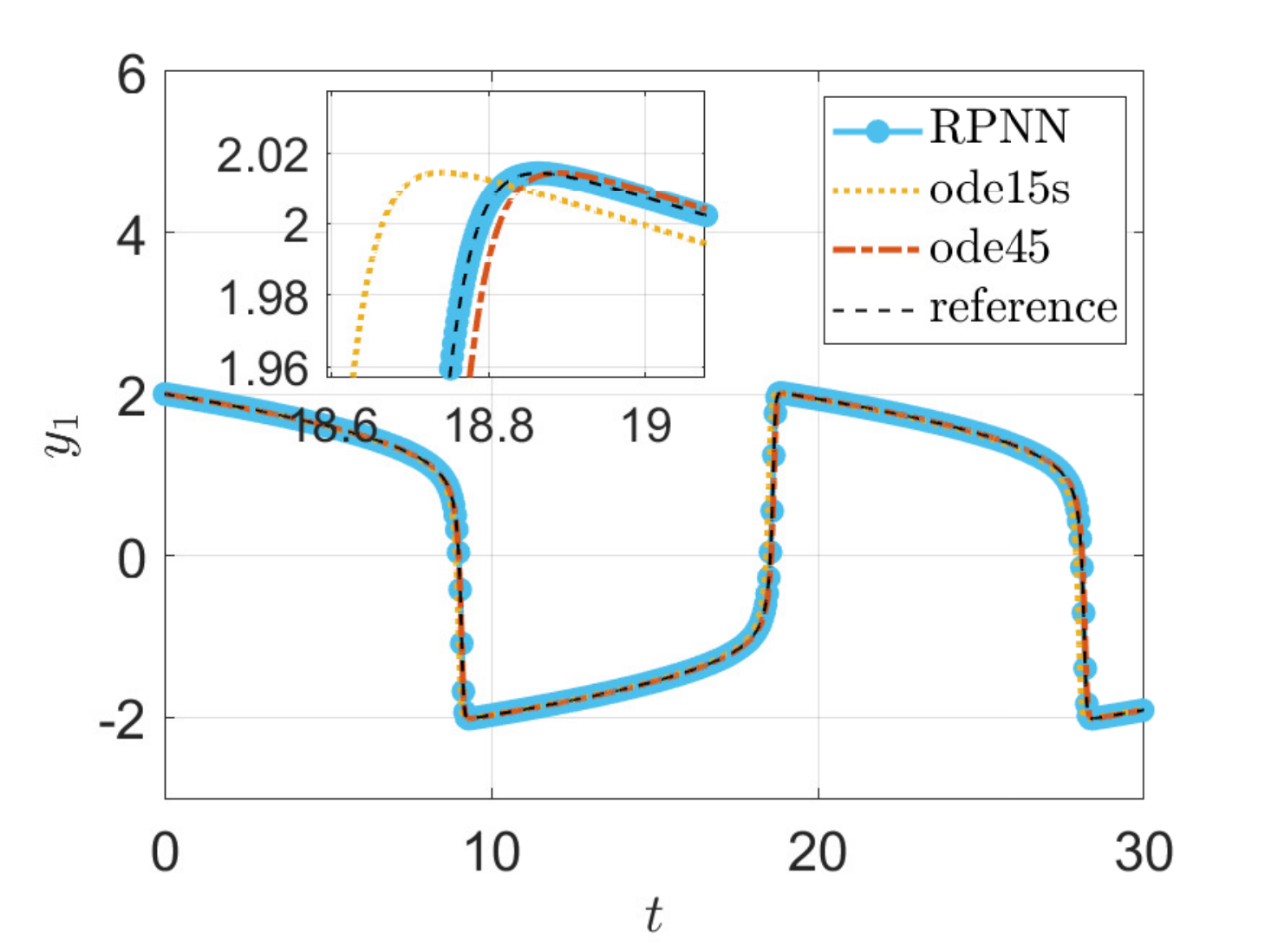}
    }
    \subfigure[]{
    \includegraphics[width=0.45 \textwidth]{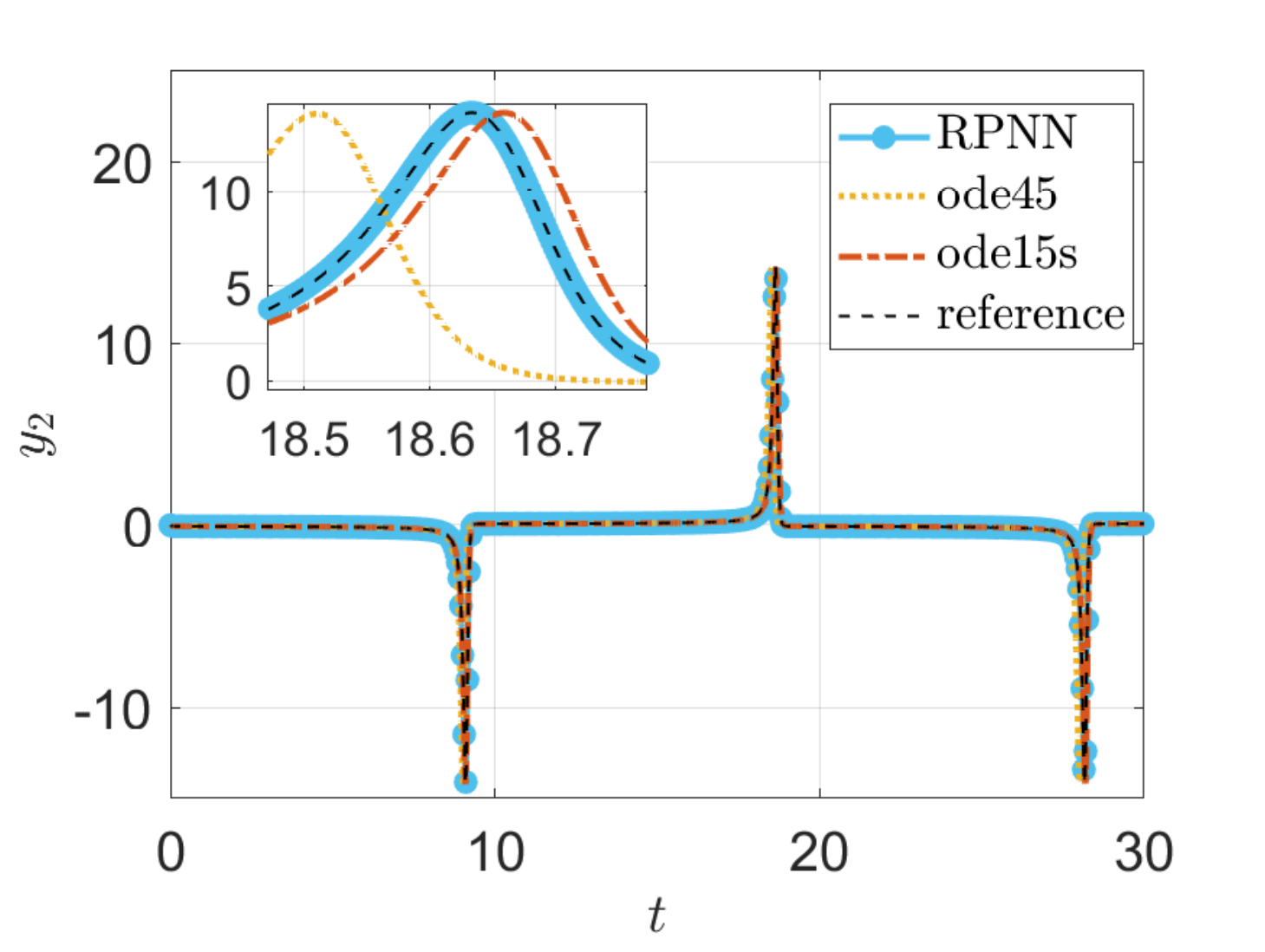}
    }
    \par \vskip -10pt
    \subfigure[]{
    \includegraphics[width=0.45 \textwidth]{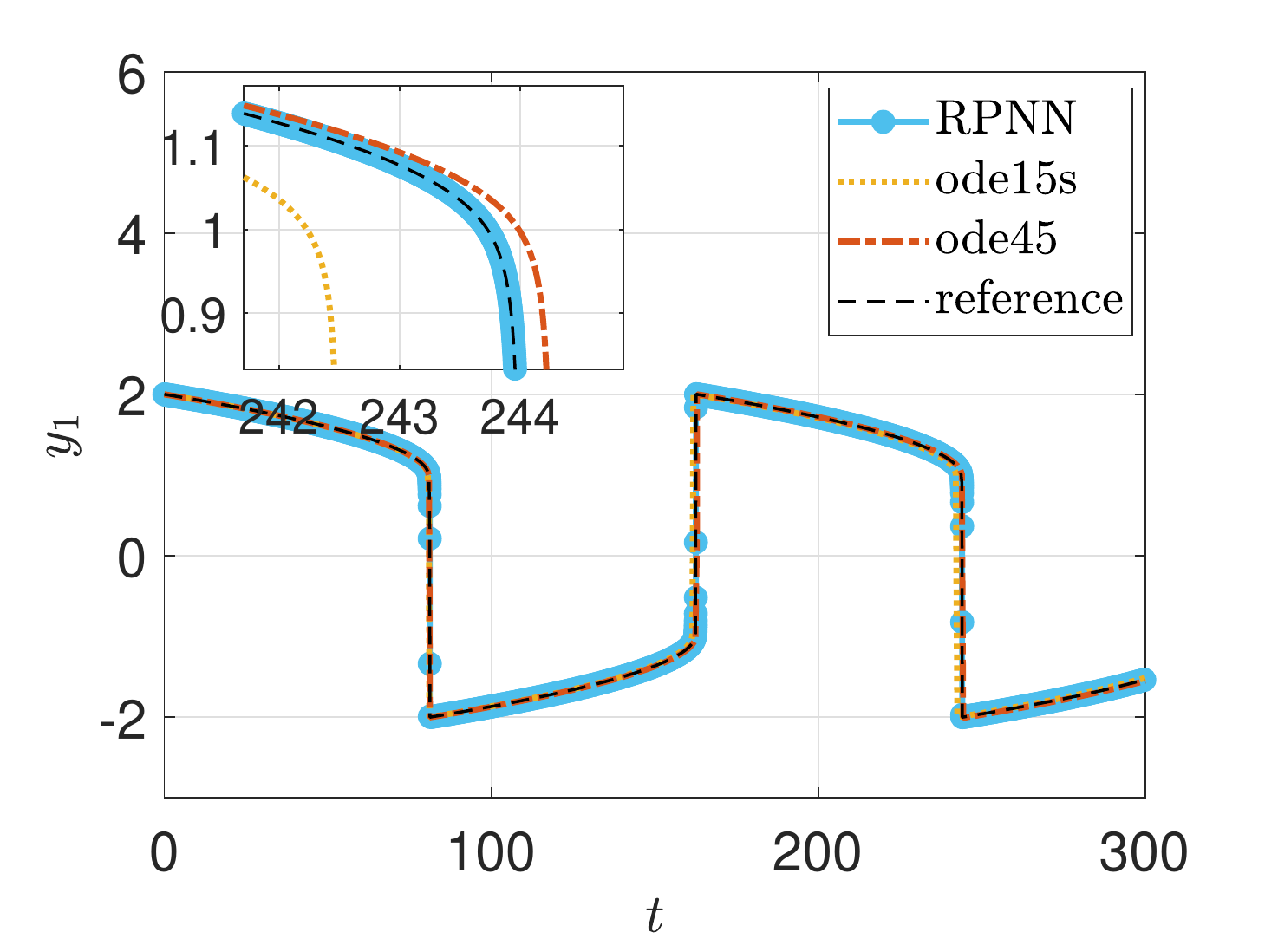}
    }
    \subfigure[]{
    \includegraphics[width=0.45 \textwidth]{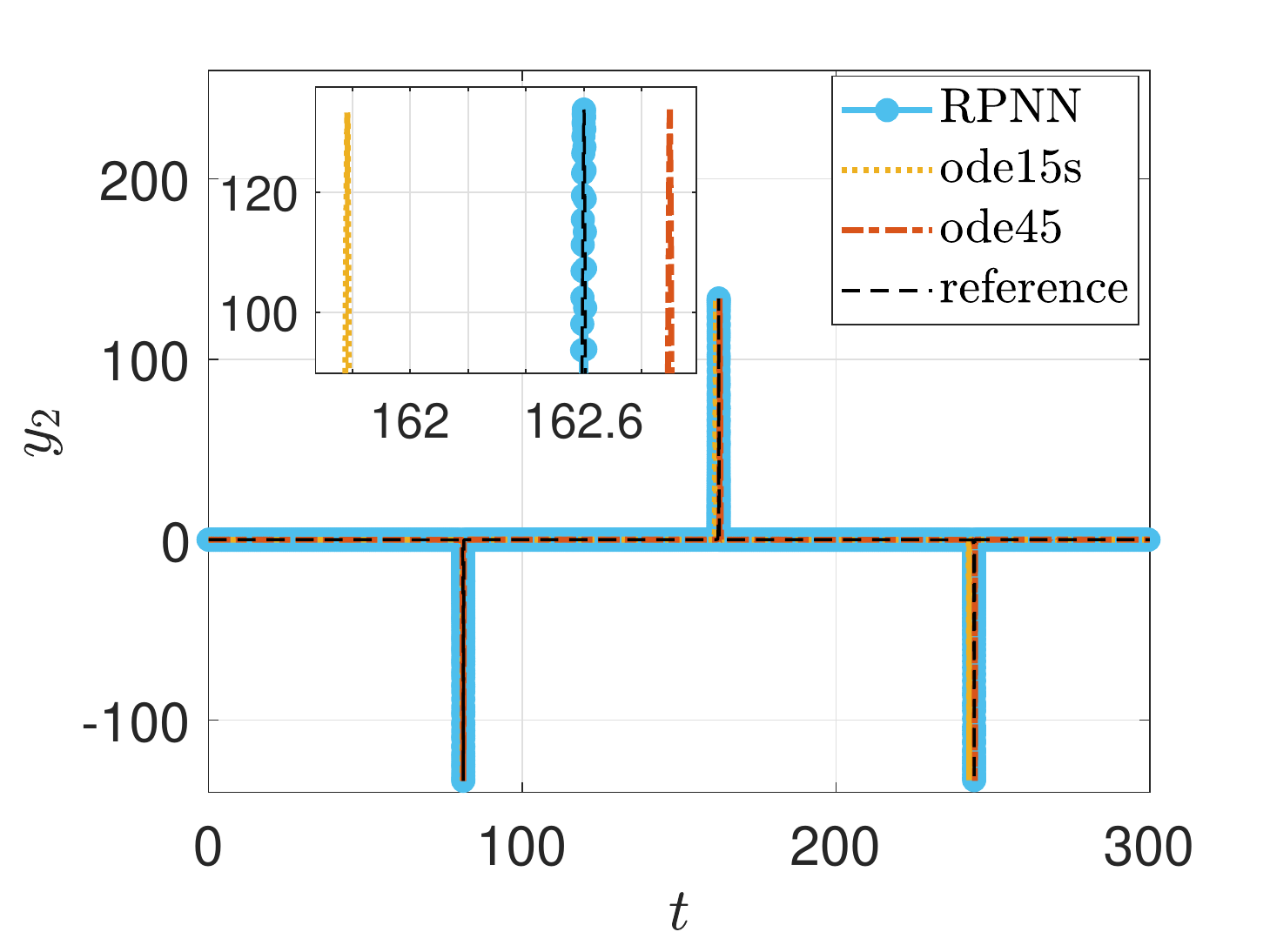}
    }
    \par \vskip -10pt
    \subfigure[]{
    \includegraphics[width=0.45 \textwidth]{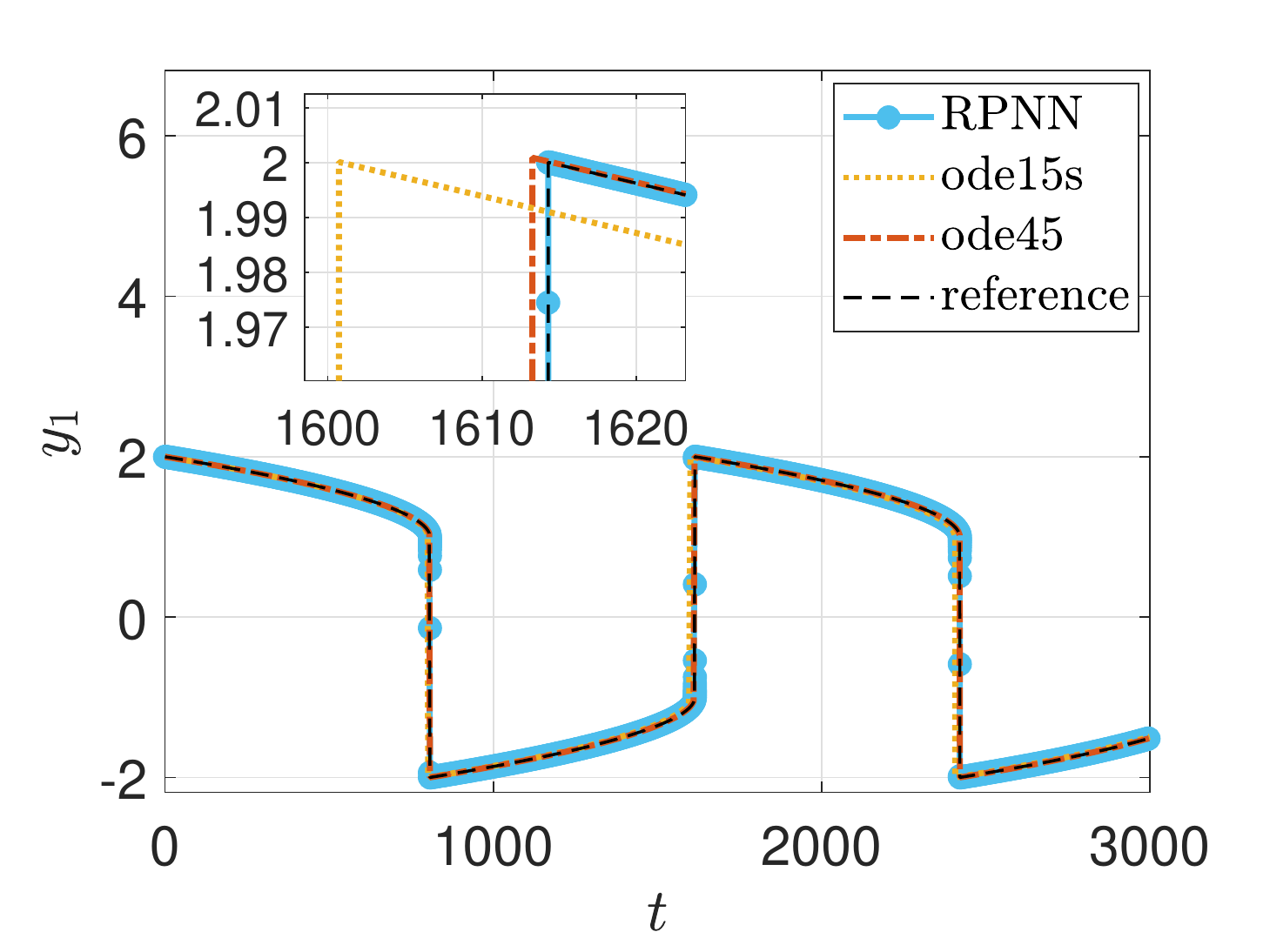}
    }
    \subfigure[]{
    \includegraphics[width=0.45 \textwidth]{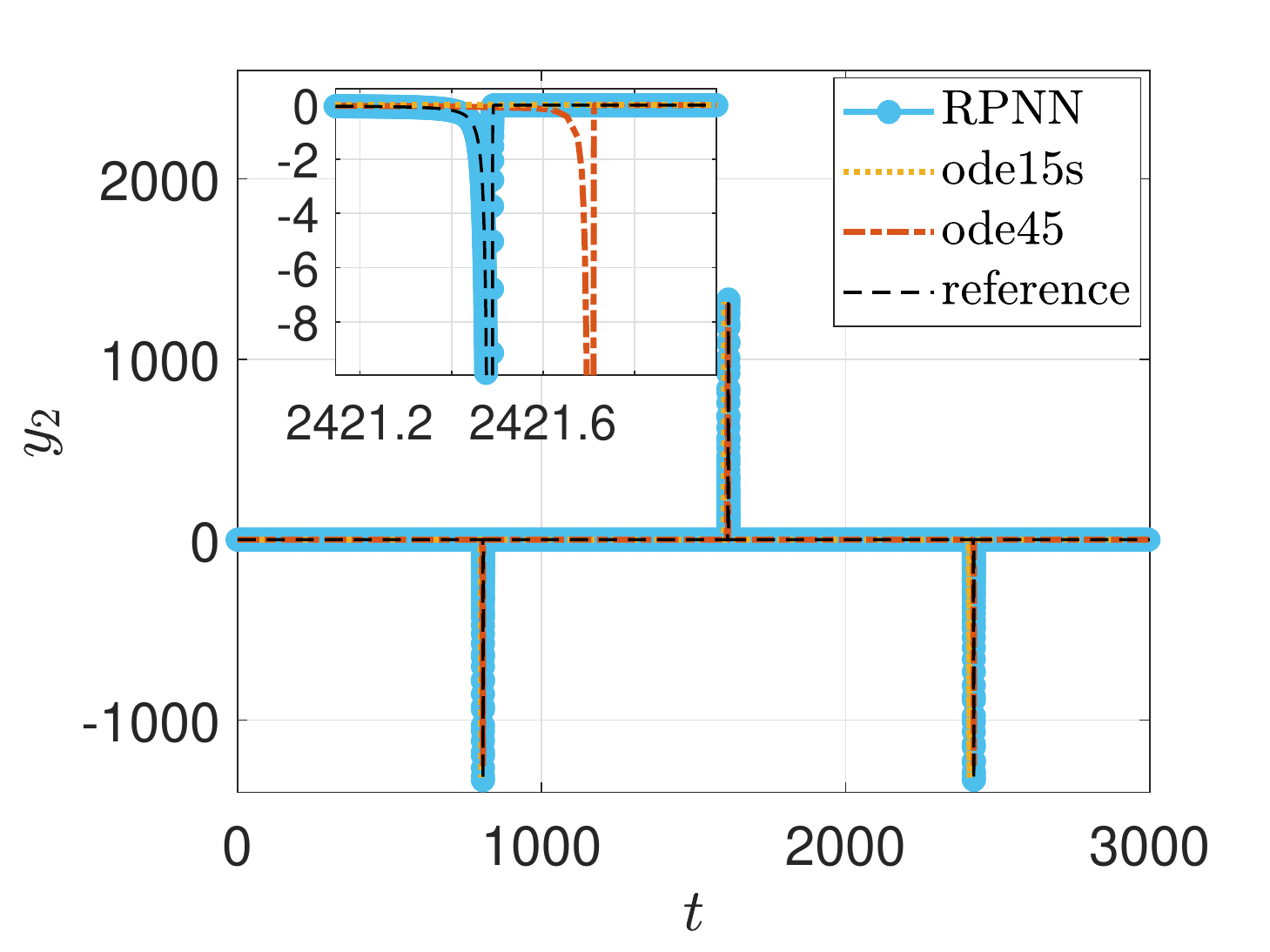}
    }
    \par \vskip -10pt
    \caption{Numerical solutions of the van der Pol problem in the interval $[0,3\mu]$ (and in $[0,30]$ for $\mu=1$), with tolerances set to $tol=1e-3$: (a)-(b) $\mu=1$, (c)-(d) $\mu=10$, (e)-(f) $\mu=100$, and (g)-(h) $\mu=1000$. The $y_1$ component is shown on the left ((a),(c),(e) and (g)) and the $y_2$ component on the right ((b),(d),(f) and (h)). The reference solution was obtained with \texttt{ode15s} with tolerances set to 1e$-$14. The insets depict a zoom around the reference solution.\label{fig:vdp_tol_3}}
\end{figure}
The resulting approximate solutions for $\mu = 1,\, 10, \, 100, \, 1000$ are depicted in Figure~\ref{fig:vdp_tol_3} for the tolerance 1e$-$03 and in Figure~\ref{fig:vdp_tol_6} for the tolerance 1e$-$06.

The corresponding approximation errors for the various values of $\mu$, in terms of $L_2$-norm, $L_{\infty}$-norm and MAE, are reported in Tables~\ref{tab:vdp_y1_accuracy}-\ref{tab:vdp_y2_accuracy}. The errors were computed using $15\,000$ equidistant grid points for $\mu=1$ and $\mu=10$, $600\,000$ points for $\mu=100$ and $30\,000\,000$ points for $\mu=1000$, as very steep gradients arise. In Table~\ref{tab:vdp_time_points}, we provide the number of points required by each method as well as the corresponding computational times (median, maximum and minimum over 10 runs), including the time required to compute the reference solution, for the various values of $\mu$. Finally, in Table ~\ref{tab:vdp_time_points_val}, we report the computational times required by each method to evaluate the solution at $10\cdot 3\mu$ equidistant points in $[0,3\mu]$ after we have trained the RPNN and solved the systems using \texttt{ode45} and \texttt{ode15s} with the selected tolerances. As we stated, the proposed method ``learns'' analytically the solution, while with the \texttt{ode45} and \texttt{ode15s} solvers the values of the solution at these points are obtained with interpolation by using the MATLAB function \texttt{deval}.

\begin{table}[h!]
\begin{center}
\caption{Van der Pol problem for $\mu=1, \, 10, \, 100, \, 1000$ in the interval $[0,3\mu]$ (and in $[0,30]$ for $\mu=1$). Absolute errors ($L_{2}$-norm, $L_{\infty}$-norm and MAE) for the component $y_1$ of the solutions computed with tolerances set to 1e$-$03 and 1e$-$06. The reference solution was obtained with \texttt{ode15s} with tolerances set to 1e$-$14.\label{tab:vdp_y1_accuracy}}
{\small
\setlength{\tabcolsep}{5pt}
\begin{tabular}{|l|l |l l l |l l l|}
\hline
\multicolumn{2}{|c|}{$y_1$} & \multicolumn{3}{c|}{$tol=$ 1e$-$03} & \multicolumn{3}{c|}{$tol=$ 1e$-$06} \\
\cline{3-8}
\multicolumn{2}{|c|}{} & $L_2$ & $L_{\infty}$ & MAE & $L_2$ & $L_{\infty}$ & MAE\\
\hline
& RPNN   & 1.25e$-$04 & 3.62e$-$06 & 7.01e$-$07 & 4.53e$-$06 & 1.10e$-$07 & 2.71e$-$08\\
$\mu=1$ & \texttt{ode45}  & 1.18e$+$00 &2.96e$-$02 & 7.10e$-$03 & 4.54e$-$04 & 1.11e$-$05 & 2.75e$-$06\\
 & \texttt{ode15s} & 1.82e$+$00 & 4.88e$-$02 & 1.01e$-$02 & 2.29e$-$03 & 5.87e$-$05 & 1.32e$-$05 \\
\hline
& RPNN   & 7.31e$-$04 & 7.09e$-$05 & 1.33e$-$06 & 1.24e$-$05 & 1.27e$-$06 & 2.26e$-$08\\
$\mu=10$ & \texttt{ode45}  & 4.06e$+$00 & 3.61e$-$01 & 7.06e$-$03 & 7.70e$-$04 & 8.73e$-$05 & 1.30e$-$06\\
 & \texttt{ode15s} & 2.30e$+$01 & 1.79e$+$00 & 4.73e$-$02 & 9.55e$-$02 & 1.02e$-$02 & 1.74e$-$04 \\
 \hline
& RPNN & 6.44e$-$02 & 9.95e$-$03 & 2.03e$-$06 & 4.97e$-$04 & 7.38e$-$05 & 1.63e$-$08\\
$\mu=100$ & \texttt{ode45}  & 4.77e$+$01 & 3.00e$+$00 & 9.69e$-$03 & 2.43e$-$03 & 3.56e$-$04 & 8.12e$-$08\\
 & \texttt{ode15s} & 2.12e$+$02 & 3.10e$+$00 & 3.45e$-$02 & 6.40e$+$00 & 9.83e$-$01 & 2.12e$-$04 \\
\hline
 & RPNN   & 4.77e$+$00 & 9.75e$-$01 & 2.13e$-$06 & 7.44e$-$02 & 1.90e$-$02 & 2.77e$-$06\\
$\mu=1000$ & \texttt{ode45}  & 3.35e$+$02 & 3.03e$+$00 & 1.79e$-$03 & 2.34e$-$01 & 4.89e$-$02 & 3.34e$-$07\\
 & \texttt{ode15s} & 1.77e$+$03 & 3.12e$+$00 & 4.64e$-$02 & 6.87e$+$01 & 2.97e$+$00 & 1.71e$-$04 \\
\hline
\end{tabular}
}
\end{center}
\end{table}
\begin{table}[h!]
\begin{center}
\caption{Van der Pol problem for $\mu=1, \, 10, \, 100, \, 1000$ in the interval $[0,3\mu]$ (and in $[0,30]$ for $\mu=1$). Absolute errors ($L_{2}$-norm, $L_{\infty}$-norm and MAE) for the component $y_2$ of the solutions computed with tolerances set to 1e$-$03 and 1e$-$06. The reference solution was obtained with \texttt{ode15s} with tolerances set to 1e$-$14.\label{tab:vdp_y2_accuracy}}
{\small
\setlength{\tabcolsep}{5pt}
\begin{tabular}{|l|l |l l l |l l l|}
\hline
\multicolumn{2}{|c|}{$y_2$} & \multicolumn{3}{c|}{$tol=$ 1e$-$03} & \multicolumn{3}{c|}{$tol=$ 1e$-$06} \\
\cline{3-8}
\multicolumn{2}{|c|}{} & $L_2$ & $L_{\infty}$ & MAE & $L_2$ & $L_{\infty}$ & MAE\\
\hline
& RPNN   & 1.81e$-$04 & 6.41e$-$06 & 9.82e$-$07 & 6.55e$-$06 & 1.99e$-$07 & 3.72e$-$08\\
$\mu=1$ & \texttt{ode45}  & 1.68e$+$00 & 5.71e$-$02 & 8.87e$-$03 & 6.55e$-$04 & 2.02e$-$05 & 3.73e$-$06\\
 & \texttt{ode15s} & 2.71e$+$00 & 8.43e$-$02 & 1.48e$-$02 & 3.49e$-$03 & 9.57e$-$05 & 2.07e$-$05 \\
\hline
& RPNN   & 7.08e$-$03 & 7.81e$-$04 & 9.77e$-$06 & 1.20e$-$04 & 1.40e$-$05 & 1.67e$-$07\\
$\mu=10$ & \texttt{ode45}  & 3.89e$+$01 & 3.95e$+$00 & 4.85e$-$02 & 7.45e$-$03 & 9.63e$-$04 & 8.91e$-$06\\
 & \texttt{ode15s} & 1.79e$+$02 & 1.37e$+$01 & 2.93e$-$01 & 9.25e$-$01 & 1.12e$-$01 & 1.27e$-$03 \\
 \hline
 & RPNN   & 3.22e$+$00 & 1.10e$+$00 & 1.46e$-$04 & 4.72e$-$02 & 7.84e$-$03 & 1.09e$-$06\\
$\mu=100$ & \texttt{ode45}  & 6.72e$+$02 & 1.34e$+$02 & 4.05e$-$02 & 2.31e$-$01 & 3.78e$-$02 & 4.88e$-$06\\
 & \texttt{ode15s} & 8.23e$+$02 & 1.33e$+$02 & 6.08e$-$02 & 5.71e$+$02 & 9.33e$+$01 & 1.32e$-$02 \\
 \hline
 & RPNN   & 4.20e$+$03 & 9.05e$+$02 & 1.40e$-$03 & 7.05e$+$01 & 2.01e$+$01 & 1.92e$-$05\\
$\mu=1000$ & \texttt{ode45}  & 1.13e$+$04 & 1.33e$+$03 & 6.08e$-$03 &  2.26e$+$02 & 5.21e$+$01 & 6.54e$-$05\\
 & \texttt{ode15s} & 1.16e$+$04 & 1.33e$+$03 & 6.22e$-$03 & 1.05e$+$04 & 1.33e$+$03 & 4.94e$-$03\\
 \hline
\end{tabular}
}
\end{center}
\end{table}
\begin{table}[ht]
\begin{center}
\caption{Van der Pol problem for $\mu=1, \, 10, \, 100, \, 1000$ in the interval $[0,3\mu]$ (and in $[0,30]$ for $\mu=1$). Computational times in seconds (median, minimum and maximum times over 10 runs) and number of points required with tolerances set to 1e$-$03 and 1e$-$06. The reference solution was obtained with \texttt{ode15s} with tolerances set to 1e$-$14.\label{tab:vdp_time_points}}
{\small
\setlength{\tabcolsep}{3pt}
\begin{tabular}{|l |l |l l l r |l l l r|}
\hline
\multicolumn{2}{|c|}{} & \multicolumn{4}{c|}{$tol=$ 1e$-$03} & \multicolumn{4}{c|}{$tol=$ 1e$-$06} \\
\cline{3-10}
\multicolumn{2}{|c|}{} & median & min & max & \multicolumn{1}{l|}{\# pts} & median & min & max & \multicolumn{1}{l|}{\# pts} \\
\hline
& RPNN & 1.25e$-$01 & 1.10e$-$01 & 1.33e$-$01 &  740  & 1.89e$-$01 & 1.74e$-$01 & 2.79e$-$01 & 1045\\
$\mu=1$ & \texttt{ode45}  & 1.98e$-$03 & 1.70e$-$03 & 2.24e$-$03 & 77 & 6.82e$-$03 & 5.05e$-$03 & 3.97e$-$02 & 263\\
& \texttt{ode15s} & 6.99e$-$03 & 6.52e$-$03 & 7.39e$-$03 & 267 & 1.71e$-$02 & 1.63e$-$02 & 2.99e$-$02 & 862\\
& reference & 2.98e$-$01 & 2.92e$-$01 & 3.22e$-$01 &  15474  & 2.98e$-$01 & 2.92e$-$01 & 3.22e$-$01 &  15474\\
\hline
& RPNN & 1.39e$-$01 & 1.32e$-$01 & 1.67e$-$01 &  779  & 2.07e$-$01 & 1.90e$-$01 & 2.22e$-$01 & 1199\\
$\mu=10$ & \texttt{ode45}  & 4.21e$-$03 & 3.92e$-$03 & 1.71e$-$02 & 199 & 7.81e$-$03 & 7.52e$-$03 & 1.34e$-$02 & 439\\
& \texttt{ode15s} & 9.15e$-$03 & 8.43e$-$03 & 2.30e$-$02 & 258 & 1.76e$-$02 & 1.68e$-$02 & 2.89e$-$02 & 782\\
& reference & 2.67e$-$01 & 2.51e$-$01 & 4.16e$-$01 & 12924  & 2.67e$-$01 & 2.51e$-$01 & 4.16e$-$01 & 12924\\
\hline
& RPNN & 2.57e$-$01 & 2.42e$-$01 & 2.87e$-$01 &  1691  & 6.93e$-$01 & 6.02e$-$01 & 7.68e$-$01 & 3183\\
$\mu=100$ & \texttt{ode45}  & 2.78e$-$01 & 2.71e$-$01 & 3.00e$-$01 & 16959 & 2.91e$-$01 & 2.83e$-$01 & 3.83e$-$01  & 17142\\
& \texttt{ode15s} & 1.18e$-$02 & 1.14e$-$02 & 1.44e$-$02 & 347 & 2.79e$-$02 & 2.70e$-$02 & 4.43e$-$02  & 1119\\
& reference & 3.96e$-$01 & 3.90e$-$01 & 5.00e$-$01 &  19868  & 3.96e$-$01 & 3.90e$-$01 & 5.00e$-$01 &  19868\\
\hline
& RPNN & 6.24e$-$01 & 5.89e$-$01 & 6.47e$-$01 &  3962  & 8.02e$+$00 & 7.85e$+$00 & 8.68e$+$00 & 15671\\
$\mu=1000$ & \texttt{ode45}  & 2.78e$+$01 & 2.74e$+$01  &  2.84e$+$01 & 1684374 & 2.76e$+$01 & 2.71e$+$01 & 2.79e$+$01 & 1685043\\
& \texttt{ode15s} & 1.57e$-$02 & 1.50e$-$02 & 2.68e$-$02 & 461 & 3.46e$-$02 & 3.32e$-$02 & 4.17e$-$02  & 1401\\
& reference & 5.26e$-$01 &    5.13e$-$01 & 5.68e$-$01 &  25569  & 5.37e$-$01 & 5.24e$-$01 & 5.46e$-$01 &  25569\\
\hline
\end{tabular}
}
\end{center}
\end{table}
\begin{table}[ht]
\begin{center}
\caption{Van der Pol problem for $\mu=1, \, 10, \, 100, \, 1000$ in the interval $[0,3\mu]$ (and in $[0,30]$ for $\mu=1$). Computational times in seconds (median, minimum and maximum times over 10 runs) for the evaluation (interpolation) of the solution in $10\cdot 3\mu$ equidistant grid points after the solution was obtained by the three schemes with tolerances set to 1e$-$03 and 1e$-$06. The reference solution was obtained with \texttt{ode15s} with tolerances set to 1e$-$14.\label{tab:vdp_time_points_val}}
{\small
\setlength{\tabcolsep}{5pt}
\begin{tabular}{|l |l |l l l |l l l|}
\hline
\multicolumn{2}{|c|}{} & \multicolumn{3}{c|}{$tol=$ 1e$-$03} & \multicolumn{3}{c|}{$tol=$ 1e$-$06} \\
\cline{3-8}
\multicolumn{2}{|c|}{} & median & min & max & median & min & max \\
\hline
& RPNN & 8.94e$-$04 & 8.45e$-$04 & 1.39e$-$03  & 1.33e$-$03 &   1.17e$-$03 & 4.15e$-$03\\
$\mu=1$ & \texttt{ode45}  & 1.47e$-$03 & 1.40e$-$03 & 1.73e$-$03 & 3.96e$-$03 & 3.75e$-$03 & 7.62e$-$03\\
& \texttt{ode15s} & 4.08e$-$03 & 3.95e$-$03 & 7.86e$-$03 & 5.12e$-$03 & 4.93e$-$03 & 5.60e$-$03\\
& reference & 1.15e$-$02 & 1.08e$-$02 & 1.27e$-$02 & 1.23e$-$02 & 1.10e$-$02 & 1.71e$-$02\\
\hline
& RPNN & 1.06e$-$03 & 9.82e$-$04 & 1.43e$-$03 & 1.62e$-$03 & 1.47e$-$03 & 9.18e$-$03\\
$\mu=10$ & \texttt{ode45}  & 3.03e$-$03  & 2.94e$-$03 &     3.42e$-$03 & 4.38e$-$03 & 4.26e$-$03 & 4.66e$-$03\\
& \texttt{ode15s} & 2.07e$-$03 & 1.98e$-$03 & 3.54e$-$03 & 3.67e$-$03 & 3.55e$-$03 & 3.75e$-$03\\
& reference & 1.12e$-$02 & 1.06e$-$02 & 1.33e$-$02 & 1.05e$-$02 & 9.74e$-$03 & 1.15e$-$02\\
\hline
& RPNN & 6.50e$-$03 & 6.12e$-$03  & 6.81e$-$03 & 1.03e$-$02 & 8.69e$-$03 & 1.43e$-$02\\
$\mu=100$ & \texttt{ode45}  & 1.23e$-$01 & 1.20e$-$01 & 1.38e$-$01 & 1.18e$-$01 & 1.13e$-$01 & 1.31e$-$01\\
& \texttt{ode15s} & 3.08e$-$03 & 2.95e$-$03 & 3.39e$-$03 & 7.11e$-$03 & 7.03e$-$03 & 7.65e$-$03\\
& reference & 9.72e$-$02 & 9.34e$-$02 & 1.16e$-$01 & 9.64e$-$02 & 9.29e$-$02 & 1.04e$-$01\\
\hline
& RPNN & 6.09e$-$02 & 5.73e$-$02 & 6.35e$-$02 & 2.52e$-$01 & 2.24e$-$01 & 2.62e$-$01\\
$\mu=1000$ & \texttt{ode45}  & 6.34e$+$01 & 6.13e$+$01 & 6.61e$+$01 & 5.91e$+$01 & 5.88e$+$01 &  6.80e$+$01\\
& \texttt{ode15s} & 1.49e$-$02 & 1.42e$-$02 & 1.77e$-$02 & 3.38e$-$02  & 3.26e$-$02 & 3.43e$-$02\\
& reference & 4.67e$-$01 & 4.34e$-$01 & 4.79e$-$01 & 4.64e$-$01 & 4.56e$-$01 & 4.75e$-$01\\
\hline
\end{tabular}
}
\end{center}
\end{table}

As shown in Tables~\ref{tab:vdp_y1_accuracy}-\ref{tab:vdp_y2_accuracy}, the proposed numerical method outperforms both the \texttt{ode45} and \texttt{ode15s} solvers in terms of numerical accuracy for all values of $\mu$ in all metrics and for both tolerances. It is worthy to note that for the larger values of $\mu$, i.e. $\mu=100, 1000$, where the solution exhibits steep gradients in the specific time interval, \texttt{ode15s} gives relatively large errors. For $\mu=100$, the $L_{\infty}$-norm of the difference between the values of $y_1$ ($y_2$) computed by \texttt{ode15s} and the corresponding components of the reference solution in the interval where a steep gradient arises is $\|y_1 -y^{ref}_1\|_{L_{\infty}} \simeq 3$ ($\|y_2-y^{ref}_2 \|_{L_{\infty}} \simeq 133$) when the tolerances are set to 1e$-$03, and $\|y_1 -y^{ref}_1\|_{L_{\infty}} \simeq 1$ ($\|y_2-y^{ref}_2 \|_{L_{\infty}} \simeq 10$) when the tolerances are set to 1e$-$06. The maximum absolute reference value is $y^{ref}_1 \simeq 2$ ($y^{ref}_2 \simeq 133$). That is, \texttt{ode15s} fails to adequately approximate the solution. This happens because the approximate solutions obtained with \texttt{ode15s} do not ``catch'' adequately the time interval where the very steep gradient arises, even when the tolerances are set to 1e$-$06 (see Figure~\ref{fig:vdp_tol_3}(e)(f) and Figure~\ref{fig:vdp_tol_6}(e)(f)). For the implementation of the proposed machine learning scheme, the corresponding value of the $L_{\infty}$-norm of the approximation error is $\|y_1 -y^{ref}_1\|_{L_{\infty}} \simeq 0.009$ ($\|y_2 -y^{ref}_2\|_{L_{\infty}} \simeq 1$) when the tolerance is set to 1e$-$03, and $\|y_1 -y^{ref}_1\|_{L_{\infty}} \simeq 1$e$-$04 ($\| y_2 -y^{ref}_2\|_{L_{\infty}} \simeq 2$e$-$07) when the tolerance is set to 1e$-$06. The proposed scheme results in solutions that follow closely also the very steep gradients that arise for even higher stiffness, as for example those arising for $\mu=1000$, while \texttt{ode15s} with the same tolerance fails to do that (see Figure~\ref{fig:vdp_tol_3}(g)(h) and Figure~\ref{fig:vdp_tol_6}(g)(h)).

A similar behaviour is observed when the value of the stiffness parameter is small: both \texttt{ode45} and \texttt{ode15s} result in considerably larger approximation errors when compared with the proposed machine learning scheme. For example, for $\mu=1$, when the tolerances are set to 1e$-$03, the $L_{\infty}$-norm of the difference between the value of $y_1$ ($y_2$) provided by the schemes and the reference solution are as follows: with \texttt{ode45} $\|y_1-y^{ref}_1 \|_{L_{\infty}} \simeq 3$e$-$02 ($\|y_2-y^{ref}_2 \|_{L_{\infty}} \simeq 6$e$-$02), with \texttt{ode15s} $\| y_1-y^{ref}_1 \|_{L_{\infty}} \simeq 5$e$-$02 ($\| y_2-y^{ref}_2 \|_{L_{\infty}} \simeq 8$e$-$02), and finally with the proposed scheme $\|y_1 -y^{ref}_1\|_{L_{\infty}} \simeq 4$e$-$06 ($\|y_2 -y^{ref}_2\|_{L_{\infty}} \simeq 6$e$-$06). When the tolerances are set equal to 1e$-$06, we have a behaviour similar to the case when they are set to 1e$-$03: the proposed method provides better approximations when compared with both \texttt{ode45} and \texttt{ode15s} (see Tables~\ref{tab:vdp_y1_accuracy}-\ref{tab:vdp_y2_accuracy}).

Regarding the execution times, those required by the RPNN method are generally larger than those of the classical methods (see Table~\ref{tab:vdp_time_points}), and the faster solver is \texttt{ode15s}, which however fails to approximate adequately the solution (see Tables~\ref{tab:vdp_y1_accuracy}-\ref{tab:vdp_y2_accuracy})). We note that for the largest value of $\mu$, i.e. $\mu=1000$, the implementation of our algorithm outperforms significantly the \texttt{ode45} solver also in terms of computational times, because the number of points that are required by \texttt{ode45} is huge (of the order of 1.7 million).
Furthermore, we underline that the computational times required by the two classical methods to compute the solutions via interpolation at $10 \cdot 3\mu$ equidistant points in the interval $[0, 3\mu]$ are larger than the ones required by the proposed method (see Table~\ref{tab:vdp_time_points_val}), since the proposed machine learning scheme upon training, and in contrast to the other two classical methods, offers analytical functions for the solutions, thus their computation in any point of the interval is straightforward.

Therefore, we conclude that our RPNN-based approach performs well in terms of numerical accuracy regardless of the level of stiffness, thus outperforming the other two methods, while the required computational times are at least comparable with those of \texttt{ode15s} and considerably smaller 
when the stiffness is large.

\subsection{Case Study 3: ROBER problem}

The ROBER model was developed to describe the kinetics of an autocatalytic reaction~\cite{Robertson1966}. The set of the reactions reads:
\begin{equation}
\begin{aligned}
A\quad & \rightarrow^{k_1}\quad C+B, \\
B+B\quad & \rightarrow^{k_2}\quad C+B,\\
B+C\quad & \rightarrow^{k_3}\quad A+C,
\end{aligned}
\label{strrober}
\end{equation}
where $A$, $B$, $C$ are chemical species and $k_1$, $k_2$ and $k_3$ are reaction rate constants. Assuming idealized conditions and the mass action law is applied for the rate functions, we have the following system of ODEs:
\begin{equation}
\begin{aligned}
y_1'& = -k_1y_1+k_2y_2y_3,\\
y_2' & = k_1y_1--k_2y_2y_3-k_3y_2^2,\\
y_3' & =k_3y_2^2,
\end{aligned}
\label{modrober}
\end{equation}
where $y_1$, $y_2$ and $y_3$ denote the concentrations of $A$, $B$ and $C$, respectively. In our simulations, we set the typical values of the parameters, i.e. $k_1=0.04$, $k_2=10^4$ and $k_3=3\times10^7 $, with $y_1(0)=1$, $y_2(0)=0$ and $y_3(0)=0$ as initial values of the concentrations, and we consider the same time interval as in the original paper \cite{Robertson1966}, i.e. $[0, 40]$. Stiffness arises from the large differences among the reaction rate constants.

\begin{figure}[p]
    \centering
    \par \vskip -10pt
    \subfigure[]{
    \includegraphics[width=0.45 \textwidth]{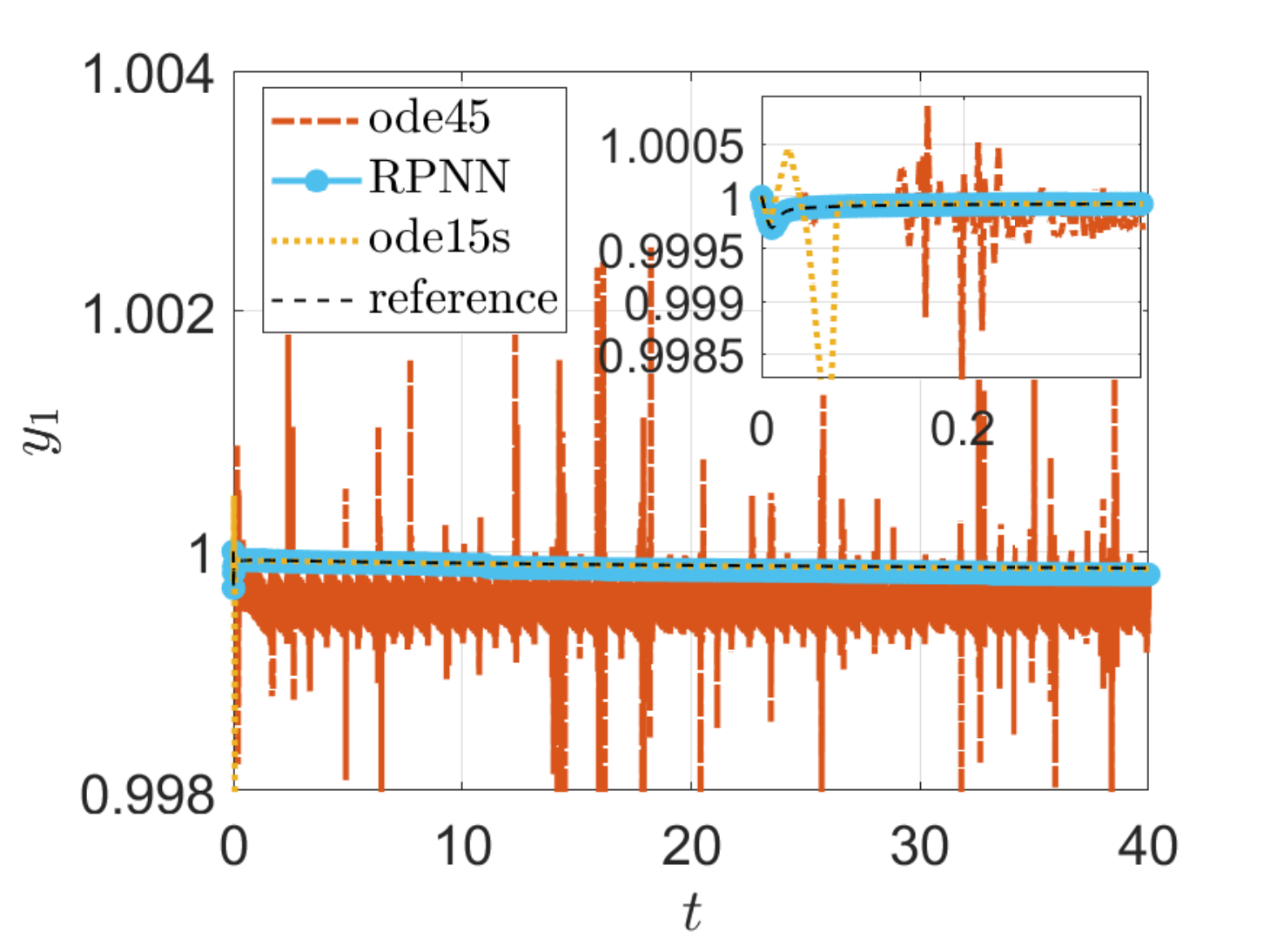}
    }
    \subfigure[]{
    \includegraphics[width=0.45 \textwidth]{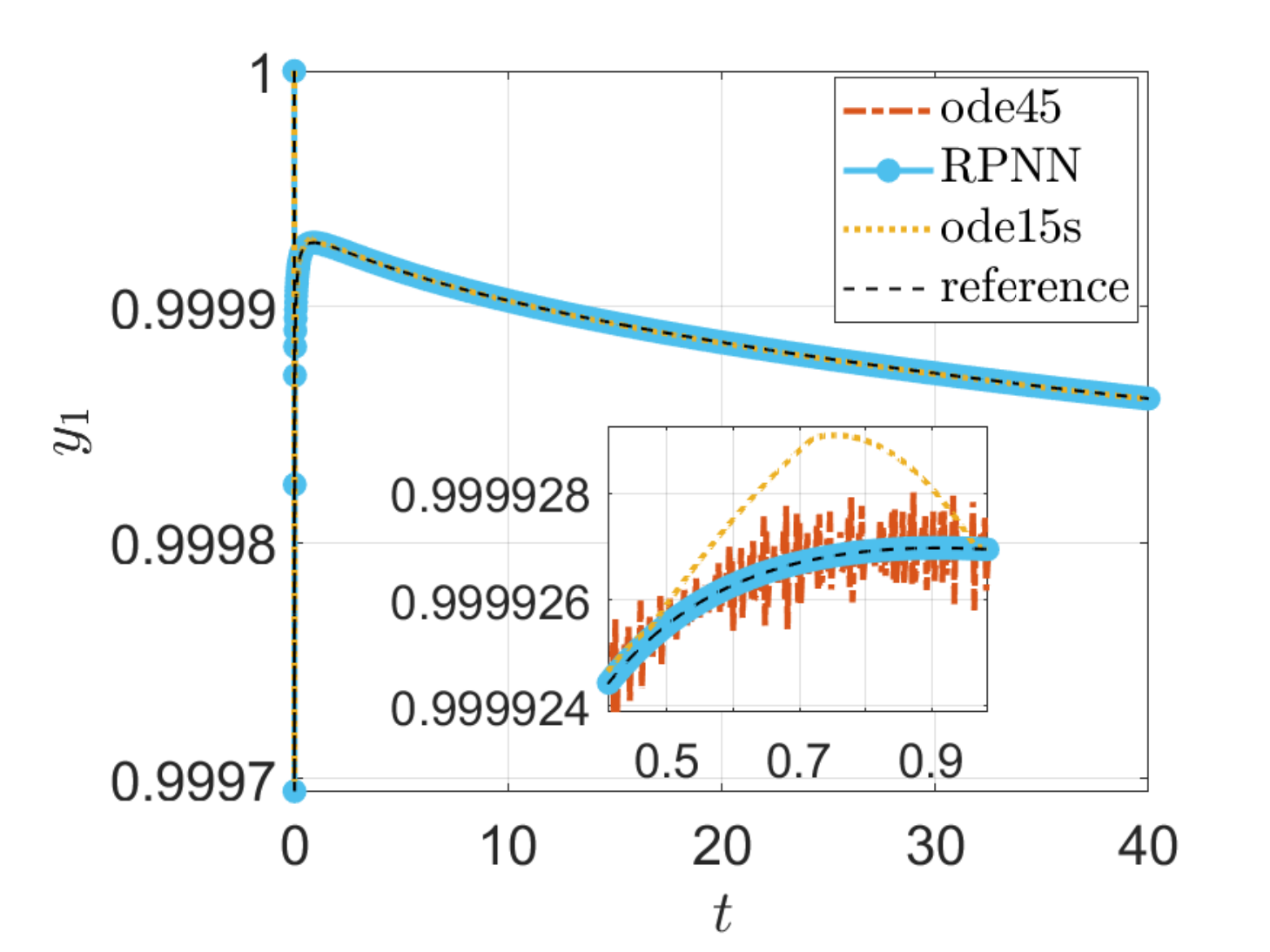}
    }
    par \vskip -10pt
    \subfigure[]{
    \includegraphics[width=0.45 \textwidth]{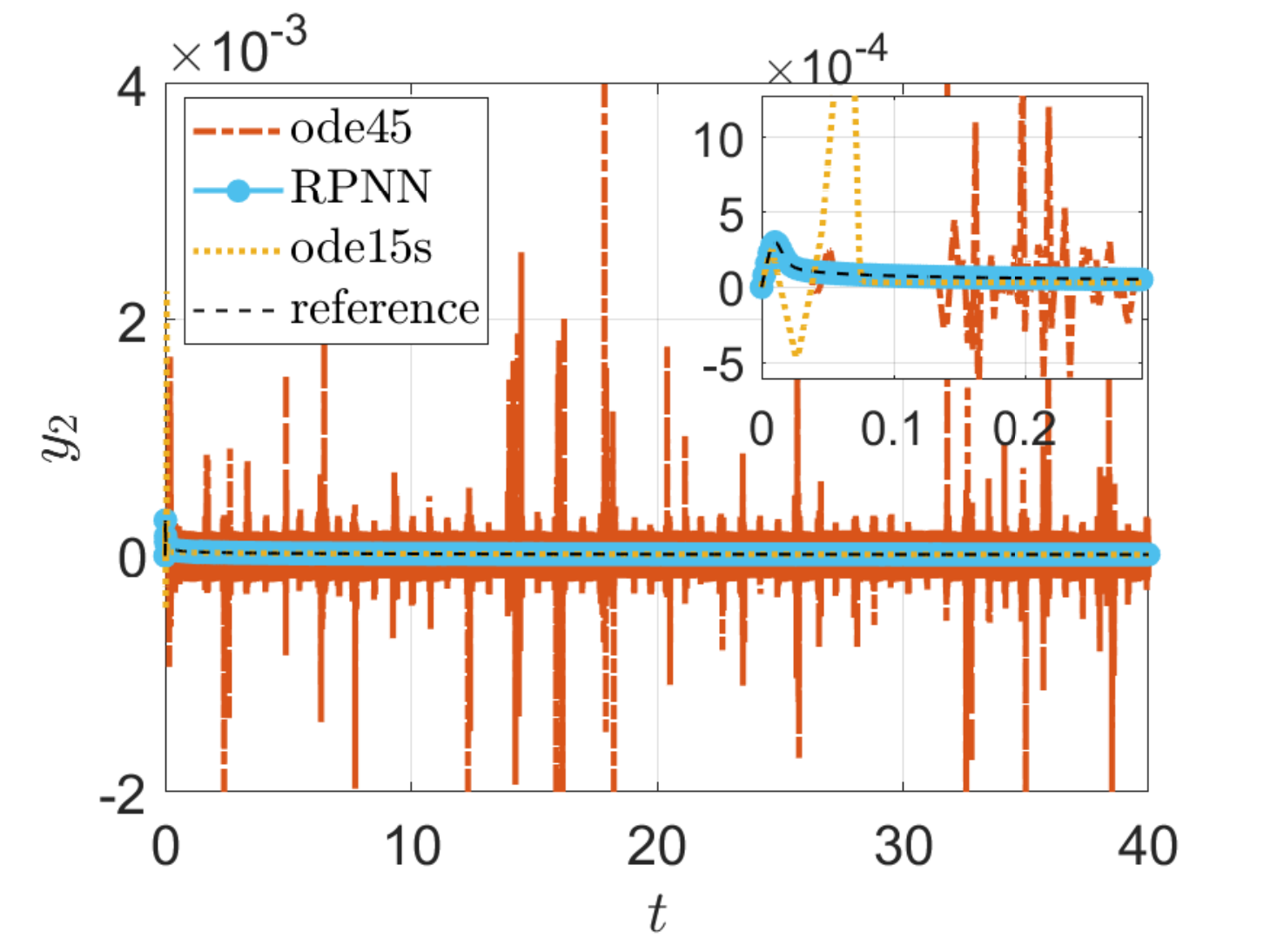}
    }
    \subfigure[]{
    \includegraphics[width=0.45 \textwidth]{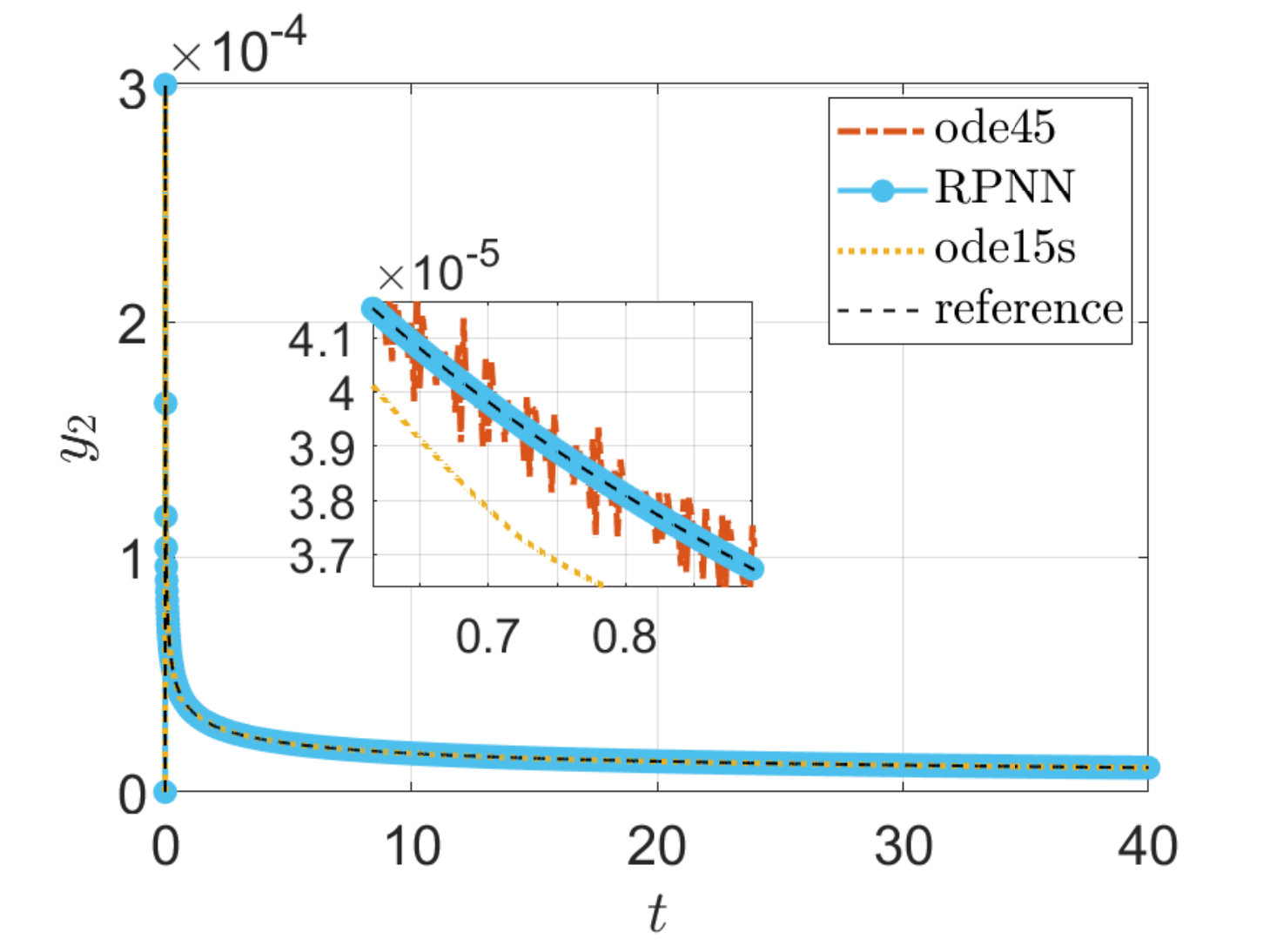}
    }
    par \vskip -10pt
    \subfigure[]{
    \includegraphics[width=0.45 \textwidth]{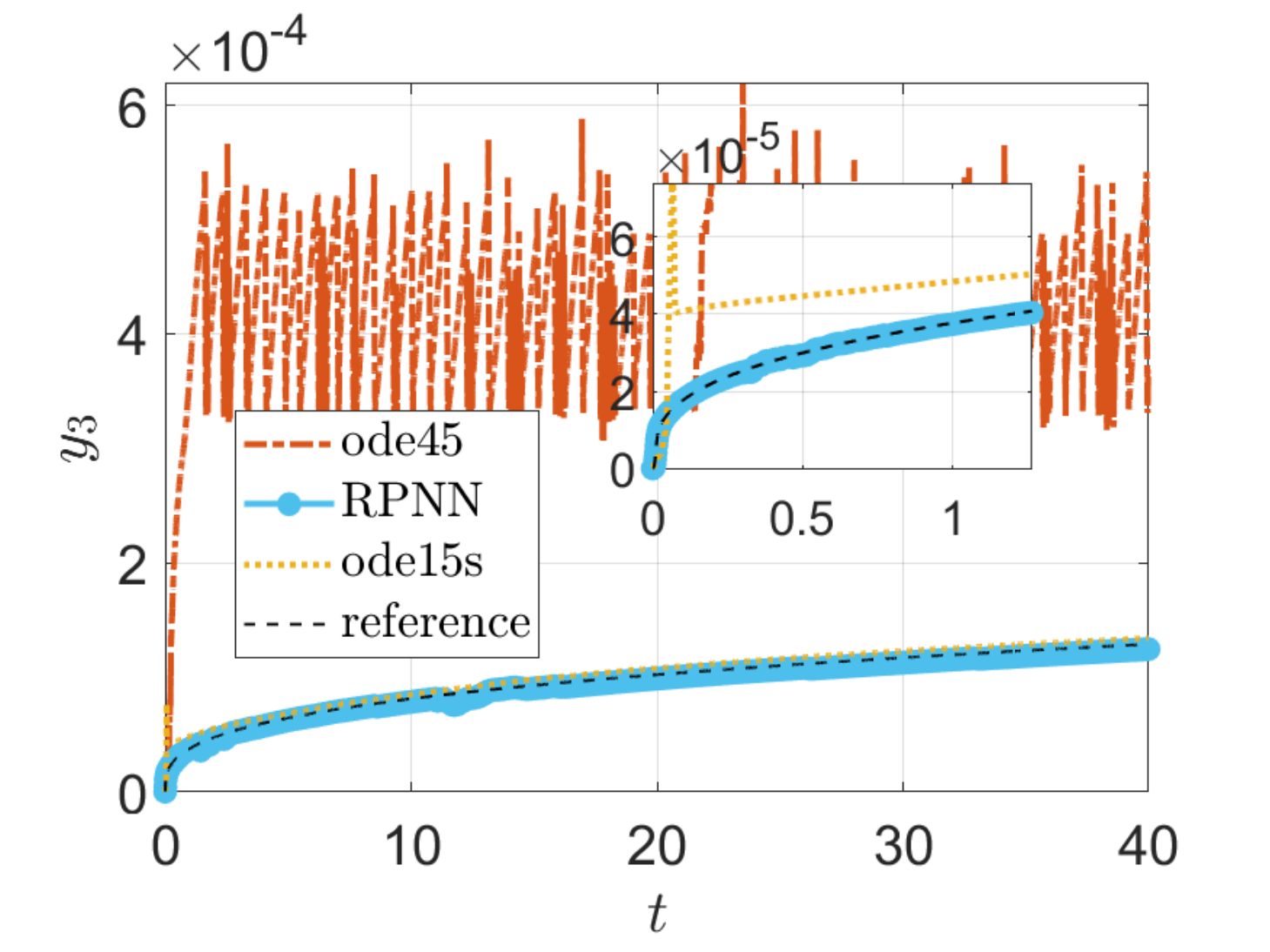}
    }
    \subfigure[]{
    \includegraphics[width=0.45 \textwidth]{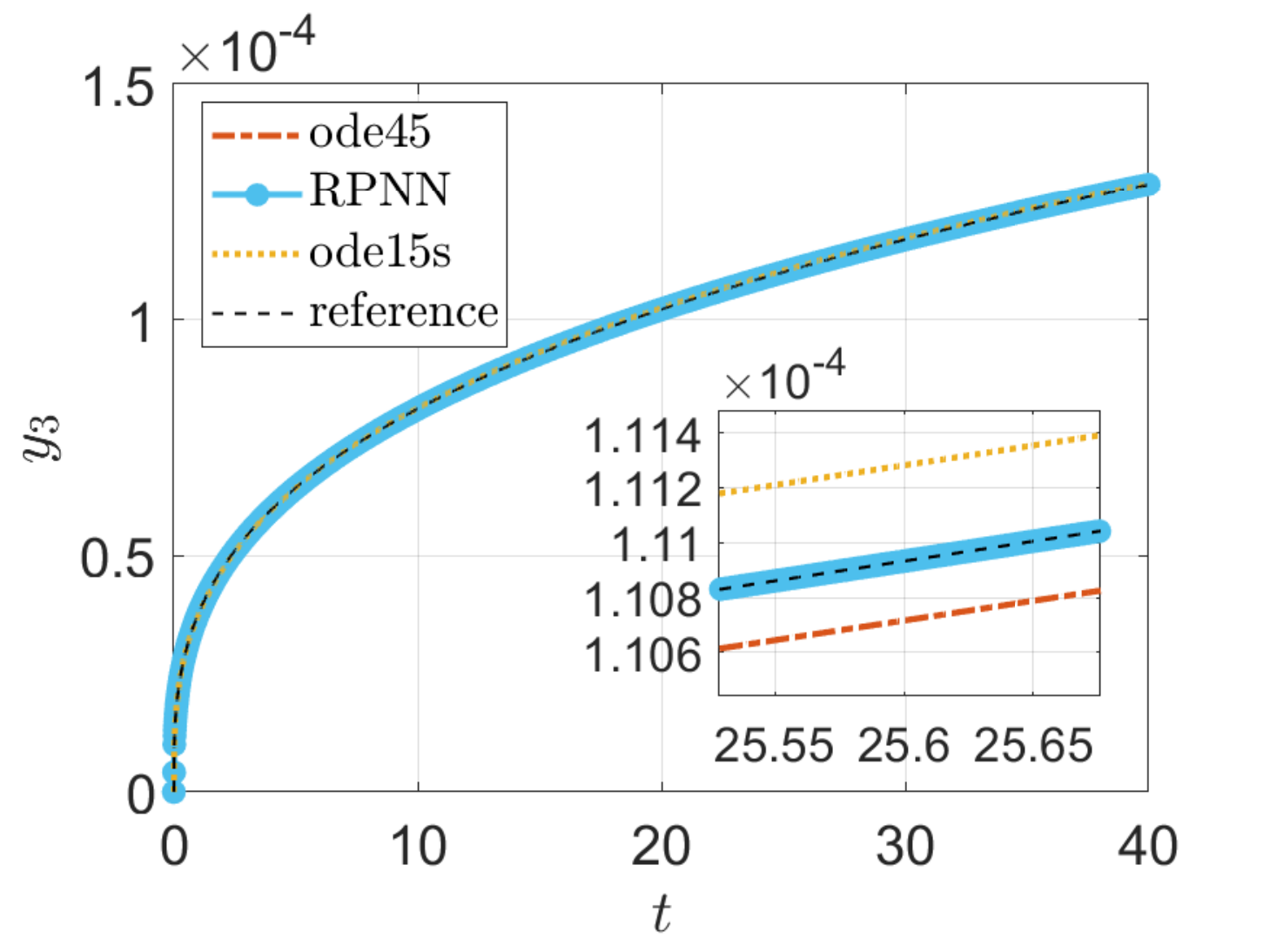}
    }
    par \vskip -10pt
    \caption{ROBER problem. Approximate solutions computed by the different methods in the interval $[0, 40]$ with tolerances set to~1e$-$03 ((a), (c) and (e)) and 1e$-$06 ((b), (d) and (f)). The reference solution was obtained with \texttt{ode15s} with tolerances set to 1e$-$14. The insets depict a zoom around the reference solution.\label{fig:ROBER}}
\end{figure}

Based on the proposed methodology, we construct an (initial) trial solution that satisfies the initial conditions:
\begin{equation*}
\begin{aligned}
\Psi_{1}(t,\bw^o_1) & = \alpha^{(0)}_1+tN_1(t,\bw^o_1),\\
\Psi_{2}(t,\bw^o_2) & = \alpha^{(0)}_2+tN_2(t,\bw^o_2),\\
\Psi_{3}(t,\bw^o_3) & = \alpha^{(0)}_3+tN_3(t,\bw^o_3),
\end{aligned}
\end{equation*}
where $\alpha^{(0)}_1=1$, $\alpha^{(0)}_2=0$ and $\alpha^{(0)}_3=0$ and the $\bp_i$'s have been neglected because they have been fixed.
\begin{table}[b]
\begin{center}
\caption{ROBER problem. Absolute errors ($L_{2}$-norm, $L_{\infty}$-norm and MAE) for the solutions computed with tolerances set to 1e$-$03 and 1e$-$06.
The reference solution was obtained with \texttt{ode15s} with tolerances set to 1e$-$14.\label{tab:ROBER_accuracy}}
{\small
    \begin{tabular}{|l|l |l l l |l l l|}
    \hline
\multicolumn{2}{|c|}{} & \multicolumn{3}{c|}{$tol=$ 1e$-$03} & \multicolumn{3}{c|}{$tol=$ 1e$-$06} \\
\cline{3-8}
        \multicolumn{2}{|c|}{} & $L_2$ & $L_{\infty}$ & MAE & $L_2$ & $L_{\infty}$ & MAE \\
        \hline
        & RPNN   & 1.04e-02 & 1.10e-04 & 6.42e-05 & 1.52e-07 & 1.83e-09 & 9.02e-10\\
        $y_1$ & \texttt{ode45}  & 5.15e-02 & 5.25e-03 & 3.20e-04 & 6.33e-05 & 1.38e-06 & 3.36e-07\\
        & \texttt{ode15s} & 5.76e-03 & 2.21e-03 & 6.21e-06 & 4.21e-05 & 4.97e-06 & 2.39e-07 \\
        \hline
        & RPNN   & 3.51e-04 & 7.91e-06 & 2.18e-06 & 1.52e-08 & 5.19e-10 & 6.65e-11\\
        $y_2$ & \texttt{ode45}  & 2.53e-02 & 4.93e-03 & 1.17e-04 & 5.83e-05 & 1.20e-06 & 3.00e-07\\
        & \texttt{ode15s} & 5.57e-03 & 2.15e-03 & 2.07e-06 & 2.58e-05 & 5.02e-06 & 7.00e-08 \\
        \hline
        & RPNN   & 1.88e-03 & 2.33e-05 & 1.24e-05 & 7.97e-08 & 1.31e-09 & 4.39e-10\\
        $y_3$ & \texttt{ode45}  & 4.66e-02 & 5.17e-04 & 3.23e-04 & 2.76e-05 & 3.38e-07 & 1.69e-07\\
        & \texttt{ode15s} & 8.39e-04 & 6.19e-05 & 5.54e-06 & 4.05e-05 & 4.91e-07 & 2.71e-07 \\
    \hline
    \end{tabular}
}
\end{center}
\end{table}
The approximate solutions obtained with tolerances 1e$-$03 and 1e$-$06 are shown in Figure~\ref{fig:ROBER}.
Furthermore, in Table~\ref{tab:ROBER_accuracy}, we report the corresponding numerical approximation accuracy obtained with the various methods, in terms of $L_2$-norm and $L_{\infty}$-norm errors and MAE, with respect to the reference solution. In order to compute these errors, we evaluated the corresponding solutions in $20\,000$ equidistant grid points in $[0,40]$. In Table~\ref{tab:ROBER_time_points}, we report the number of points and the computational times required by each method, including the time for computing the reference solution.
\begin{table}[ht]
\begin{center}
\caption{ROBER problem. Computational times and number of points required in the interval [0, 40] by RPNN, \texttt{ode45} and \texttt{ode15s} with tolerances~1e$-$03 and~1e$-$06. The reference solution was computed by \texttt{ode15s} with tolerances equal to 1e$-$14.\label{tab:ROBER_time_points}}
{\small
\setlength{\tabcolsep}{3pt}
\begin{tabular}{|l|lllr|lllr|}
\hline
& \multicolumn{4}{c|}{$tol=$ 1e$-$03} & \multicolumn{4}{c|}{$tol=$ 1e$-$06} \\
\cline{2-9}
  & median & min & max & \multicolumn{1}{l|}{\# pts} & median & min & max & \multicolumn{1}{l|}{\# pts}\\
  \hline
  RPNN & 1.52e$-$01 & 1.37e$-$01 & 2.53e$-$01 &  586  & 9.00e$-$02 & 8.32e$-$02 & 1.11e$-$01 & 300\\
  \texttt{ode45}  & 3.47e+00 & 3.35e+00 & 3.65e+00 & 151395 & 6.41e$-$01 & 6.31e$-$01 & 6.96e$-$01 & 34857\\
  \texttt{ode15s} & 2.23e$-$03 & 2.11e$-$03 & 2.60e$-$03 & 34 & 2.71e$-$03 & 2.49e$-$03 & 3.34e$-$03 & 61\\
  reference & 2.54e$-$02 & 2.40e$-$02 & 2.66e$-$02 & 1007  & 2.54e$-$02 & 2.40e$-$02 & 2.66e$-$02 & 1007\\
  \hline
\end{tabular}
}
\end{center}
\end{table}

As it is shown in Figure~\ref{fig:ROBER}, for both tolerances the proposed machine learning method achieves more accurate solutions than \texttt{ode45} and \texttt{ode15s} (actually, as depicted in Figure~\ref{fig:ROBER}, \texttt{ode45} fails to converge, resulting in high frequency oscillations around the reference solution). The proposed method outperforms \texttt{ode15s} and \texttt{ode45} in all metrics (see Table~\ref{tab:ROBER_accuracy}). The number of points required by the proposed machine learning approach is significantly smaller than the number of points required by \texttt{ode45}. On the other hand, the computing times of the proposed method are generally larger than the ones concerning \texttt{ode15s} (Table \ref{tab:ROBER_time_points}), but yet comparable with the ones required by the reference solution. 
\begin{table}[ht]
\begin{center}
\caption{ROBER problem in the interval $[0,40]$. Computational times in seconds (median, minimum and maximum times over 10 runs) for the evaluation (interpolation) of the solution in $1000$ equidistant grid points after the solution was obtained by the three schemes with tolerances set to 1e$-$03 and 1e$-$06. The reference solution was obtained with \texttt{ode15s} with tolerances set to 1e$-$14.}
\label{tab:ROBER_time_points_val}
{\small
\begin{tabular}{|l |l l l |l l l|}
\hline
& \multicolumn{3}{c|}{$tol=$ 1e$-$03} & \multicolumn{3}{c|}{$tol=$ 1e$-$06} \\
\cline{2-7}
& median & min & max & median & min & max \\
\hline
RPNN & 2.00e$-$03 & 1.19e$-$03  & 1.64e$-$02  & 9.52e$-$04 & 8.38e$-$04 & 7.34e$-$03\\
\texttt{ode45}  & 1.49e$-$01  & 1.42e$-$01 & 1.86e$-$01  & 4.52e$-$02 & 4.43e$-$02 & 4.62e$-$02\\
\texttt{ode15s} & 5.68e$-$04  & 5.41e$-$04  & 1.15e$-$03  & 8.21e$-$04 & 8.00e$-$04 & 1.69e$-$03\\
reference & 6.90e$-$03 & 6.03e$-$03 & 7.82e$-$03  & 6.32e$-$03 & 5.85e$-$03 & 7.09e$-$03\\
\hline
\end{tabular}
}
\end{center}
\end{table}
However, as in the van der Pol problem, this is paid off by the fact that our method provides an approximate solution in the form of a function that can be evaluated at points different from the collocation ones, and if it is requested to evaluate the solution in a ``dense'' grid of points, then the proposed machine learning scheme is faster than \texttt{ode15s}. Indeed, in Table \ref{tab:ROBER_time_points_val}, we report the computational times required by all the methods to compute the solution at $1000$ equidistant points after obtaining the solutions with RPNN, \texttt{ode45} and \texttt{ode15s}. As it is shown, the proposed scheme results in computational times comparable with those of \texttt{ode15s} and smaller computational times than \texttt{ode45}.

\subsection{Case Study 4: HIRES problem}

The ``High Irradiance RESponse'' (HIRES) problem proposed by \cite{schafer1975new} is made up of 8 nonlinear stiff ODEs:
\begin{equation}
\begin{split}
    y_1'&=-1.71y_1+0.43y_2+8.32y_3+0.0007 \\
    y_2'&=1.71y_1-8.75y_2\\
    y_3'&=-10.03y_3+0.43y_4+0.035y_5\\
    y_4'&= 8.32y_2+1.71y_3-1.12y_4\\
    y_5'&= -1.745y_5+0.43y_6+0.43y_7\\
    y_6'&= -280y_6y_8+0.69y_4+1.71y_5-0.43y_6+0.69y_7\\
    y_7'&= 280y_6y_8-1.81y_7\\
    y_8'&= -280y_6y_8+1.81y_7
\end{split}
\label{eq:HIRES}
\end{equation}
The initial condition is given by $\boldsymbol{y}(0) = [ 1 \;\, 0 \;\, 0 \;\, 0 \;\, 0 \;\, 0 \;\, 0 \;\, 0.0057 ]^T$, while the solution is sought in the interval $[0, 321.8122]$;
\begin{figure}[p]
    \centering
    par \vskip -30pt
    \subfigure[]{
    \includegraphics[width=0.45 \textwidth]{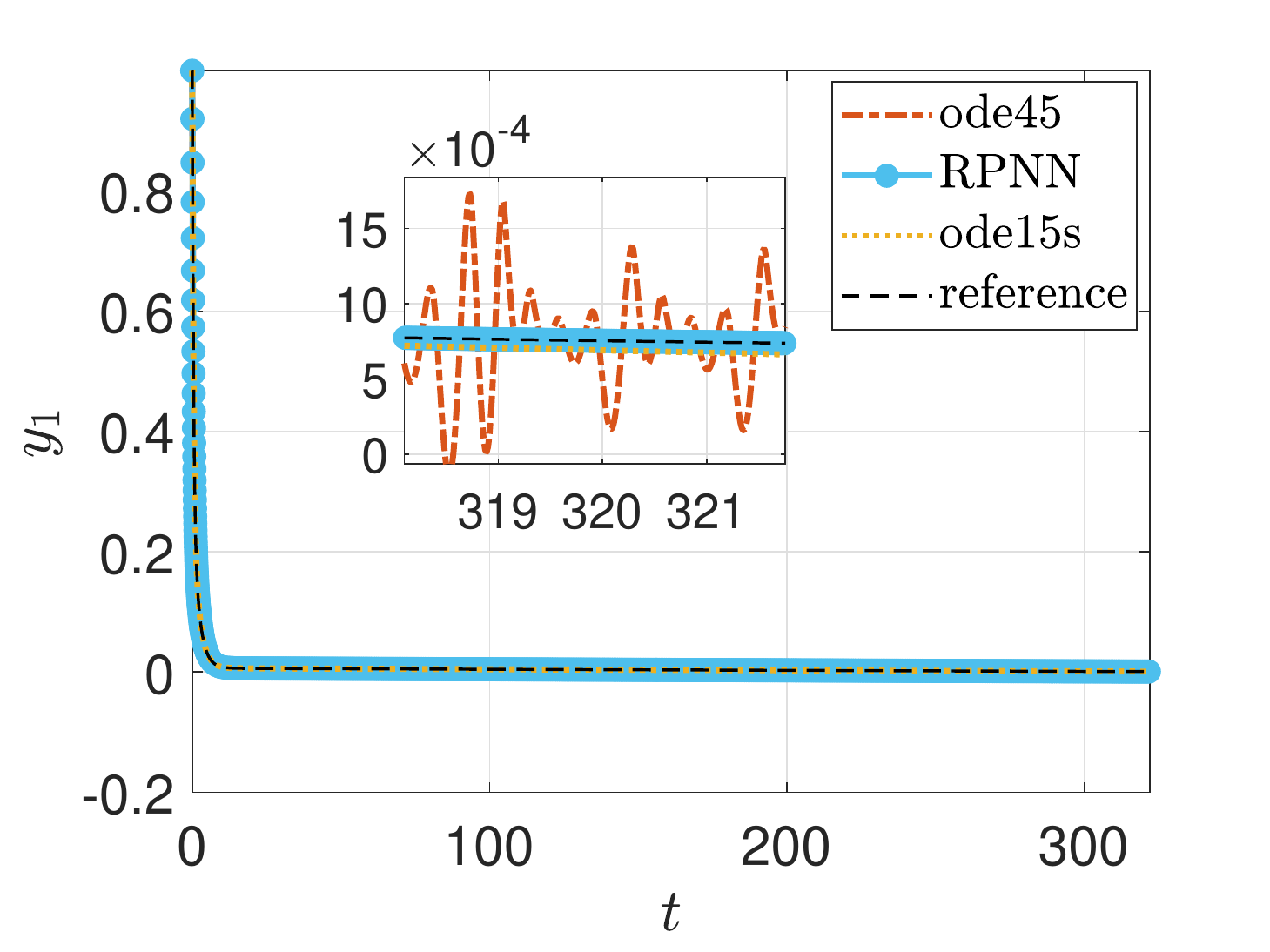}
    }
    \subfigure[]{
    \includegraphics[width=0.45 \textwidth]{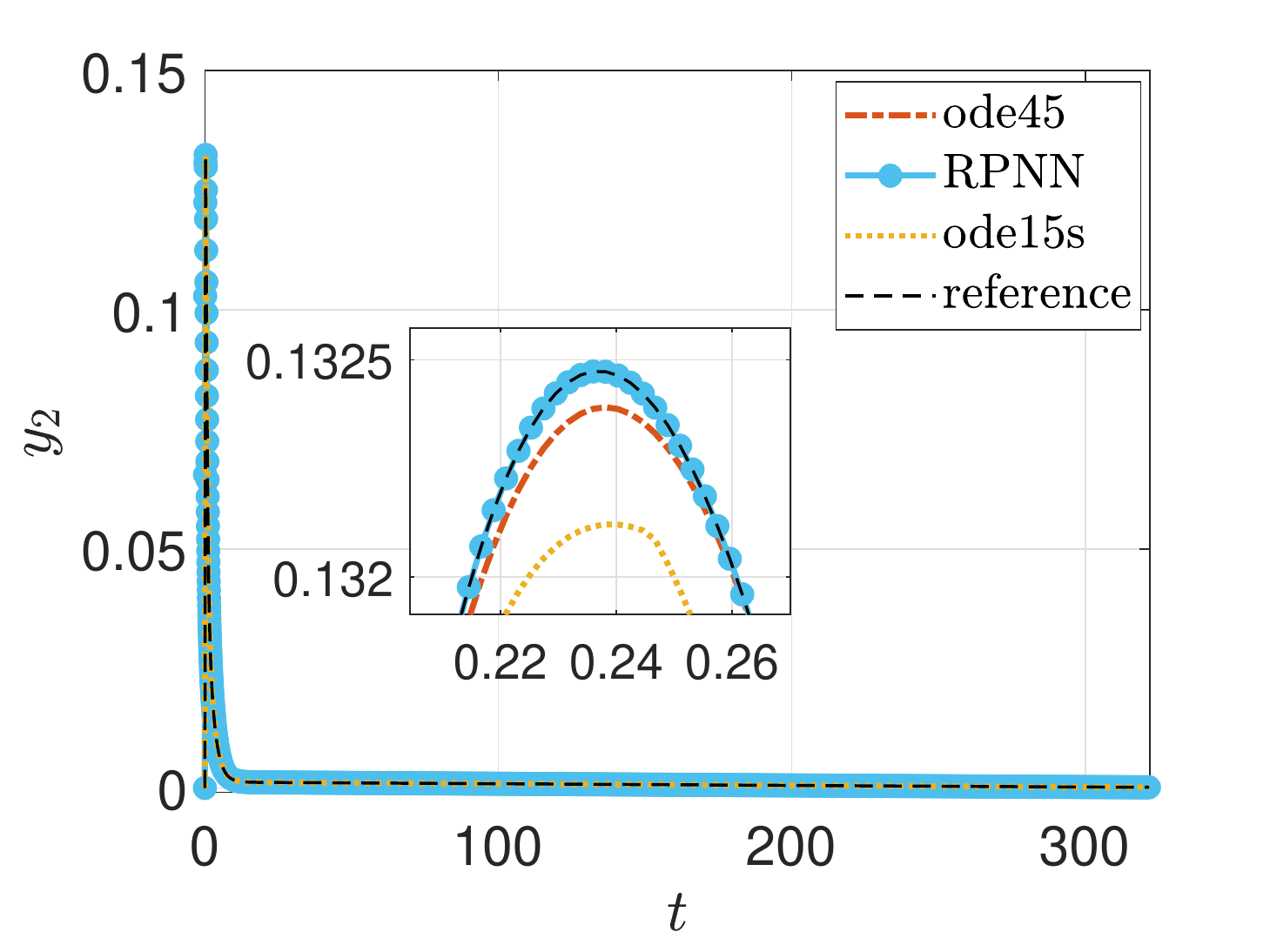}
    }
    par \vskip -10pt
    \subfigure[]{
    \includegraphics[width=0.45 \textwidth]{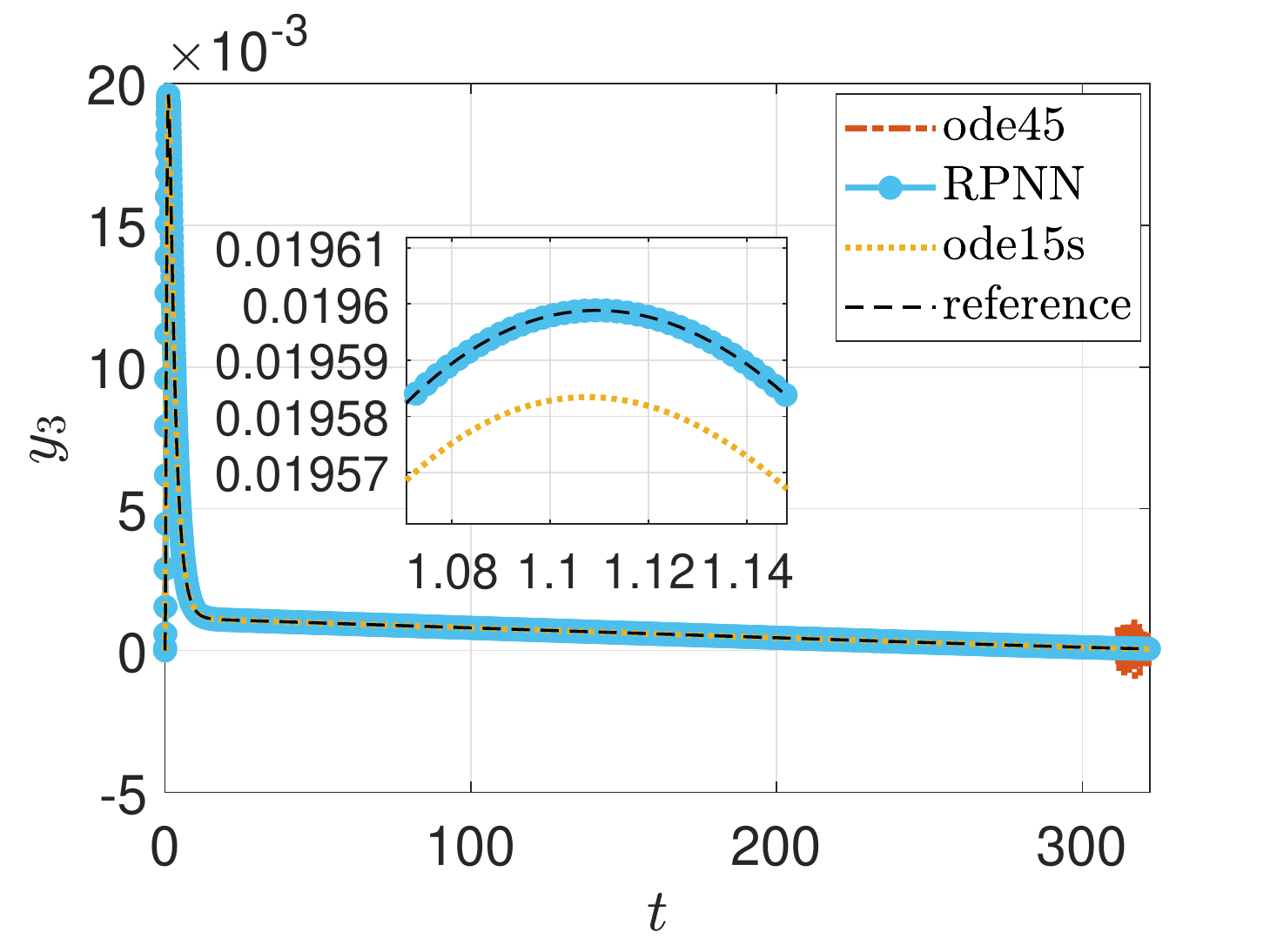}
    }
    \subfigure[]{
    \includegraphics[width=0.45 \textwidth]{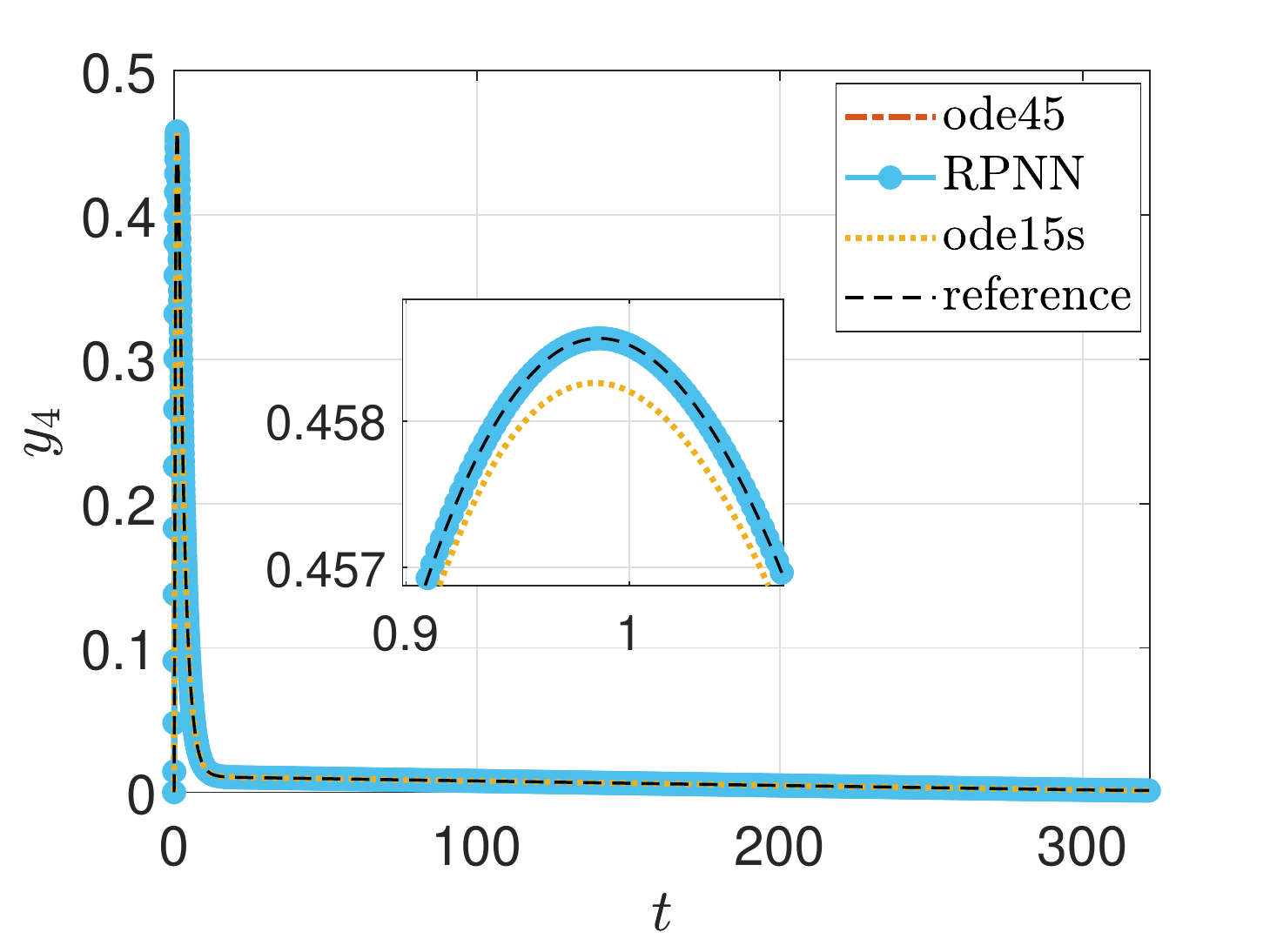}
    }
    par \vskip -10pt
    \subfigure[]{
    \includegraphics[width=0.45 \textwidth]{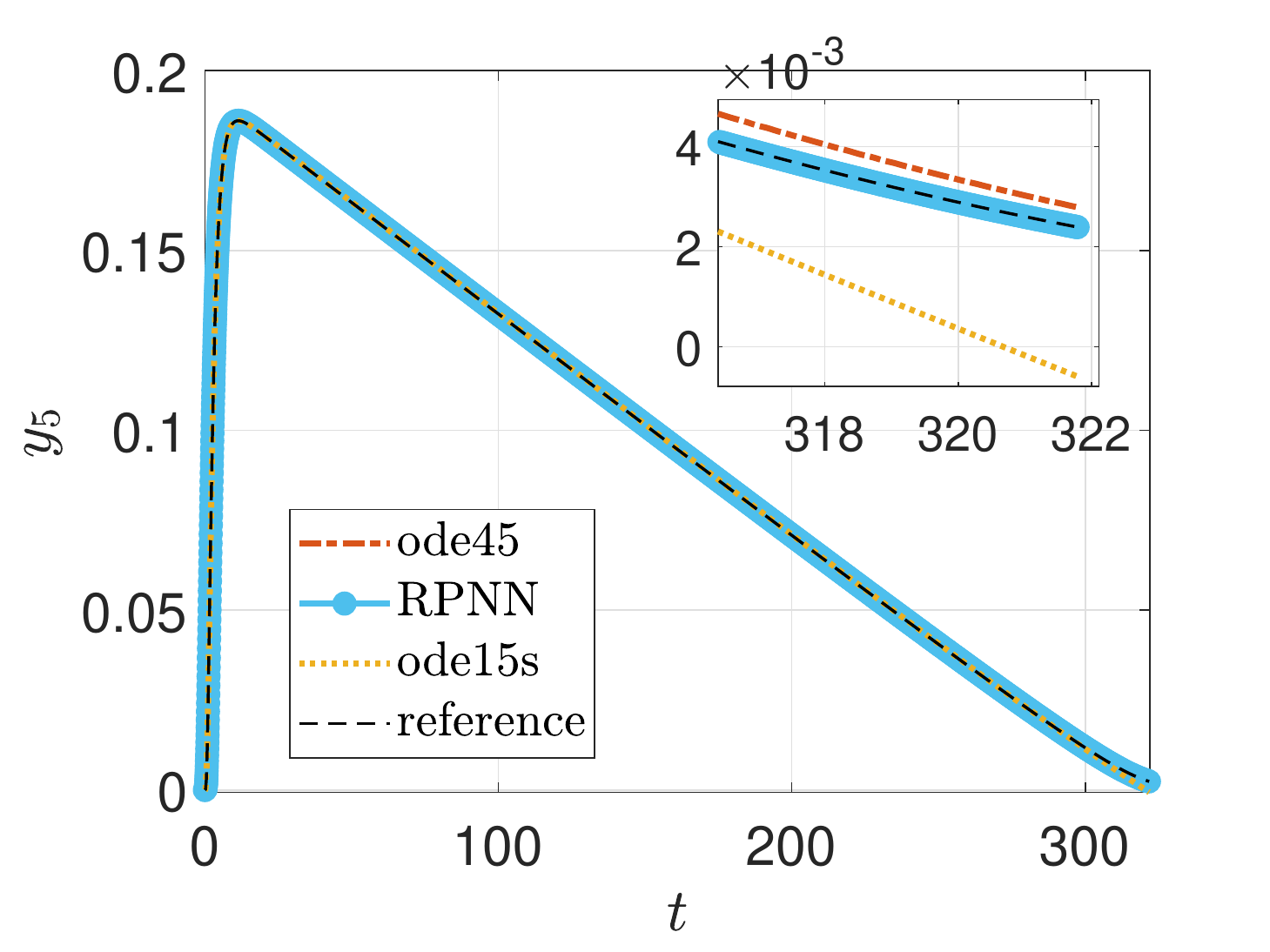}
    }
    \subfigure[]{
    \includegraphics[width=0.45 \textwidth]{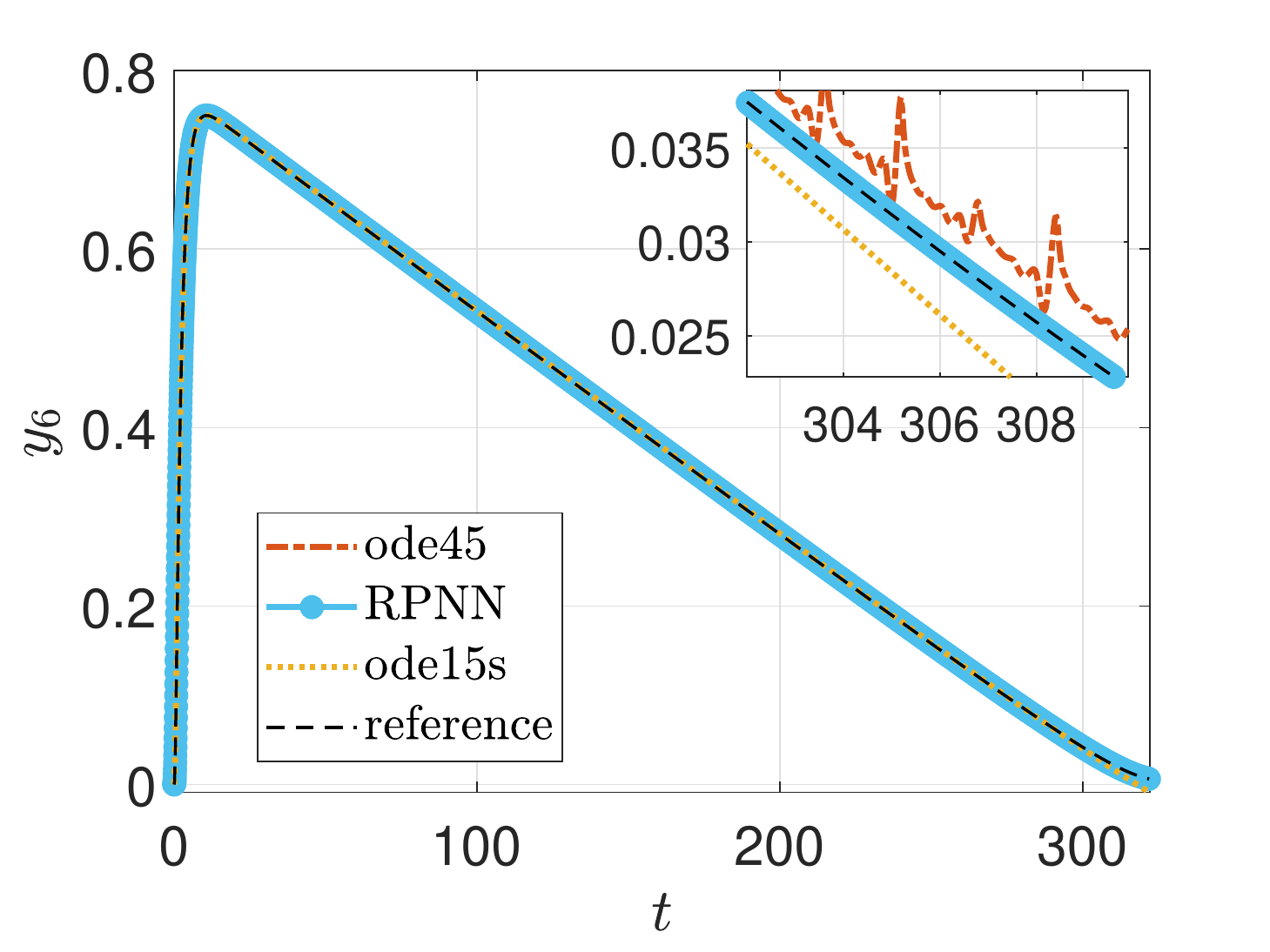}
    }
    par \vskip -10pt
    \subfigure[]{
    \includegraphics[width=0.45 \textwidth]{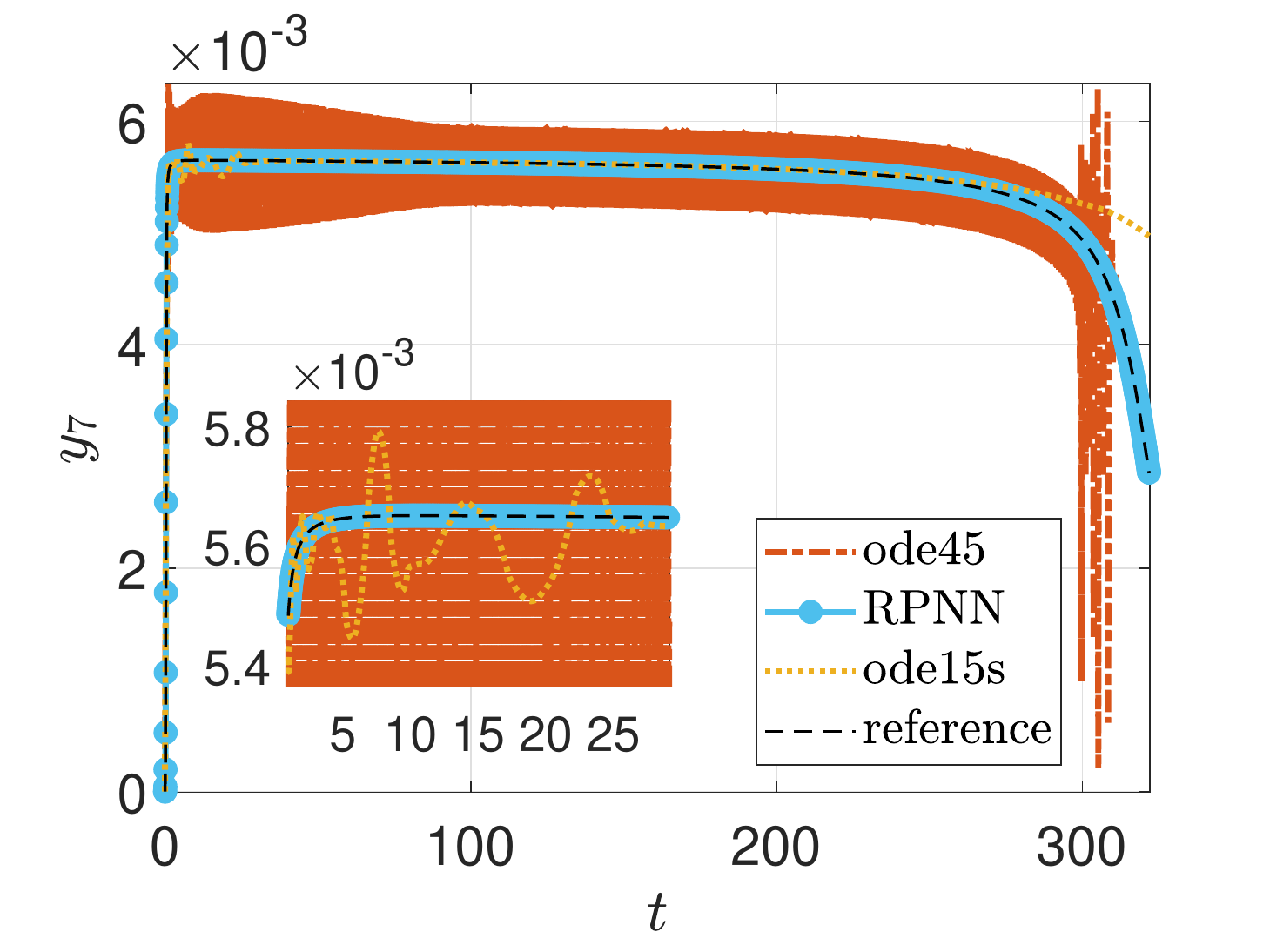}
    }
    \subfigure[]{
    \includegraphics[width=0.45 \textwidth]{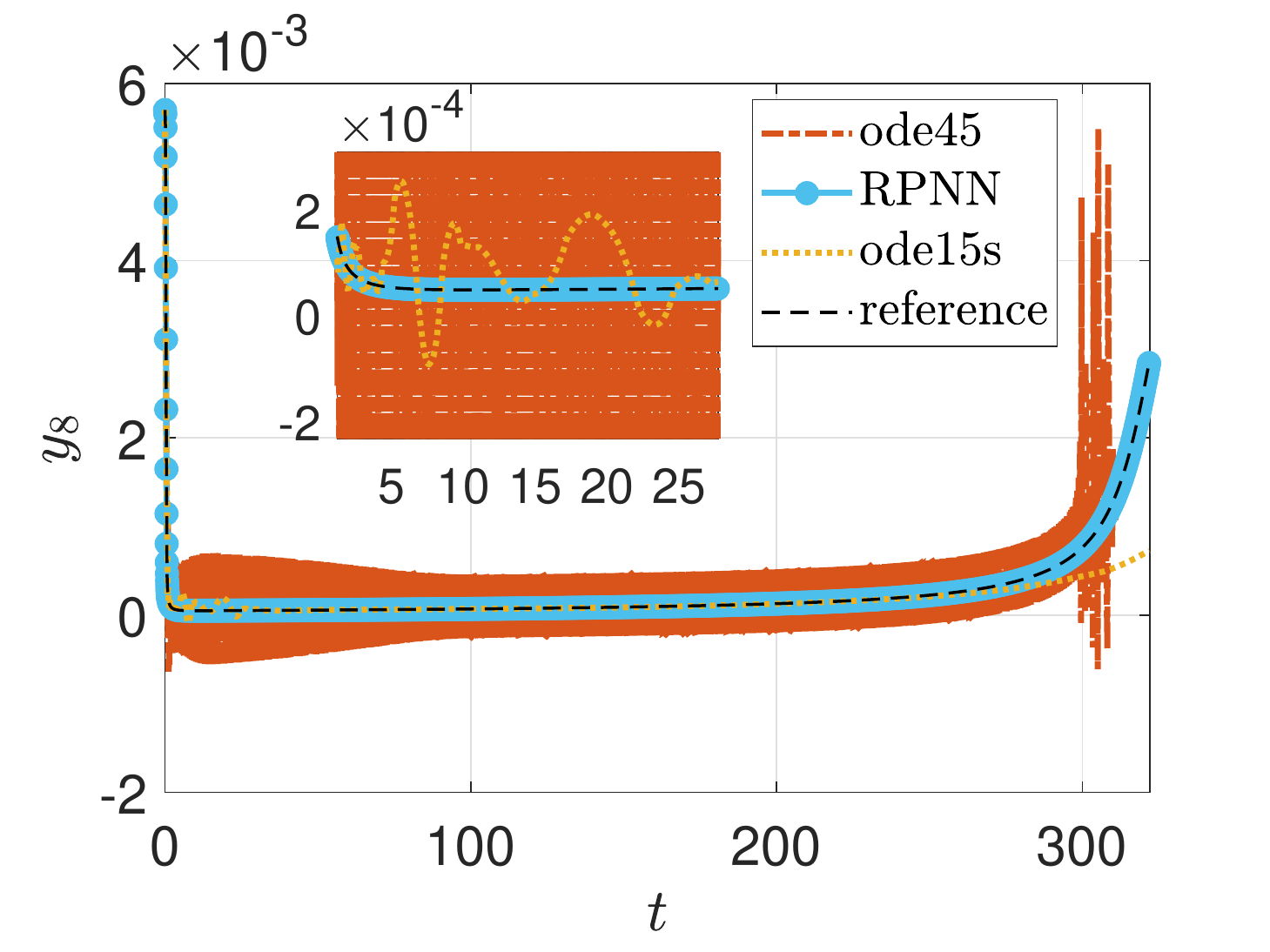}
    }
    par \vskip -10pt
    \caption{HIRES problem. Approximate solutions computed by the different methods in the interval $[0, 321.8122]$ with tolerances set to~1e$-$03. In panels (a) to (h) the eight components $y_i$ of the solution are depicted. The reference solution was obtained with \texttt{ode15s} with tolerances set to 1e$-$14. The insets depict a zoom around the reference solution.\label{fig:Hires_tol3}}
\end{figure}
\begin{figure}[p!]
    \centering
    par \vskip -30pt
    \subfigure[]{
    \includegraphics[width=0.45 \textwidth]{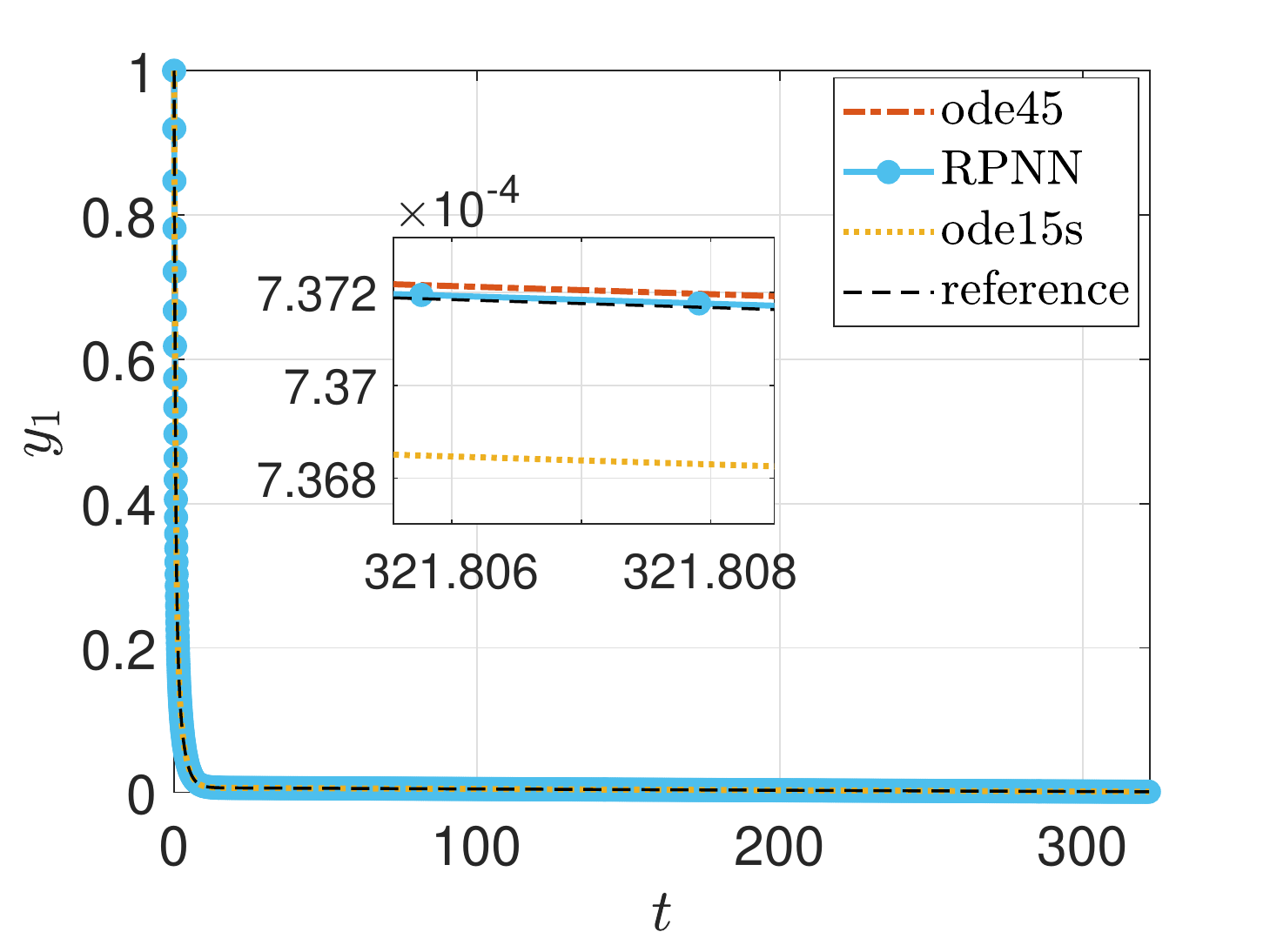}
    }
    \subfigure[]{
    \includegraphics[width=0.45 \textwidth]{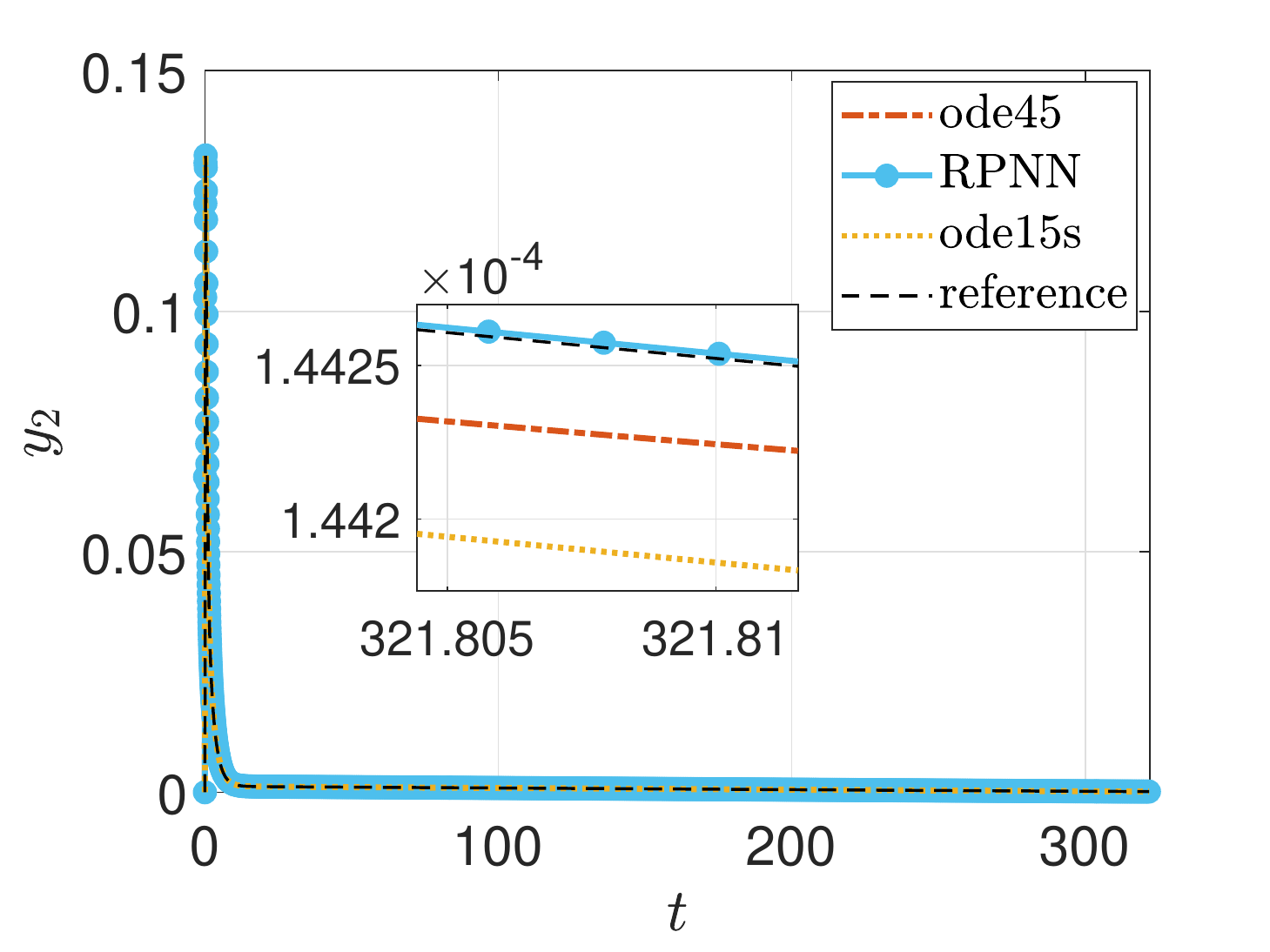}
    }
    par \vskip -10pt
    \subfigure[]{
    \includegraphics[width=0.45 \textwidth]{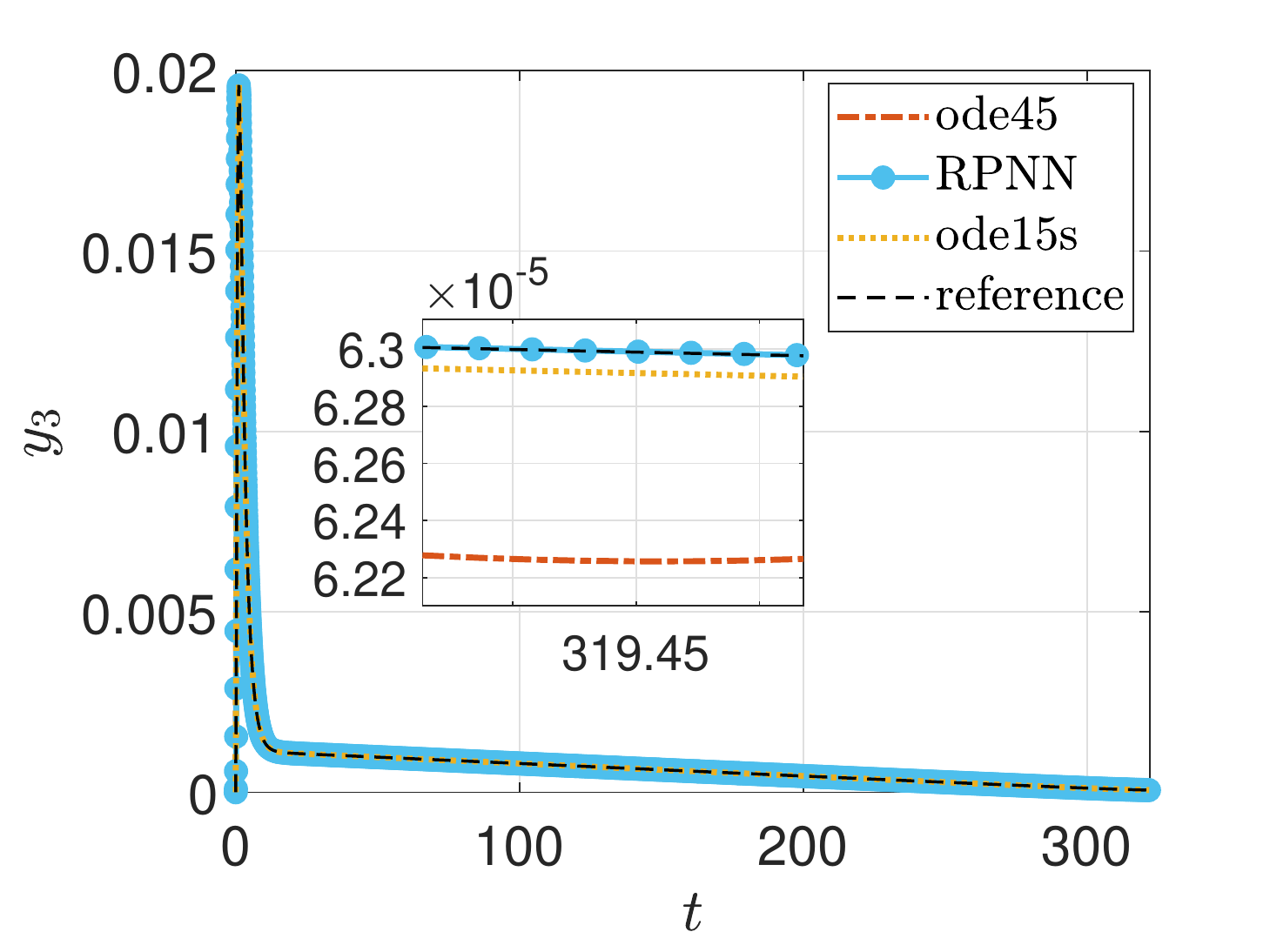}
    }
    \subfigure[]{
    \includegraphics[width=0.45 \textwidth]{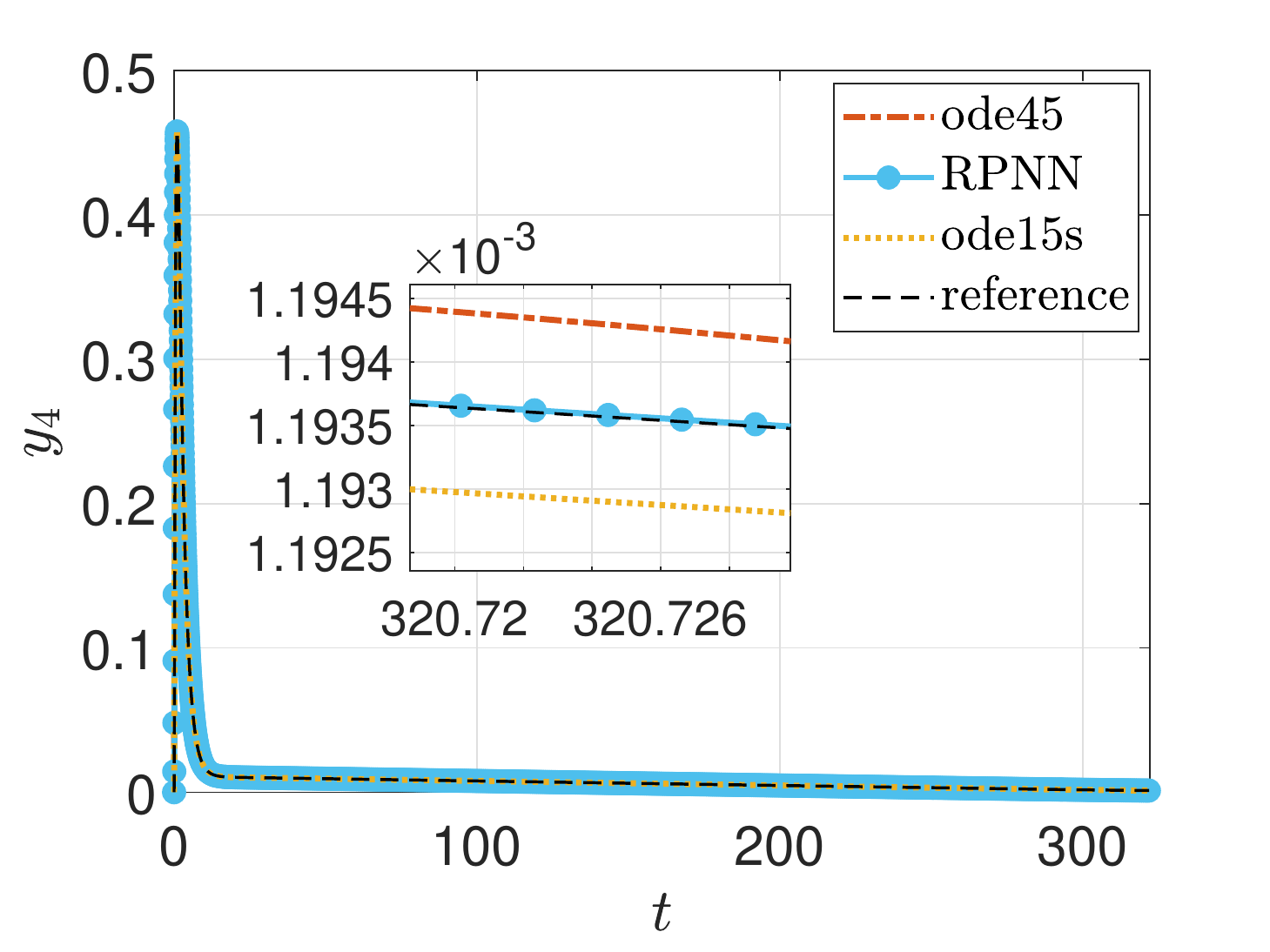}
    }
    par \vskip -10pt
    \subfigure[]{
    \includegraphics[width=0.45 \textwidth]{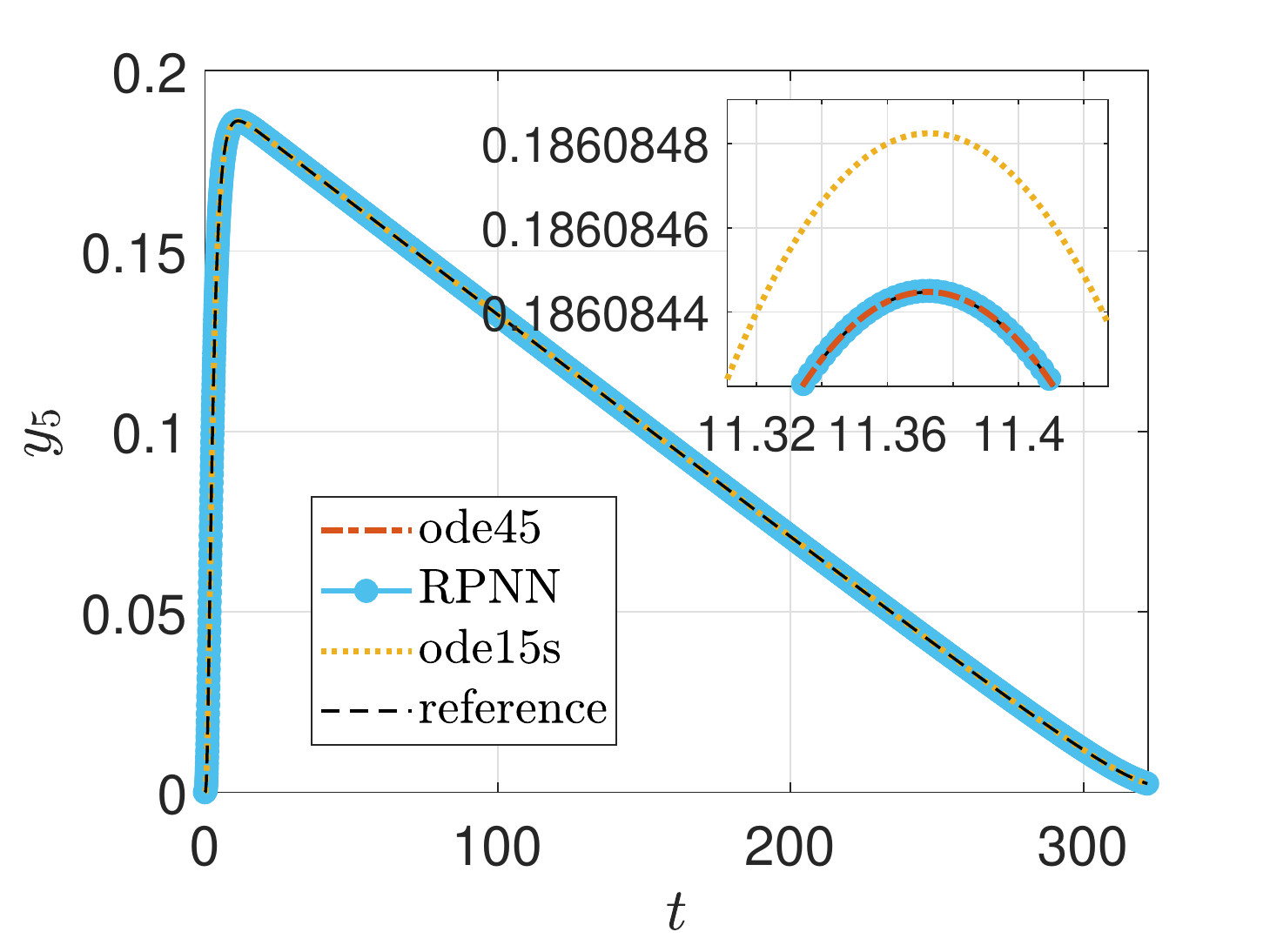}
    }
    \subfigure[]{
    \includegraphics[width=0.45 \textwidth]{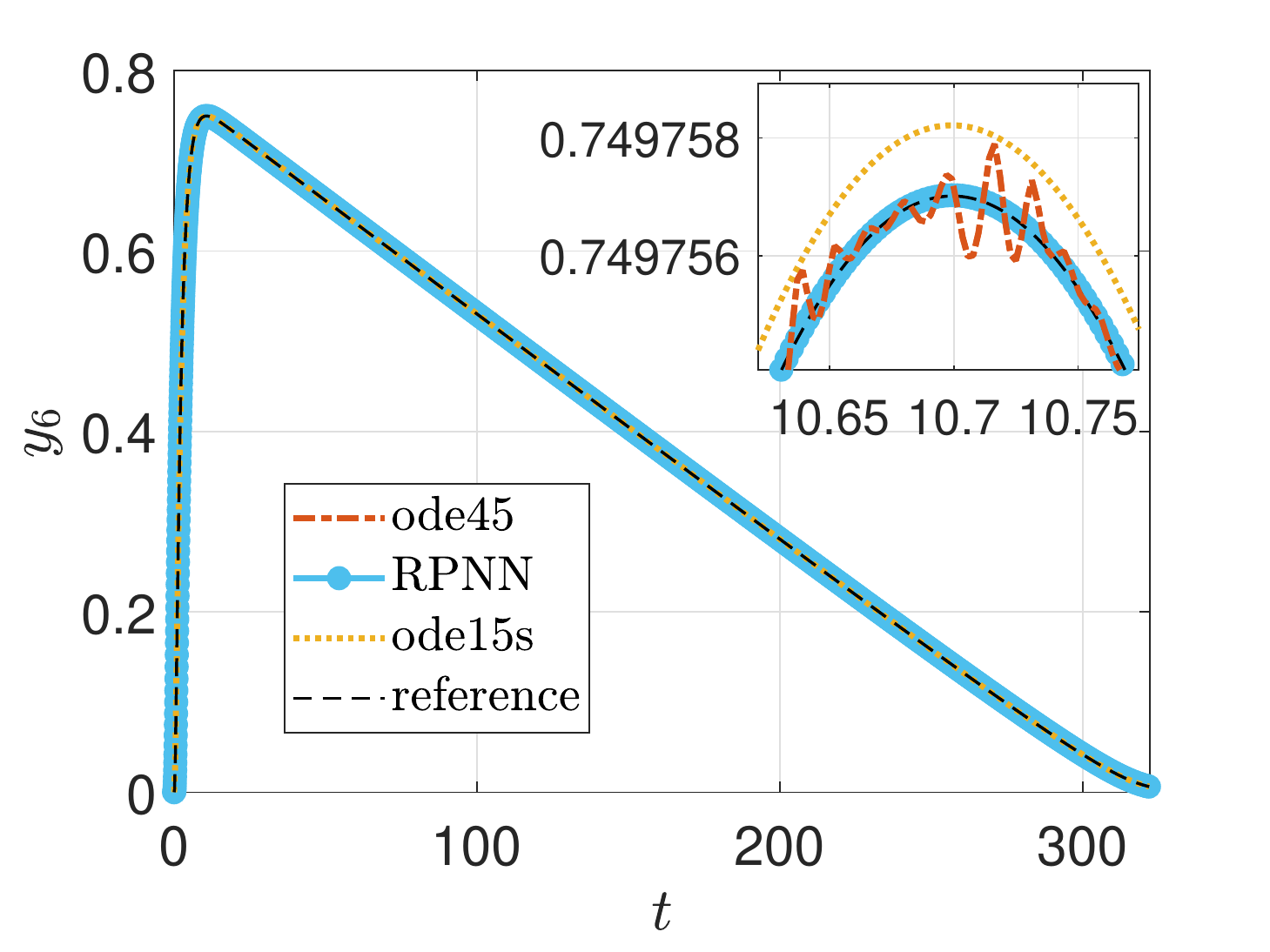}
    }
    par \vskip -10pt
    \subfigure[]{
    \includegraphics[width=0.45 \textwidth]{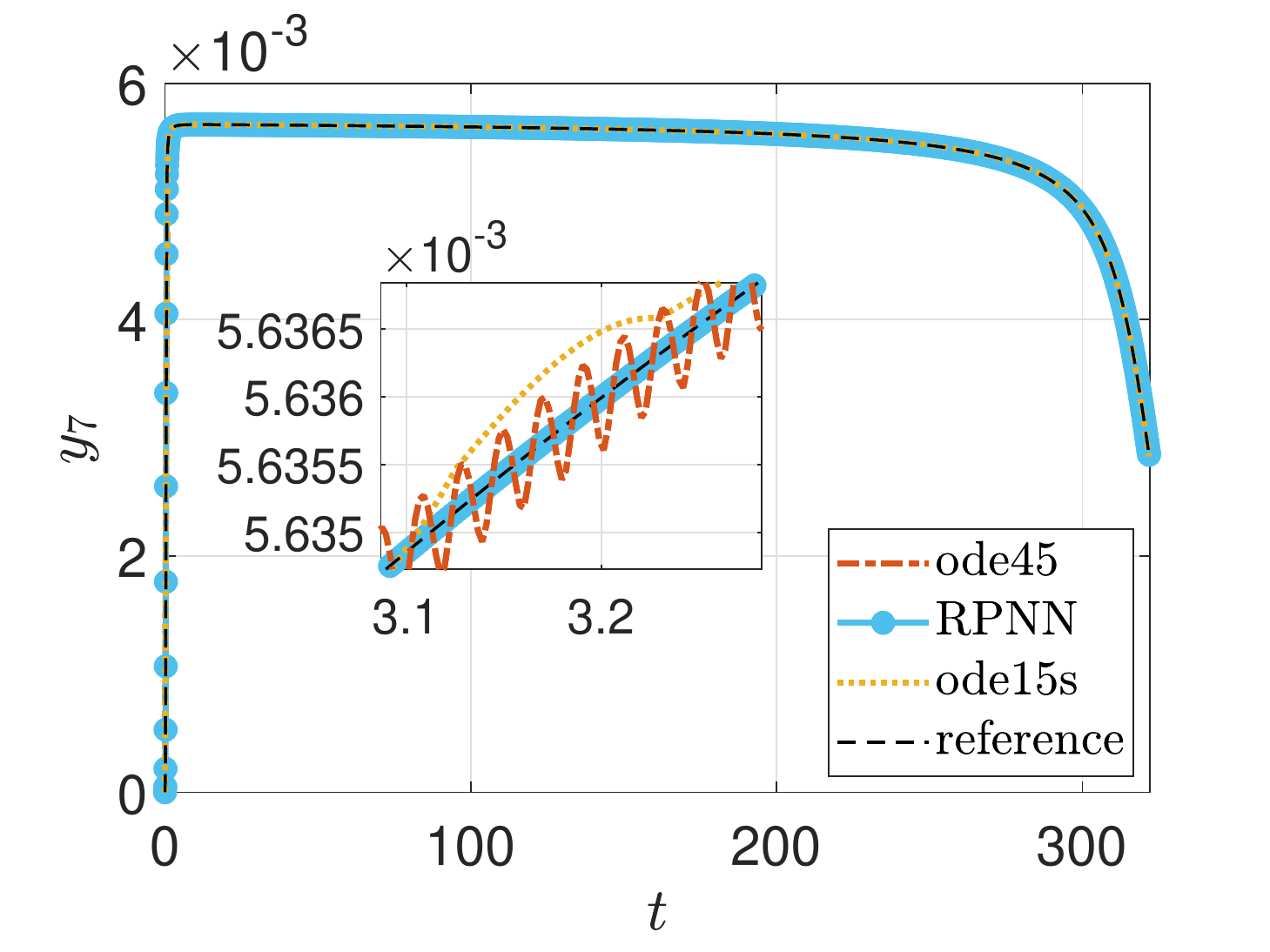}
    }
    \subfigure[]{
    \includegraphics[width=0.45 \textwidth]{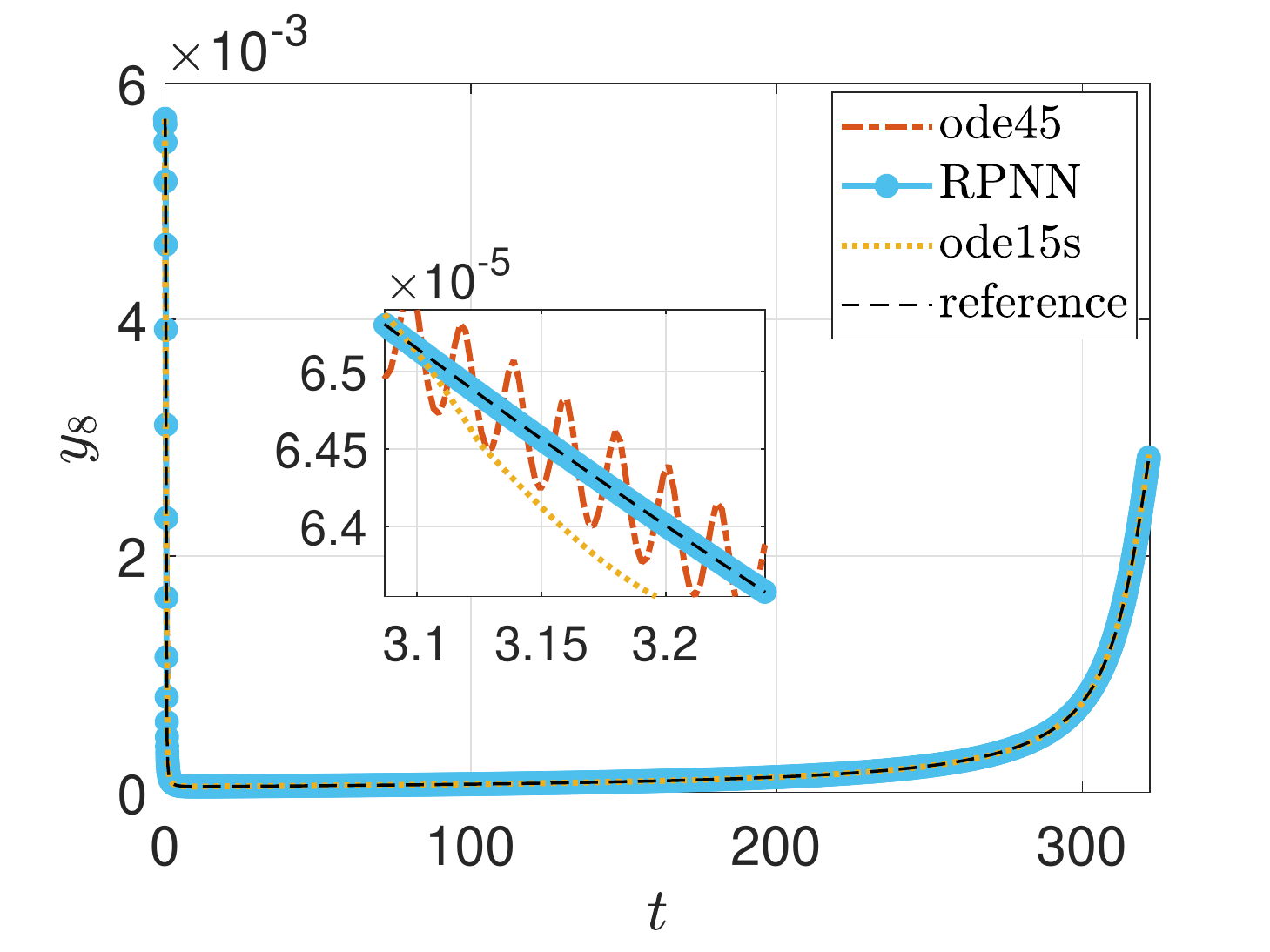}
    }
    par \vskip -10pt
    \caption{HIRES problem. Approximate solutions computed by the different methods in the interval $[0, 321.8122]$ with tolerances set to~1e$-$06. In panel from (a) to (h) are depicted the eight components $y_i$ of the solution. The reference solution was obtained with \texttt{ode15s} with tolerances set to 1e$-$14. The insets depict a zoom around the reference solution.\label{fig:Hires_tol6}}
\end{figure}

The approximate solutions obtained with tolerances 1e$-$03 are shown in Figure~\ref{fig:Hires_tol3} and the ones obtained with tolerances 1e$-$06 in Figure~\ref{fig:Hires_tol6}. In Table~\ref{tab:Hires_accuracy}, we report the corresponding numerical approximation accuracy obtained with the various methods, in terms of $L_2$-norm and $L_{\infty}$-norm of the error and of MAE, with respect to the reference solution. In order to compute these errors, we evaluated the corresponding solutions in 150000 equidistant grid points in $[0, 321.8122]$. In Table~\ref{tab:Hires_time_points}, we report the number of points and the computational times for each method, including the time for computing the reference solution.

\begin{table}[ht]
\begin{center}
\caption{HIRES problem. Absolute error ($L_2$-norm, $L_{\infty}$-norm and MAE) for the solutions computed with tolerances set to 1e$-$03 and 1e$-$06. The reference solution was obtained with \texttt{ode15s} with tolerances set to 1e$-$14.\label{tab:Hires_accuracy}}
    {\small
    \begin{tabular}{|l|l |l l l |l l l|}
        \hline
\multicolumn{2}{|c|}{} & \multicolumn{3}{c|}{$tol=$ 1e$-$03} & \multicolumn{3}{c|}{$tol=$ 1e$-$06} \\
\cline{3-8}
        & & $L_2$ & $L_{\infty}$ & MAE & $L_2$ & $L_{\infty}$ & MAE\\
        \hline
        & RPNN   & 1.25e$-$05 & 2.15e$-$07 & 2.73e$-$08 & 4.90e$-$06 & 2.54e$-$08 & 9.36e$-$09\\
        $y_1$ & \texttt{ode45}  & 2.64e$-$02 & 1.08e$-$03 & 1.14e$-$05 & 2.82e$-$05 & 1.07e$-$06 & 1.03e$-$08\\
        & \texttt{ode15s} & 1.12e$-$02 & 3.03e$-$04 & 1.14e$-$05 & 7.15e$-$05 &  9.69e$-$07 & 9.82e$-$08 \\
        \hline
        & RPNN   & 2.78e$-$06 & 1.70e$-$07 & 5.53e$-$09 & 9.65e$-$07 & 7.10e$-$09 & 1.84e$-$09\\
        $y_2$ & \texttt{ode45}  & 2.61e$-$02 & 1.05e$-$03 & 9.78e$-$06 & 2.78e$-$05 & 1.05e$-$06 & 1.03e$-$08\\
        & \texttt{ode15s} & 6.16e$-$03 & 9.83e$-$04 & 2.80e$-$06 & 2.05e$-$05  &  2.09e$-$06  &  2.14e$-$08 \\
        \hline
        & RPNN   & 2.37e$-$06 & 2.97e$-$08 & 5.25e$-$09 & 9.61e$-$07 & 5.08e$-$09    & 1.84e$-$09\\
        $y_3$ & \texttt{ode45}  & 2.65e$-$02 & 1.06e$-$03 & 9.82e$-$06 & 2.83e$-$05 & 1.07e$-$06 & 1.03e$-$08\\
        & \texttt{ode15s} & 1.37e$-$03 & 2.60e$-$05 & 1.75e$-$06 & 1.30e$-$05 & 1.80e$-$07 & 1.75e$-$08 \\
        \hline
        & RPNN   & 2.21e$-$05 & 2.68e$-$07 & 4.81e$-$08 & 8.58e$-$06 & 4.47e$-$08 & 1.64e$-$08\\
        $y_4$ & \texttt{ode45}  & 2.81e$-$02 & 1.16e$-$03 & 1.35e$-$05 & 2.99e$-$05 & 1.13e$-$06 & 1.11e$-$08\\
        & \texttt{ode15s} & 2.47e$-$02 & 7.21e$-$04 & 2.33e$-$05 & 1.40e$-$04 & 2.18e$-$06 & 2.02e$-$07 \\
        \hline
        & RPNN   & 4.10e$-$04 & 1.67e$-$06 & 8.99e$-$07 & 1.70e$-$04 & 8.59e$-$07 & 3.26e$-$07\\
        $y_5$ & \texttt{ode45}  & 6.23e$-$02 & 6.58e$-$04 & 9.15e$-$05 & 2.65e$-$07 & 6.03e$-$09 & 6.13e$-$10\\
        & \texttt{ode15s} & 1.57e$-$01 & 2.99e$-$03 & 2.12e$-$04 & 2.20e$-$03 & 2.00e$-$05 & 2.51e$-$06 \\
        \hline
        & RPNN & 1.67e$-$03 & 6.76e$-$06 & 3.67e$-$06 & 6.94e$-$04 & 3.50e$-$06 &    1.34e$-$06\\
        $y_6$ & \texttt{ode45}  & 2.84e$-$01 & 6.62e$-$03 & 4.59e$-$04 & 1.53e$-$04 & 1.20e$-$06 & 2.83e$-$07\\
        & \texttt{ode15s} & 7.72e$-$01 & 1.49e$-$02 & 9.48e$-$04 & 9.11e$-$03 & 9.11e$-$05 & 1.05e$-$05 \\
        \hline
        & RPNN   & 3.74e$-$05 & 6.97e$-$07 & 5.96e$-$08 & 1.99e$-$05 & 3.00e$-$07 & 3.11e$-$08\\
        $y_7$ & \texttt{ode45}  & 1.17e$-$01 & 4.51e$-$03 & 2.35e$-$04 & 1.55e$-$04 & 1.21e$-$06 & 2.86e$-$07\\
        & \texttt{ode15s} & 1.11e$-$01 & 2.12e$-$03 & 8.63e$-$05 & 6.70e$-$04 & 9.44e$-$06 & 8.04e$-$07 \\
        \hline
        & RPNN   & 3.70e$-$05 & 7.36e$-$07 & 2.87e$-$08 & 2.02e$-$05 & 3.43e$-$07 & 1.59e$-$08\\
        $y_8$ & \texttt{ode45}  & 1.17e$-$01 & 4.51e$-$03 & 2.35e$-$04 & 1.55e$-$04 & 1.21e$-$06 & 2.86e$-$07\\
        & \texttt{ode15s} & 1.11e$-$01 & 2.12e$-$03 & 8.63e$-$05 & 6.70e$-$04 & 9.44e$-$06 & 8.04e$-$07 \\
        \hline
    \end{tabular}
    }
\end{center}
\end{table}

As it is shown in Figures~\ref{fig:Hires_tol3}-\ref{fig:Hires_tol6}, the proposed machine learning method achieves more accurate solutions than \texttt{ode45} and \texttt{ode15s}, generally outperforming these two solvers in all metrics (see Table~\ref{tab:Hires_accuracy}). The solution computed by \texttt{ode45} oscillates around the reference solution, thus resulting in poor performance when compared with the other methods.
\begin{table}[ht]
\begin{center}
\caption{HIRES problem. Computational times in seconds (median, minimum and maximum times over 10 runs) and Number of points required in the interval [0, 321.8122] by RPNN, \texttt{ode45} and  \texttt{ode15s} with tolerances~1e$-$03 and~1e$-$06. The reference solution was computed by \texttt{ode15s} with tolerances equal to 1e$-$14.\label{tab:Hires_time_points}}
    {\small
    \setlength{\tabcolsep}{3pt}
    \begin{tabular}{|l |l l l r |l l l r|}
        \hline
        & \multicolumn{4}{c|}{$tol=$ 1e$-$03} & \multicolumn{4}{c|}{$tol=$ 1e$-$06} \\
        \cline{2-9}
        & median & min & max & \multicolumn{1}{l|}{\# pts} & median & min & max & \multicolumn{1}{l|}{\# pts}\\
        \hline
        RPNN & 3.82e$-$01 & 3.75e$-$01 & 4.04e$-$01 &  280  & 5.99e$-$01   & 5.38e$-$01 & 8.58e$-$01 &  340\\
        \texttt{ode45}  & 2.65e$-$01 & 2.62e$-$01 & 2.73e$-$01 &  10338  & 2.94e$-$01 & 2.76e$-$01 & 4.59e$-$01 &  10382\\
        \texttt{ode15s} & 3.05e$-$03 & 2.95e$-$03 & 3.47e$-$03 &  61  & 7.21e$-$03 & 6.54e$-$03 & 3.29e$-$02 &  201\\
        reference & 8.87e$-$02 & 8.75e$-$02 & 1.06e$-$01 &  3274 & 8.87e$-$02 & 8.75e$-$02 & 1.06e$-$01 &  3274\\
        \hline
    \end{tabular}
    }
\end{center}
\end{table}
For example, when the tolerances are set to 1e$-$03, the $L_{\infty}$-norm of the difference between the value of $y_7$ ($y_8$) provided by the schemes and the reference solution are as follows: with \texttt{ode45} $\| y_7-y^{ref}_7 \|_{L_{\infty}} \simeq 5$e$-$03 ($\| y_8-y^{ref}_8 \|_{L_{\infty}} \simeq 5$e$-$03), with \texttt{ode15s} $\| y_7-y^{ref}_7 \|_{L_{\infty}} \simeq 2$e$-$03 ($\| y_8-y^{ref}_8 \|_{L_{\infty}} \simeq 5$e$-$03), and with the proposed scheme $\| y_1 -y^{ref}_7\|_{L_{\infty}} \simeq 7$e$-$07 ($\| y_2 -y^{ref}_8\|_{L_{\infty}} \simeq 7$e$-$07). When the tolerances are set equal to 1e$-$06, we have a similar behaviour: the proposed method provides better approximations than \texttt{ode45} and \texttt{ode15s}.

The number of points used by our approach is significantly smaller that the number of points used by \texttt{ode45}. On the other hand, the computing times of the proposed scheme are larger than those of \texttt{ode15s} (Table~\ref{tab:Hires_time_points}). However, as for the other problems, if the solution has to be evaluated at many points, then our method is faster than the other two methods.
\begin{table}[ht]
\begin{center}
\caption{HIRES problem in the interval [0, 321.8122]. Computational times in seconds (median, minimum and maximum times over 10 runs) for the evaluation (interpolation) of the solution in $8000$ equidistant grid points after the solution was obtained by the three schemes with tolerances 1e$-$03 and 1e$-$06. The reference solution was obtained with \texttt{ode15s} with tolerances set to 1e$-$14.\label{tab:Hires_time_points_val}}
{\small
\begin{tabular}{|l |l l l |l l l|}
\hline
& \multicolumn{3}{c|}{$tol=$ 1e$-$03} & \multicolumn{3}{c|}{$tol=$ 1e$-$06} \\
\cline{2-7}
& median & min & max & median & min & max \\
\hline
RPNN & 4.06e$-$03 & 3.80e$-$03  & 7.96e$-$03  & 4.28e$-$03 & 3.95e$-$03 & 4.88e$-$03\\
\texttt{ode45}  & 3.22e$-$01 & 3.20e$-$01 & 3.31e$-$01  & 3.53e$-$01 & 3.48e$-$01 & 3.72e$-$01\\
\texttt{ode15s} & 5.32e$-$03 & 4.58e$-$03  & 6.77e$-$03  & 1.15e$-$02 & 1.06e$-$02 & 5.74e$-$02\\
reference & 7.93e$-$02 & 7.75e$-$02 & 1.00e$-$01  & 8.79e$-$02 & 8.17e$-$02 & 9.54e$-$02\\
\hline
\end{tabular}
}
\end{center}
\end{table}
Indeed, in Table~\ref{tab:Hires_time_points_val}, we report the times required by all methods to compute the solution at $8000$ equidistant points. We see that the proposed machine learning scheme outperforms all the other methods in terms of computational times for both tolerances.

\section{Discussion\label{sec:discussion}}

We proposed a machine learning algorithm for the solution of ODEs with a focus on stiff ODEs, which combines the speed and the generalization advantages of RPNNs. In this algorithm, only the parameters from the hidden to the output layer have to be determined, and every hidden unit is ``responsible'' of learning the behaviour of the underlying physical model around a center point.
Theoretical results stemming from the Johnson and Lindenstrauss Theorem 
and universal approximation theorems that have been proved for ELMs 
guarantee the approximation capabilities of our RPNN, despite the much simpler way of obtaining the weights of the network.

The results show that the proposed machine learning approach is able to provide reasonably accurate numerical approximations to the solutions of four benchmark stiff ODE problems. Our algorithm outperforms \texttt{ode15s} in terms of numerical accuracy in cases where steep gradients arise in the solution, while it is more robust than \texttt{ode45}, which in turn is generally less accurate and needs more points to converge, or even fails for stiff problems.

The computational times are generally larger (but comparable for all practical purposes) than those required by the \texttt{ode45} and \texttt{ode15s} solvers. However, our method results in smaller computational times than \texttt{ode45} and \texttt{ode15s} when the solution has to be computed in ``dense'' sets of points, because the method provides an analytical expression for the approximation of the solution upon training with a finite set of collocation points. Furthermore, in our opinion the larger computational times are also due to the fact that our home-made code is not optimized (its optimization is out of the scope of the current work).

Possible future work includes the implementation of our RPNN-based scheme for the solution of large-scale systems of stiff ODEs as they arise from biological and/or chemical kinetics problems, a comparison with other well-established ``classical'' methods for solving medium-to-large-scale stiff problems, a systematic investigation and selection of the hyper-parameters of the network ~\cite{dong2021modified,dong2021computing}, a comparison with other physics-informed machine learning approaches, such as Deep Learning (see, e.g., \cite{wu2018physics,raissi2018numerical,raissi2019physics,chen2021physics,lu2021deepxde}) and other ELM-based schemes (see, e.g., \cite{dong2021local,dong2021modified,dwivedi2020physics}), and finally a stability analysis and error propagation of the scheme.

\section*{Acknowledgments}
E.G. is supported by a 3-year scholarship from the Universit\`a degli Studi di Napoli Federico II, Italy, and
G.F. is supported by a 4-year scholarship from the Scuola Superiore Meridionale, Universit\`a degli Studi di Napoli Federico II, Italy.
This work is also partially supported by the Gruppo Nazionale per il Calcolo Scientifico - Istituto Nazionale di Alta Matematica (GNCS-INdAM), Italy, and by the Italian program ``Fondo Integrativo Speciale per la Ricerca (FISR)'' - B55F20002320001.

\section*{Statements and Declarations}
\subsection*{Authors' Contribution Statement}
C.S. conceived the idea, developed the methodology, wrote the first draft of the manuscript and supervised the project. E.G. and G.F. contributed equally to the work: they developed the code, helped shape research and analysis, performed the computational experiments
and wrote the first draft of the section on the numerical results. F.C. and D.d.S. developed, analysed and wrote the numerical aspects of this work, provided critical feedback, shaped research and revised the manuscript. F.C., D.d.S and C.S. designed the numerical experiments and analysed the results. All the authors discussed the results and contributed to the revision of the final manuscript.
\subsection*{Competing/Financial Interests}
The authors have no relevant financial or non-financial interests to disclose.
The authors have no competing interests to declare that are relevant to the content of this article.
\subsection*{Data Availability}
Not applicable.
\subsection*{Code availability}
The code is available to the reviewers upon request and will be uploaded to github if the manuscript is accepted.
\bibliographystyle{spmpsci}      
\bibliography{references}  

\end{document}